\documentclass{article}
\usepackage[german, english]{babel}
\usepackage{amsmath,amssymb,latexsym}

\newcommand{\infixand}{\text{ and }}
\newcommand{\nin}{\not\in}
\newcommand{\nobracket}{}
\newcommand{\nocomma}{}
\newcommand{\tmem}[1]{{\em #1\/}}
\newcommand{\tmmathbf}[1]{\ensuremath{\boldsymbol{#1}}}
\newcommand{\tmop}[1]{\ensuremath{\operatorname{#1}}}
\newcommand{\tmrsup}[1]{\textsuperscript{#1}}
\newcommand{\tmsamp}[1]{\textsf{#1}}
\newcommand{\tmstrong}[1]{\textbf{#1}}

\begin{document}

\

\begin{center}
  \ \ {\tmstrong{ON THE MARKOFF EQUATION }}
\end{center}

\begin{center}
  \ {\tmstrong{Norbert Riedel}}
\end{center}

\begin{tabular}{l}
  
\end{tabular}

\begin{tabular}{l}
  
\end{tabular}

{\tmstrong{{\tmem{Abstract}}}}{\tmem{ }}A triple $(a, b, c)$ of positive
integers is called a Markoff triple if it satisfies the Diophantine equation
\begin{eqnarray*}
  a^2 + b^2 + c^2 = a b c . &  & 
\end{eqnarray*}
Recasting the Markoff tree, whose vertices are Markoff triples, in the
framework of integral upper unitriangular 3x3 matrices, it will be shown that
the largest member of such a triple determines the other two uniquely. This
answers a question which has been open for 100 years. The solution of this
problem will be obtained in the course of a broader investigation of the
Markoff equation by means of 3x3 matrices.

\begin{tabular}{l}
  
\end{tabular}

\begin{tabular}{l}
  
\end{tabular}

{\tmem{{\tmstrong{Introduction}}}}

Markoff numbers, the solutions of the Markoff Diophantine equation, have
captured the imagination of mathematicians for over a century. Rooted in A.A.
Markoff's late 19th century work on binary quadratic forms and their
connection to the top hierarchy of the worst approximable (quadratic) numbers
by rationals, these numbers have found their place in seemingly unrelated
endeavors of mathematical activity, such as 4-dimensional manifolds ([HZ]),
quantum field theory ([CV]), hyperbolic geometry ([Se]), combinatorics ([Po]),
group and semi group theory ([Co],[Re]). Two in-depth treatments of the
classical aspects of the theory ([Ca], [CF]) bracket almost four decades. One
problem that has resisted a conclusive solution so far is the question whether
the largest number of a Markoff triple determines uniquely the other two. F.G.
Frobenius posed this question in 1913 ([F]). It was restated most recently by
M.Waldschmidt in ([W]). A brief discussion of the uniqueness question is
included in the exposition of Markoff's theory by E. Bombieri [Bo]. Over the
past twenty years various proofs were obtained showing the uniqueness of
dominant Markoff numbers which are powers of primes (again, see [Bo] for a
survey of the relevant literature). Most of these contributions, however, seem
to be superseded by a result which was published by B. Stolt in 1952 ([St],
Theorem 9; see also the discussion in Section 9). The primary objective in the
present work is to show that the answer is affirmative throughout, as
expressed by the following theorem, which is equivalent to a conjecture by A.
N. Tyurin in complex geometry, stating that a representative exceptional
bundle on the complex projective plane is uniquely determined by its rank. For
details see A. N. Rudakov's article [Ru].

$\begin{array}{l}
  
\end{array}$

{\tmstrong{Theorem}} \ Given two triples of positive integers, $(a_1
\nobracket$, $b_{1,} c_1$) and $(a_2 \nobracket$, $b_2$, $c_2$), such that

\ \ \ \ \ \ \ \ \ \ \ \ \ $a_k < b_k$ $<$ $c_k$ , \ \ and \ \ $a^2_k$ + $b_{
k}^2$ + $c^2_k$ \ = $a_k b_k c_k$ , \ \ $k \begin{array}{l}
  \epsilon
\end{array} \{ 1, 2 \}$ ,

it follows that \ $c_1$ = $c_2$ \ implies \ $a_1$ = $a_2$ \ and \ $b_1$ =
$b_2$.

\begin{tabular}{l}
  
\end{tabular}

However, since the techniques and formulae leading up to the proof of this
statement are far broader than the primary objective itself, a great deal of
effort will be dedicated to issues relating to, but not necessarily
indispensable for the proof. Hopefully, this broadened approach to the issues
involved will contribute to an enhanced understanding of the ideas and the
formalism which are so particular to the Markoff equation. We start by
encoding every Markoff triple in a (upper) triangular 3x3 matrix, with 1's in
the diagonal, and then move on to determine an explicit form for the
``isomorphs'' of these matrices. More specifically, given any pair of such
matrices, the connectedness of the Markoff tree gives rise to an integral
unimodular matrix transforming one into the other, in the same vein as
equivalent quadratic forms are related. An integral nilpotent rank 2 matrix,
which is associated (essentially uniquely) with each of the aforementioned
unitriangular matrices, gives rise to a one-variable parametrization of all
``automorphs'' of those triangular matrices. All of this will be covered in
Section 1 through Section 3. The parametrization of the ``automorphs''
obtained in Section 3 will lead in Section 4 to a diophantine matrix
equations, whose solutions are closely related to integers $n$ for which the
number $- 1$ is a quadratic residue modulo $n$. This in turn will lead to a
canonical matrix factorization of these solutions, which is particular to the
Markoff property. \ In Section 5 we will embark on a closer analysis of the
matrix which is at the center of the factorization obtained in Section 4. In
Section 6 we will draw some number theoretic conclusions which will \ lead to
further insight into the nature of cycles of reduced indefinite binary
quadratic forms containing Markoff forms. In particular we will show, that
such a cycle contains two symmetric forms, and furthermore, how Markoff
numbers can be characterized by means of this property. This central result
will be instrumental in the proof of the Theorem. In Section 7 we will prove
the Theorem. The proof will be carried out in two steps. In the first step we
will employ unique prime ideal factorization in quadratic number fields to
conclude that any two given Markoff forms associated with a common Markoff
number have to be (properly or improperly) equivalent. In the second step we
will employ the formalism developed in Section 6 to show how the uniqueness of
a pair of symmetric forms associated with a given Markoff number entails the
unicity claim of the Theorem. Working out the specific composition of the
forms involved in the first step of the proof in terms of Gauss' bilinear
substitutions (which have recently found their reincarnation in the so-called
Bhargava cubes), leads to some further insight into the connection between
Markoff numbers and certain principal forms. This topic will be taken up in
Section 8. Section 9 contains the brief discussion of a norm form equation
which depends on a given Markoff number and the affiliated discriminant only,
highlighting its connection with the uniqueness question. In Section 10 and
Section 11 there will be a discussion of recursions producing data affiliated,
and to some degree determined by a given Markoff number. Specifically, in
Section 10 we deal with canonical decompositions of the discriminant into sums
of two squares, while Section 11 engages in a discussion of the algebraic
framework in terms of 3x3 matrices for the quadratic residues.

Finally, here is a guide for the reader who wishes to focus exclusively on
the arguments providing a self-contained proof of the Theorem, while
dispensing with the buildup of the matrices instrumental in the proof. Since a
combination of Lemma 4.11 and Lemma 4.12 \ ensures the crucial divisibility
property employed in the proof of Proposition 5.1, it suffices to go through
these two lemmas before proceeding directly to Section 5.

\

\

\begin{tabular}{l}
  
\end{tabular}

\begin{tabular}{l}
  
\end{tabular}{\tmem{{\tmstrong{1 Markoff tree and triangular 3x3 matrices }}}}

\begin{tabular}{l}
  
\end{tabular}

Since the matrix manipulations employed in the first four sections render the
more common version of the Markoff equation
\[ \mathfrak{a}^2 +\mathfrak{b}^2 +\mathfrak{c}^2 =
   3\mathfrak{a}\mathfrak{b}\mathfrak{c}, \mathfrak{a}, \mathfrak{b},
   \mathfrak{c} \begin{array}{l}
     \epsilon
   \end{array} \mathbb{N} \]
impractical, we shall use in those sections mostly the alternative form
\[ a^2 + b^2 + c^2 = a b c, \]
where $a = 3\mathfrak{a}, b = 3\mathfrak{b}, c = 3\mathfrak{c}$. It is also
common to represent the three numbers as the components of a triple, arranged
in increasing order from the left to the right, for instance. This arrangement
is unsuitable for the objectives in the present section. While still referring
to this arrangement as a Markoff triple, and the largest number as the
dominant member, we will supplement this notion by the following, denoting by
$\tmmathbf{M}_n (\mathbb{Z})$ ( $\tmmathbf{M}_n^+$($\mathbb{Z}$) ) the set of
n$\times$n \ matrices whose entries are integers (non negative integers).

\begin{tabular}{l}
  
\end{tabular}

{\tmstrong{1.1 Definition}} \ A Markoff triple matrix, or MT-matrix, is a
matrix in $\tmmathbf{M}^+_3 (\mathbb{Z})$ of the form
\[ \left(\begin{array}{ccc}
     1 & a & b\\
     0 & 1 & c\\
     0 & 0 & 1
   \end{array}\right), \]
where \ $a^2$+$b^2$+$c^2$= $a b c$ , \ and max$\{ a, b, c \} \begin{array}{l}
  \epsilon
\end{array} \{ a, c \}$.

\begin{tabular}{l}
  
\end{tabular}

For each Markoff triple, with the exception of (3, 3, 3) and (3, 3, 6), there
are exactly four MT-matrices. We shall use the notation
\[ M (a, b, c) = \left(\begin{array}{ccc}
     1 & a & b\\
     0 & 1 & c\\
     0 & 0 & 1
   \end{array}\right) \]
for arbitrary entries $a, b, c$. Throughout this work, a matrix followed by an
upper right exponent $t$ denotes the corresponding transpose matrix.

\begin{tabular}{l}
  
\end{tabular}

{\tmstrong{1.2 Proposition}} \ For any two MT-matrices $M (a_1$, $b_1$, $c_1$)
and \ $M (a_2$, $b_2$, $c_2$) there exists

$N \begin{array}{l}
  \epsilon
\end{array} \tmop{SL} (3, \mathbb{Z})$ such that

a)
\[ N^t M (a_2, b_2, c_2) N = M (a_1, b_1, c_1), \]
b)
\[ N \left(\begin{array}{c}
     c_1\\
     - b_1\\
     a_1
   \end{array}\right) = \left(\begin{array}{c}
     c_2\\
     - b_2\\
     a_2
   \end{array}\right), N^t  \left(\begin{array}{c}
     c_2\\
     a_2 c_2 - b_2\\
     a_2
   \end{array}\right) = \left(\begin{array}{c}
     c_1\\
     a_1 c_1 - b_1\\
     a_1
   \end{array}\right) \]
\begin{tabular}{l}
  
\end{tabular}

{\tmstrong{Proof}} a) If
\[ P (x) = \left(\begin{array}{ccc}
     0 & - 1 & 0\\
     1 & x & 0\\
     0 & 0 & 1
   \end{array}\right), Q (y) = \left(\begin{array}{ccc}
     1 & 0 & 0\\
     0 & y & 1\\
     0 & - 1 & 0
   \end{array}\right), \]
then \ $P (x), Q (y) \begin{array}{l}
  \epsilon
\end{array} \tmop{SL} (3, \mathbb{Z}) $ for $x, y \begin{array}{l}
  \epsilon
\end{array} \mathbb{Z}$, and
\begin{eqnarray*}
  P (a)^t M (a, b, c) P (a) = M (a, c, a c - b) &  & 
\end{eqnarray*}
\[ Q (c)^t M (a, b, c) Q (c) = M (a c - b, a, c) . \]
If $M (a, b, c)$ is a MT-matrix, then the matrices on the right hand side are
also MT-matrices, and both are associated with the same neighbor of the
Markoff triple corresponding to the MT-matrix on the left hand side. Here the
word neighbor refers to two adjacent Markoff triples in the so-called Markoff
tree. By the very definition of MT-matrices the Markoff triple associated with
the right hand side is further removed from the root of the tree than the
corresponding triple on the left hand side. Furthermore, application of
transposition and conjugation by
\[ \mathcal{J}= \left(\begin{array}{ccc}
     0 & 0 & 1\\
     0 & 1 & 0\\
     1 & 0 & 0
   \end{array}\right) \]
to the two identities above leads to new identities:
\[ Q (a)^t M (c, b, a) Q (a) = M (a c - b, c, a), \]
\[ P (c)^t M (c, b, a) P (c) = M (c, a, a c - b) . \]
So, on the right hand side of these four identities combined, we obtain
exactly the four MT-matrices associated with a common Markoff triple. It
follows that, through repeated applications of the four identities, the
claimed statement is true in case $a_1$= $b_1$= $c_1$= 3. Notice that it is
vital that there is only one MT-matrix associated with the root of the Markoff
tree! The claim in the general case now follows immediately by combining the
special case applied to $M (a_1$, $b_1$, $c_1$) and to M($a_2$, $b_2$, $c_2$)
separately.

b) It suffices to note that the matrices $P (x)$ and$Q (y)$ have the claimed
property for the appropriate Markoff numbers $x$ and $y$, and so the claimed
identities follow from part a). \ $\Box$

\begin{tabular}{l}
  
\end{tabular}

{\tmstrong{Remarks}} 1) The first two of the identities in the proof of
Proposition 1.2 give rise to the definition of neighbors in a binary tree with
MT-matrices serving as vertices. The Markoff tree, which is not entirely
binary, can be recovered form this tree simply by identifying the four
MT-matrices with the Markoff triple they are associated with.

2) If
\[ N^t M (3, 3, 3) N = M (a, b, c), N^t  \left(\begin{array}{c}
     3\\
     6\\
     3
   \end{array}\right) = \left(\begin{array}{c}
     c\\
     a c - b\\
     a
   \end{array}\right), N \begin{array}{l}
     \epsilon
   \end{array} \tmop{SL} (3, \mathbb{Z}), \]
then
\[ N^{- 1} M (- 3, 6, - 3) (N^{- 1})^t = M (- a, a c - b, - c) . \]
Therefore, if
\[ \tilde{N}  = \left(\begin{array}{ccc}
     1 & 0 & 0\\
     0 & - 1 & 0\\
     0 & 0 & 1
   \end{array}\right) (N^{- 1})^t \left(\begin{array}{ccc}
     1 & 0 & 0\\
     0 & - 1 & 0\\
     0 & 0 & 1
   \end{array}\right) \]
then
\[ (\tilde{N})^t M (3, 6, 3)  \tilde{N} = M (a, a c - b, c), \tilde{N}
   \left(\begin{array}{c}
     c\\
     b - a c\\
     a
   \end{array}\right) = \left(\begin{array}{c}
     3\\
     - 6\\
     3
   \end{array}\right) . \]
Since
\[ P (3)^t M (3, 6, 3) P (3) = Q^t (3) M (3, 6, 3) Q (3) = M (3, 3, 3), \]
it follows that, given any two Markoff triples, any permutation of the first,
($a_1$, $b_1$, $c_1$) say, and any permutation of the second, ($a_2$, $b_2$,
$c_2$) say, there exists $N \begin{array}{l}
  \epsilon
\end{array}$SL(3,{\tmsamp{$\mathbb{Z})$}}, such that
\[ N^t M (a_2, b_2, c_2) N = M (a_1, b_1, c_1) . \]
3) Markoff triples have also been associated with triples of integral
unimodular matrices, exploiting two of the so-called Fricke identities. For an
in-depth survey of this approach, mostly due to H. Cohn, see [Pe]. The
connection between that approach and the present one is as follows: Let
\[ A_0 = \left(\begin{array}{cc}
     2 & 1\\
     1 & 1
   \end{array}\right) \tmop{and} \begin{array}{l}
     
   \end{array} B_0 = \left(\begin{array}{cc}
     1 & 1\\
     1 & 2
   \end{array}\right) . \]
We say that ($A_0$, $A_0 B_0$, $B_0$) is an admissible triple. New admissible
triples can be generated \ out of given ones by the rule, that if ($A, A  B$,
$B$) is an admissible triple, then so are$(A, A^2 B$, $\tmop{AB}$) and ($A B,
A B^2$, $B$). Fricke's identities ensure that the corresponding triple of
traces associated with an admissible triple solves the Markoff equation.
Moreover, the lower left entry of each matrix is one-third of its trace. So,
once again with the notion of neighbor defined in a natural way, the
admissible triples represent nothing but the vertices of the Markoff tree.
However, since ($\tmop{Tr} (A_0)$, Tr($A_0 B_0$), Tr($B_0$)) =(3, 6, 3), the
first Markoff triple (3, 3, 3) is missing from the picture. As pointed out in
the proof of Proposition 1.2, its availability in the present approach is
crucial, due to the fact that it is the only Markoff triple for which all
components are equal. Exploiting the fact that a matrix solves its own
characteristic equation, one can easily see that each matrix in an admissible
triple can be written as a linear combination of the matrices $A_0$, $A_0 B_0$
and $B_0$ with integral coefficients. If $a_2$=$b_2$=$c_2$=3 in Proposition
1.2, and if $N \tmop{is} \tmop{the} \tmop{matrix}$exhibited in its proof, then
the coefficient vectors for the admissible triple associated with $(c_1, a_1
c_1 - b_1, a_1)$ are exactly the columns of the matrix $N$ in the order of
their appearance. The 1's \ in the diagonal of the matrix $M (a_1$, $b_1$,
$c_1$) reflect the unimodularity of the 2$\times 2$ matrices in the
corresponding admissible triple. Other choices for the basis $A_0$, $A_0 B_0$
and $B_0$ appear in the literature, mostly motivated by the desire to connect
them to the continued fraction expansion of the quadratic irrationals, which
are at the core Markoff's original work. That all these choices are connected
via a single integral nilpotent 3$\times 3$ matrix, and that this matrix holds
the key to the uniqueness question of the Markoff triples, is one of the
central observations in the present work.

\begin{tabular}{l}
  
\end{tabular}

\begin{tabular}{l}
  
\end{tabular}

{\tmem{{\tmstrong{2. Markoff triples and nilpotent matrices}}}}

\begin{tabular}{l}
  
\end{tabular}

The statement of Proposition 1.2 raises the issue of ``automorphs'', to borrow
a notion from the theory of quadratic forms. More specifically, what can be
said about the matrices $N \begin{array}{l}
  \epsilon
\end{array} \tmop{SL} (3, \mathbb{Z})$ which leave $M$ invariant, i.e.
\[ N^t M (a, b, c) N = M (a, b, c) ? \]
There are two natural candidates that could serve as generators. While
defining them, we will temporarily relinquish the requirement that $a, b
\tmop{and} c$ are in $\mathbb{Z}$. A commutative ring will do. Let
\[ H (a, b, c) = M (a, b, c)^{- 1} M  (a, b, c)^t . \]
If possible, we will suppress the arguments.

\begin{tabular}{l}
  
\end{tabular}

{\tmstrong{2.1 Proposition}} a) $H^t M H$=$M$

b) If $N$ is invertible and $N^t M (a_2$, $b_2$, $c_2$)$N = M (a_1$, $b_1$,
$c_1$), then
\[ N^{- 1} H (a_2, b_2, c _2) N = H (a_1, b_1, c_1) . \]
\begin{tabular}{l}
  
\end{tabular}

{\tmstrong{Proof}} a)
\[ H^t M H = M (M^{- 1})^t M M^{- 1} M^t = M . \]
b) Writing
\[ M_k = M (a_k, b_k, c_k), H_k = M_k^{- 1} M_k^t, k \begin{array}{l}
     \epsilon
   \end{array} \{ 1, 2 \}, \]
$N^t M^t_2 N$=$M_1$ implies
\[ N^t M_2 N = M_1^t \tmop{and} N^{- 1} M_2^{- 1} (N^t)^{- 1} = M_{1 }^{- 1},
\]
so,
\[ N^{- 1} H_2 N = N^{- 1} M_2^{- 1} M_2^t N = N^{- 1} M_2^{- 1} (N^t)^{- 1}
   N^t M_2^t N = M_1^{- 1} M_1^t = H_1 \]
$\Box$

\begin{tabular}{l}
  
\end{tabular}

The explicit form of $H$ is
\[ H (a, b, c) = \left(\begin{array}{ccc}
     1 - (a^2 + b^2 - a b c)  & a c^2 - b c - a & a c - b\\
     a - b c & 1 - c^2  & - c\\
     b & c & 1
   \end{array}\right) \]
Its characteristic polynomial is given by \
\[ \det (H - \lambda E) = - (\lambda - 1)^3 - d (\lambda - 1)^2 - d (\lambda
   - 1), d = a^2 + b^2 + c^2 - a b c \]
\begin{tabular}{l}
  
\end{tabular}

{\tmstrong{Remark}} \ The matrix $H$has a place in quantum field theory
([CV]). More specifically $H$ (or rather its inverse), is the monodromy matrix
for the so-called $\tmop{CP}^2$ $\sigma$-model. This is a model with N=2
superconformal symmetry and Witten index n=3.

\begin{tabular}{l}
  
\end{tabular}

The other candidate is related to a matrix $R \epsilon \tmmathbf{M}_3
(\mathbb{Z})$ which solves the matrix equation

(2.1)
\[ R^t M + M R = 0 \]
This matrix is unique up to a multiplicative constant. We may choose
\[ R = \left(\begin{array}{ccccc}
     a^2 + b^2 - a b c & 2 a + b c - a c^2 &  &  & 2 b - a c\\
     b c - 2 a & c^2 - a^2  &  &  & 2 c - a b\\
     a c - 2 b & - 2 c - a b + a^2 c &  &  & a b c - b^2 - c^2
   \end{array}\right) \]
Its characteristic polynomial is
\[ \det (R - \lambda E) = - \lambda^3 + d (d - 4) \lambda, d = a^2 + b^2 + c^2
   - a b c \]
In the context of real numbers we can state the following:

\begin{tabular}{l}
  
\end{tabular}

{\tmstrong{2.2 Proposition }}a) For any $x \begin{array}{l}
  \epsilon
\end{array} \mathbb{R}$, \ ($e^{x R}$)$^t M e^{x R}$=$M$.

b) If (a ,b ,c) is a Markoff triple, then the adjugate matrix of $R$, i.e. the
transpose of the cofactor matrix, is
\[ R^{\tmop{adj}} = R^2 = 4 \left(\begin{array}{c}
     c\\
     - b\\
     a
   \end{array}\right) \left(\begin{array}{c}
     c, a c - b, a
   \end{array}\right) \]
\begin{tabular}{l}
  
\end{tabular}

{\tmstrong{Proof}} \ a) Since $(R^t)^k M$=(-1)\tmrsup{$k$} $M R^k$ for all $k
\begin{array}{l}
  \epsilon
\end{array} \mathbb{N}$,
\[ (e^{x R})^t M e^{x R} = \sum^{\infty}_{k, l = 0} \frac{1}{k!} \frac{1}{l!}
   x^{k + l} (R^t)^k M R^l = \sum^{\infty}_{k, l = 0} \frac{1}{k!}
   \frac{1}{l!} (- 1)^k x^{k + l} M R^{k + l} = M e^{- R} e^R = M . \]

b) This can of course be shown through straightforward calculations of the
nine minors of $R$, involving repeated applications of the Markoff property. A
more conceptional proof, however, is the following. First one observes that
\[ \mathfrak{J}^2 =\mathfrak{J}^{\tmop{adj}} \tmop{for} \mathfrak{J}=
   \left(\begin{array}{ccc}
     0 & 1 & 0\\
     0 & 0 & 1\\
     0 & 0 & 0
   \end{array}\right) . \]
Next, the operation of adjugation of a nonsingular matrix commutes obviously
with any similarity transformation. Perturbing a singular matrix into a
nonsingular one, and then letting that perturbation approach the original
matrix shows that adjugation and similarity transformations commute for
singular matrices as well. This, once again, settles the claim b).

$\Box$

\begin{tabular}{l}
  
\end{tabular}

{\tmstrong{Remark \ }}In reference \,to Remark 3 in Section1, the conjugation
of $N$ by $e^{- \frac{x}{6} R}$ corresponds to the conjugation of the
components of the related admissible triple by the matrix
\[ \left(\begin{array}{cc}
     1 & x\\
     0 & 1
   \end{array}\right) . \]

The matrices $H$ and $R$ commute, and so they share common eigenvectors. Let
us briefly consider $R$ in the context of the ring $P_{\mathbb{Z}}$[X], the
polynomials with integral coefficients. There are exactly two cases in which
$R$ is nilpotent, namely $d$=0 and \ $d$=4. The case $d$=0 leads us to Markoff
triples, while the case $d$=4 leads us to triples of Tchebycheff polynomials:
For the root of the tree we choose the triple (X, X, 2), X being the free
variable. Beginning at the root, we obtain three adjacent (but not necessarily
distinct) triples out of a given one, $(P_1, P_2, P_3)$ say, as follows.
\[ (P_2 P_3 - P_1, P_2, P_3), (P_1, P_1 P_3 - P_2, P_3), (P_1, P_2, P_1 P_{2
   .} - P_3) \]
The polynomials thus obtained are monic polynomials which are mutually
orthogonal with respect to a certain probability measure derived from
classical potential theory. The triples of integers representing the degrees
of these polynomials form the vertices of the so-called ``Euclid tree''. While
the kinship between the cases $d$=0 and $d$=4 goes well beyond the shared
nilpotence of $R$, a fact which has been exploited by Zagier in [Z] with
profit in deriving an asymptotic bound for Markoff numbers through comparison
of the two cases, the uniqueness question, which is the subject of the present
investigation, has clearly a negative answer in the case $d$=4. The crucial
difference between these two cases is the fact that, while $R$ is of rank 2 in
the case $d$=0, it is of rank 1 in the case $d$=4. Notice also that, while $H
- E$ is nilpotent for $d$=0, it still has two equal but non-vanishing
eigenvalues for $d$=4.

From now on we will be exclusively concerned with Markoff triples.
Let{\tmstrong{ }}
\[ S = H - E, \]
where $E$ denotes the unit matrix.

\begin{tabular}{l}
  
\end{tabular}

{\tmstrong{2.3 Proposition}} a) $H = e^{- \frac{R}{2}}$=$E - \frac{1}{2} R +
\frac{1}{8} R^2, R = 3 E - 4 H + H^2$

b)
\[ S^2 = \left(\begin{array}{c}
     c\\
     - b\\
     a
   \end{array}\right) \left(\begin{array}{c}
     c, a c - b, a
   \end{array}\right) \]

The proof is obtained through straightforward manipulations, involving
repeated employment of the Markoff property. Proposition 2.3 shows that we are
essentially dealing with a single nilpotent matrix of rank 2. It will follow
from our subsequent discussion that all ``automorphs'' have the form \ $e^{s
R}$ for a suitable rational parameter $s$. Since the matrix $R$ has some mild
redundancies, thus making manipulations a bit more lengthy, and since these
redundancies are not shared by the matrix $S$, we will be working in the
sequel with $S$ only.

Before we are going to embark on the parametrization of all ``automorphs'' of
the matrices $\begin{array}{l}
  M (a, b, c)
\end{array}$via the Jordan normal form, yielding rational matrices which are
crucially non-integral, we digress briefly to present a normal form for $R$
which highlights the integrality of the matrix $R$.

\begin{tabular}{l}
  
\end{tabular}

{\tmstrong{2.4 Proposition}} a) For each Markoff triple $(a, b, c)$ there
exists a matrix $\mathcal{W} (a, b, c) \begin{array}{l}
  \epsilon
\end{array} \tmop{GL} (3, \mathbb{Z})$ such that

(2.2)
\[ \mathcal{W} (a, b, c)^{- 1} \frac{1}{3} R\mathcal{W} (a, b, c) =
   \left(\begin{array}{ccc}
     0 & 2 & 1\\
     0 & 0 & 2\\
     0 & 0 & 0
   \end{array}\right) \]
(2.3)
\[ \mathcal{W} (a, b, c) \left( \begin{array}{l}
     1\\
     0\\
     0
   \end{array} \right) = \left(\begin{array}{c}
     \mathfrak{a}\\
     -\mathfrak{b}\\
     \mathfrak{c}
   \end{array}\right), \det (\mathcal{W}(a, b, c)) = - 1 \]
b) If $N$ is the matrix constructed in Proposition 1.2 for two Markoff triples

$(a_1, b_1, c_1) \tmop{and} (a_2, b_2, c_2)$ then

(2.4)
\[ N =\mathcal{W} (a_2, b_2, c_2) \mathcal{W} (a_1, b_1, c_1)^{- 1} \]
\begin{tabular}{l}
  
\end{tabular}

{\tmstrong{Proof}} a) If $(a, b, c)$=(3,3,3), then
\[ \mathcal{W} (3, 3, 3) = \left(\begin{array}{ccc}
     1 & - 2 & 0\\
     - 1 & 1 & 1\\
     1 & 0 & - 1
   \end{array}\right) \]
has the claimed properties. Now, letting $N$ be the matrix constructed in
Proposition 1.2 such that
\[ N^t M (a , b , c ) N = M (3, 3, 3), \]
then, by Proposition 2.1 b) and Proposition 2.3 a)
\[ N^{- 1} R (a, b, c) N = R (3, 3, 3) . \]
We define
\[ \mathcal{W} (a, b, c) = N\mathcal{W} (3, 3, 3) . \]
By Proposition 1.2 b), and since $N \epsilon \tmop{SL} (3, \mathbb{Z})$, the
matrix $\mathcal{W} (a, b, c)$ has the claimed properties.

b) This is an immediate consequence of the construction of the matrix $N$ in
Proposition 1.2 a). $\Box$

\begin{tabular}{l}
  
\end{tabular}

{\tmstrong{Remark }}The brevity of the proof of Proposition 2.4 obscures the
significance of what's going on here. Some background information might
elucidate the issues, especially when this normal form is being compared to
the way in which the use of the Jordan normal form unfolds in Section 3.
First, the normal form enunciated in Proposition 2.4 is essentially unique. In
order to clarify this statement one has to place the search for such a normal
form on a more systematic footing. Specifically, the proper context for doing
so is the Smith normal form for integral matrices. (See for instance [Ne],
Chapter II, for an in depth exposition of this subject). In the case of a 3x3
integral matrix there are exactly 3 determinantal divisors, ${d_1} , d_2,
d_3$: \ $d_1$ is the greatest common divisor of all nine matrix entries, $d_2$
is the greatest common divisor of all nine entries in the corresponding
adjugate matrix, and $d_3$ is the is the determinant of the given matrix. For
a nonsingular integral 3x3 matrix the Smith normal form is then given by
diag$(d_1, \frac{d_2}{d_1}, \frac{d_3}{d_2})$, and for a singular integral 3x3
matrix of rank 2 it is given by diag$(d_1, \frac{d_2}{d_1}, 0)$ \ The diagonal
entries in the Smith normal form are called the invariant factors of the
matrix. The Smith normal form is known to be invariant under left as well as
right multiplication by unimodular integral matrices. In our case one can see
that the matrix $\mathcal{R}= \frac{1}{3} R$ has the Smith normal form
diag$(1, 4, 0)$. Indeed, the first determinantal divisor of $\mathcal{R}$ is
equal to 1 (for instance the greatest common divisor of the entries (1,1),
(1,2) and (3,3) is equal to 1), by Proposition 2.2(b) the second determinantal
divisor of $\mathcal{R}$ is equal to 4, and finally, det$(\mathcal{R}) = 0$.
By [Ne], Theorem III.12 there exists a matrix $\mathcal{W} \epsilon \tmop{GL}
(3, \mathbb{Z})$, such that
\[ \mathcal{W}^{- 1} \mathcal{R}\mathcal{W}= \left(\begin{array}{ccc}
     0 & \alpha & \beta\\
     0 & 0 & \gamma\\
     0 & 0 & 0
   \end{array}\right) \begin{array}{l}
     \epsilon
   \end{array} \tmmathbf{M}_3 (\mathbb{Z}) . \]
By [AO], section 11, we may assume that
\[ \alpha > 0, \gamma > 0, \tmop{and} 0 \leqslant \beta < \gcd (\alpha,
   \gamma), \]
rendering this choice unique in the sense that any two similar matrices of
this form must be identical (see [AO], section 12).The only non-vanishing
cofactor in the upper triangular matrix is the one with index (1,3). Since the
second determinantal divisor of $\mathcal{R}$ is equal to 4, we conclude that
$\alpha \beta = 4$. This in turn entails that
\[ \mathcal{W}^{- 1} \mathcal{R}\mathcal{W} \epsilon \{
   \left(\begin{array}{ccc}
     0 & 2 & 1\\
     0 & 0 & 2\\
     0 & 0 & 0
   \end{array}\right), \left(\begin{array}{ccc}
     0 & 1 & 0\\
     0 & 0 & 4\\
     0 & 0 & 0
   \end{array}\right), \left(\begin{array}{ccc}
     0 & 4 & 0\\
     0 & 0 & 1\\
     0 & 0 & 0
   \end{array}\right) \} . \]
In order to show that only the first matrix can occur, we note that, on the
one hand both, $\mathcal{R}$

and $\left(\begin{array}{ccc}
  0 & 2 & 1\\
  0 & 0 & 2\\
  0 & 0 & 0
\end{array}\right)$ are similar (with respect to $\tmop{GL} (3, \mathbb{Z})$)
to the negative of their transpose matrix.

In the first case this follows from (2.1), and in the latter case we have
\[ \left(\begin{array}{ccc}
     0 & 0 & - 1\\
     0 & 1 & 0\\
     - 1 & 1 & 0
   \end{array}\right) \left(\begin{array}{ccc}
     0 & 2 & 1\\
     0 & 0 & 2\\
     0 & 0 & 0
   \end{array}\right) \left(\begin{array}{ccc}
     0 & 1 & - 1\\
     0 & 1 & 0\\
     - 1 & 0 & 0
   \end{array}\right) = - \left(\begin{array}{ccc}
     0 & 0 & 0\\
     2 & 0 & 0\\
     1 & 2 & 0
   \end{array}\right) . \]
On the other hand, conjugating the matrix $\left(\begin{array}{ccc}
  0 & 1 & 0\\
  0 & 0 & 4\\
  0 & 0 & 0
\end{array}\right)$ to the negative of its transpose,
\[ X \left(\begin{array}{ccc}
     0 & 1 & 0\\
     0 & 0 & 4\\
     0 & 0 & 0
   \end{array}\right) X^{- 1} = - \left(\begin{array}{ccc}
     0 & 0 & 0\\
     1 & 0 & 0\\
     0 & 4 & 0
   \end{array}\right), \]
leads to a matrix of the from
\[ X = \left(\begin{array}{ccc}
     0 & 0 & - 4 x\\
     0 & x & y\\
     - 4 x & - y & z
   \end{array}\right) ; x, y, z \begin{array}{l}
     \epsilon
   \end{array} \mathbb{Z}. \]
Since the determinant of this matrix is $- 16 x^3$, $X$ can never be chosen to
be unimodular. Since
\[ \left(\begin{array}{ccc}
     0 & 4 & 0\\
     0 & 0 & 1\\
     0 & 0 & 0
   \end{array}\right) =\mathcal{J} \left( \begin{array}{c}
     \begin{array}{l}
       0\\
       0\\
       0
     \end{array}
   \end{array} \begin{array}{l}
     1\\
     0\\
     0
   \end{array} \begin{array}{l}
     0\\
     4\\
     0
   \end{array} \right)^{\tmmathbf{t}} \mathcal{J}, \]
the matrix $\left(\begin{array}{ccc}
  0 & 4 & 0\\
  0 & 0 & 1\\
  0 & 0 & 0
\end{array}\right)$ too is not similar to the negative of its transpose via an
integral unimod-ular matrix. Finally, since the property of a matrix to be
similar to the negative of its transpose is invariant under similarity
transformations, there has to exist a matrix $\mathcal{W} \begin{array}{l}
  \epsilon
\end{array} \tmop{GL} (3, \mathbb{Z})$ such that
\[ \mathcal{W}^{- 1} \mathcal{R}\mathcal{W}= \left(\begin{array}{ccc}
     0 & 2 & 1\\
     0 & 0 & 2\\
     0 & 0 & 0
   \end{array}\right), \]
as claimed in Proposition 2.4. Up to this point (2.2) and (2.3) follow, except
for the sign of the determinant of the matrix $\mathcal{W}$. What is not
obtainable through this line of reasoning, however, is (2.4). The best one can
get is
\[ N =\mathcal{W} (a_2, b_2, c_2) \mathcal{W} (a_1, b_1, c_1)^{- 1} +
   \varepsilon \left(\begin{array}{c}
     c_2\\
     - b_2\\
     a_2
   \end{array}\right)  \left(\begin{array}{c}
     c_1, a_1 c_1 - b_1, a_1
   \end{array}\right), \]
with an unspecified integer $\varepsilon$. Notice that the value of
$\varepsilon$ does not change if any of the matrices $\mathcal{W}$ is
mutltiplied from the left by a matrix with determinant 1 which commutes with
the corresponding nilpotent matrix $R$.

\begin{tabular}{l}
  \\
  
\end{tabular}

{\tmem{{\tmstrong{3 Determination of automorphs}}}}

\begin{tabular}{l}
  
\end{tabular}

There are two objectives in this section. First we seek to develop a
one-parameter characterization of the ``automorphs'' introduced in the
previous section. Second, this should be done in a way so that the formulas
that will be the point of departure in the the next section emerge in a
natural way from the ones we obtain in the present one. In order to meet those
two requirements, we shall, for the most part in this section, work with two
Markoff triples which share a common member. First we construct a matrix $T$
which conjugates $S$ to its Jordan form. Starting with an eigenvector for
$S^t$ yields
\[ S \left(\begin{array}{c}
     c\\
     a c - b\\
     a
   \end{array}\right) = a c \left(\begin{array}{c}
     a c^2 - b c - a\\
     - c^2\\
     c
   \end{array}\right) \]
Notice that the vector on the right hand side is nothing but the second column
of $S$ multiplied by $a c$. Applying $S \tmop{to} \tmop{its}$second column
yields by virtue of the Markoff property
\[ (a c - b) \left(\begin{array}{c}
     c\\
     - b\\
     a
   \end{array}\right), \]
which is in the kernel of $S$. So, if we define
\[ T = \left(\begin{array}{ccc}
     c & a c (a c^2 - b c - a) & a c (a c - b) c\\
     a c - b & a c (- c^2) & a c (a c - b) (- b)\\
     a & a c c & a c (a c - b) a
   \end{array}\right) \]
then we have
\[ S T = T \left(\begin{array}{ccc}
     0 & 0 & 0\\
     1 & 0 & 0\\
     0 & 1 & 0
   \end{array}\right) . \]

Furthermore,
\[ \det (T) = - [a c (a c - b)]^3 \]
In order to manage the manipulations involving this matrix efficiently, we
will use a suitable factorization. If
\[ A = \left(\begin{array}{cccccc}
     0 &  & c (a c - b) - a &  &  & c\\
     1 &  & - c^2 &  &  & - b\\
     0 &  & c &  &  & a
   \end{array}\right), \]
\[ B = \left(\begin{array}{ccc}
     a c & 0 & 0\\
     0 & 1 & 0\\
     0 & 0 & 1
   \end{array}\right) \]
\[ C = \left(\begin{array}{ccc}
     1 & 0 & 0\\
     0 & 1 & 0\\
     1 & 0 & 1
   \end{array}\right), \]

\[ D = \left(\begin{array}{cccccc}
     1 &  &  & 0 &  & 0\\
     0 &  &  & a c &  & 0\\
     0 &  &  & 0 &  & a c (a c - b)
   \end{array}\right), \]

then $T = A B C D$. Moreover
\[ A^{- 1} = - \frac{1}{(a c - b)^2 } \left(\begin{array}{ccc}
     - c (a c - b)  & - (a c - b)^2  & - a (a c - b)\\
     - a & 0 & c\\
     c & 0 & a - c (a c - b)
   \end{array}\right) \]

\[ = \frac{1}{(a c - b)^2} F K L, \]

where
\[ F = \left(\begin{array}{ccccc}
     a c - b &  & 0 &  & 0\\
     0 &  & 1 &  & 0\\
     0 &  & 0 &  & 1
   \end{array}\right) \]
\[ K = \left(\begin{array}{cccccc}
     c &  &  & 1 &  & a\\
     a &  &  & 0 &  & - c\\
     - c &  &  & 0 &  & c (a c - b) - a
   \end{array}\right) \]
\[ L = \left(\begin{array}{cccc}
     1 &  & 0 & 0\\
     0 &  & a c - b & 0\\
     0 &  & 0 & 1
   \end{array}\right) \]
We shall also need the matrix
\[ U = M T = V B C D, \]
where
\[ V = \left(\begin{array}{ccc}
     a & - a & c\\
     1 & 0 & m\\
     0 & c & a
   \end{array}\right) \]
\[ V^{- 1} = \frac{1}{(a c - b)^2} \left(\begin{array}{ccc}
     c (a c - b) & - b (a c - b)  & a (a c - b)\\
     a & - a^2  & a (a c - b) - c\\
     - c & a c & - a
   \end{array}\right) \]

Now consider two Markoff triples (at this point not necessarily distinct) with
a common \ member $m$. We assume that
\[ m = a_1 c_1 - b_1 = a_2 c_2 - b_2, \]
where $a_k$, $b_k$ and $c_k$ are the components of the unique neighbor closer
to the root of the Markoff tree, for $k = 1$ and $k = 2$, respectively. This
arrangement accommodates all vertices of the Markoff tree except for the root.
In order to make use of the matrices introduced above in the present context,
we adopt the convention of attaching an index 1 or 2 to their names, depending
on the Markoff triple in reference. Let
\[ \hat{N} = T_2 T_1^{- 1}, r = \frac{a_1 c_1}{a_2 c_2} . \]
Then
\[ \det (r \hat{N}) = 1 \]

By Proposition 1.2 there exists a matrix $N \begin{array}{l}
  \epsilon
\end{array} \tmop{SL} (3, \mathbb{Z}) \tmop{such} \tmop{that}$

(3.1)
\[ N^t M (a_2, b_2, c_2) N = M (a_1, b_1, c_1) . \]
By Proposition 2.1(b)
\[ N^{- 1} S_2 N = S_1 \]
Since $\hat{N}^{- 1} S_2 \hat{N}$=$S_1$, it follows that $N \hat{N}^{- 1}$ and
$S_2$ commute. Since $S_2$ has rank 2, this implies that there exist rational
numbers $s$ and $t$, such that

(3.2)
\[ N  = r (E + s  S_2 + t S_2^2) \hat{N} = r T_2 \left(\begin{array}{ccc}
     1 & 0 & 0\\
     s & 1 & 0\\
     t & s & 1
   \end{array}\right) T_1^{- 1} \begin{array}{l}
     \epsilon
   \end{array} \tmmathbf{M}_3 (\mathbb{Z}) . \]
Substituting (3.2) into (3.1) yields the identity

(3.3)
\[ r (T_1^t)^{- 1} \left(\begin{array}{ccc}
     1 & s & t\\
     0 & 1 & s\\
     0 & 0 & 1
   \end{array}\right) T_2^t = r^{- 1} (M (a_2, b_2, c_2) T_2 
   \left(\begin{array}{ccc}
     1 & 0 & 0\\
     s & 1 & 0\\
     t & s & 1
   \end{array}\right) (M (a_1, b_1, c_1) T_1)^{- 1})^{- 1} = \]
\[ r^{- 1} U_1 \left(\begin{array}{cccccc}
     1 &  & 0 &  &  & 0\\
     - s &  & 1 &  &  & 0\\
     s^2 - t &  & - s &  &  & 1
   \end{array}\right) U_2^{- 1} \]
We are now going to evaluate the three terms in (3.2). Writing $F, L$ in place
of $F_1$, $L_1$, respectively,
\[ r m^2 \hat{N} = r A_2 B_2 C_2 D_2 D_1^{- 1} C_1^{- 1} B_1^{- 1} F K_1 L =
\]
\[ A_2 \left(\begin{array}{ccccc}
     1 &  & 0 &  & 0\\
     0 &  & 1 &  & 0\\
     \frac{1}{a_2 c_2} - \frac{1}{a_1 c_1} &  & 0 &  & 1
   \end{array}\right) F K_1 L = \]

\[ \left(\begin{array}{cccc}
     c_2 (\frac{1}{a_2 c_2} - \frac{1}{a_1 c_1}) m & c_2 m - a_2 & c_2 & \\
     (1 - b_2 (\frac{1}{a_2 c_2} - \frac{1}{{a_1}  c_1})) m & - c_2^2  & - b_2
     & \\
     a_2 (\frac{1}{a_2 c_2} - \frac{1}{a_1 c_1}) m & c_2  & a_2 & 
   \end{array}\right) K_1 L = \]

\[ \Gamma_0 + m \Gamma_1 + m^2 \Gamma_2, \]
where,
\[ \Gamma_0 = \left(\begin{array}{ccccc}
     - (a_1 a_2 + c_1 c_2)  &  & 0 &  & - (a_1 c_2 - c_1 a_{ 2})\\
     - (a_1 c_2 - c_1 a_{ 2}) c_2  &  & 0 &  &  (a_1 c_2 - c_1 a_{ 2}) a_2\\
     a_1 c_2 - c_1 a_{ 2}  &  & 0 &  & - (a_1 a_2 + c_1 c_2)
   \end{array}\right) \]
\[ + m \left(\begin{array}{ccc}
     a_1 c_2 & 0 & 0\\
     0 & 0 & 0\\
     0 & 0 & c_1 a_2
   \end{array}\right) \]
\[ \Gamma_1 = (\frac{1}{a_2 c_2} - \frac{1}{a_1 c_1})  \left(\begin{array}{c}
     c_2\\
     - b_2\\
     a_2
   \end{array}\right)  \left(\begin{array}{c}
     c_1, m, a_1
   \end{array}\right) \]
\[ \Gamma_2 = \left(\begin{array}{ccc}
     0 & 0 & 0\\
     0 & 1 & 0\\
     0 & 0 & 0
   \end{array}\right) \]
Since

\[ m^2  \left(\begin{array}{ccc}
     0 & 0 & 0\\
     1 & 0 & 0\\
     0 & 1 & 0
   \end{array}\right) T^{- 1} = a^{- 1} c^{- 1}  \left(\begin{array}{ccc}
     0 & 0 & 0\\
     m & 0 & 0\\
     0 & 1 & 0
   \end{array}\right) K L = a^{- 1} c^{- 1} L \left(\begin{array}{ccccc}
     0 &  & 0 & 0 & \\
     c &  & 1 & a & \\
     a &  & 0 & - c & 
   \end{array}\right) L, \]
we get for the second term
\[ m r S_2 \hat{N} = m r T_2 \left(\begin{array}{ccc}
     0 & 0 & 0\\
     1 & 0 & 0\\
     0 & 1 & 0
   \end{array}\right) T_1^{- 1} = A_2 \left(\begin{array}{ccccc}
     0 &  & 0 & 0 & \\
     c_1 &  & 1 & a_1 & \\
     a_1 &  & 0 & - c_1 & 
   \end{array}\right) L = \]
\[ \Omega_0 + m \Omega_1, \]
where \
\[ \Omega_0 = \left(\begin{array}{ccccc}
     a_1 c_2 - c_1 a_2 &  & 0 &  & - (a_1 a_2 + c_1 c_2)\\
     - (a_1 a_2 + c_1 c_2) c_2 &  & 0 &  & - (a_1 c_2 - c_1 a_2) c_2\\
     a_1 a_2 + c_1 c_2 &  & 0 &  & a_1 c_2 - c_1 a_2
   \end{array}\right), \]
\[ \Omega_1 = \left(\begin{array}{ccc}
     0 & - a_2  & 0\\
     a_1 & - c_2^2  & - c_1\\
     0 & c_2 & 0
   \end{array}\right) + c_2 \left(\begin{array}{ccc}
     c_1 & m & a_1\\
     0 & 0 & 0\\
     0 & 0 & 0
   \end{array}\right) \]
Finally, for the third term
\[ r S_2^2 \hat{N} = \Phi^t = \left(\begin{array}{c}
     c_2\\
     - b_2\\
     a_2
   \end{array}\right)  \left(\begin{array}{c}
     c_1, m, a_1
   \end{array}\right) . \]

In order to manipulate the identity (3.3) we shall need a similar
decomposition involving the matrix $U$.
\[ r^{- 1} m^2 U_1 U_2^{- 1} = V_1 \left(\begin{array}{ccccccc}
     1 &  &  & 0 &  &  & 0\\
     0 &  &  & 1 &  &  & 0\\
     - (\frac{1}{a_2 c_2} - \frac{1}{a_1 c_1}) &  &  & 0 &  &  & 1
   \end{array}\right) \left(\begin{array}{ccc}
     c_2 m & - b_2 m & a_2 m\\
     a_2 & - a_2^2  & a_2 m - c_2\\
     - c_2 & a_2 c_2  & - a_2 
   \end{array}\right) \]
\[ = \Theta_0 + m \Theta_1 + m^2 \Theta_2, \]
where
\[ \Theta_0 = \left(\begin{array}{ccc}
     - (a_1 a_2 + c_1 c_2)  & - (a_1 c_2 - c_1 a_{ 2}) c_2  & a_1 c_2 - c_1
     a_{ 2}\\
     0 & 0 & 0\\
     - (a_1 c_2 - c_1 a_{ 2}) &  (a_1 c_2 - c_1 a_{ 2}) a_2  & - (a_1 a_2 +
     c_1 c_2)
   \end{array}\right) \]
\[ + m \left(\begin{array}{cccc}
     a_1 c_2 &  & 0 & 0\\
     0 &  & 0 & 0\\
     0 &  & 0 & c_1 a_2
   \end{array}\right), \]
\[ \Theta_1 = - (\frac{1}{a_2 c_2} - \frac{1}{a_1 c_1}) 
   \left(\begin{array}{c}
     c_1\\
     m\\
     a_1
   \end{array}\right)  \left(\begin{array}{c}
     c_2, - b_2, a_2
   \end{array}\right), \]
\[ \Theta_2 = \Gamma_2 = \left(\begin{array}{ccc}
     0 & 0 & 0\\
     0 & 1 & 0\\
     0 & 0 & 0
   \end{array}\right) . \]
Since
\[ m^2 \left(\begin{array}{ccc}
     0 & 0 & 0\\
     1 & 0 & 0\\
     0 & 1 & 0
   \end{array}\right) U^{- 1} = a^{- 1} c^{- 1} \left(\begin{array}{ccccc}
     0 &  & 0 &  & 0\\
     c m &  & - b m &  & a  m\\
     a &  & - a^2 &  & a m - c
   \end{array}\right), \]
we get
\[ r^{- 1} m U_1 \left(\begin{array}{ccc}
     0 & 0 & 0\\
     1 & 0 & 0\\
     0 & 1 & 0
   \end{array}\right) U_2^{- 1} = V_1 \left(\begin{array}{ccc}
     0 & 0 & 0\\
     {c_2}  & - b_2  & a_2\\
     a_2 & - a_2^2  & a_2 m - c_2 
   \end{array}\right) \]
\[ = \Lambda_0 + m \Lambda_1, \]
where
\[ \Lambda_0 = \left(\begin{array}{ccc}
     - (a_1 c_2 - c_1 a_{ 2})  &  (a_1 c_2 - c_1 a_{ 2}) a_2 & - (a_1 a_2 +
     c_1 c_2) \\
     0 & 0 & 0\\
     a_1 a_2 + c_1 c_2  & - (a_1 a_2 + c_1 c_2) a_2  & - (a_1 c_2 - c_1 a_{
     2}) 
   \end{array}\right), \]
\[ \Lambda_1 = \left(\begin{array}{ccc}
     0 & - a_1 & 0\\
     a_2 & - a^2_2  & - c_2\\
     0 & c_1 & 0
   \end{array}\right) + a_2 \left(\begin{array}{ccc}
     0 & 0 & c_1\\
     0 & 0 & m\\
     0 & 0 & a_1
   \end{array}\right) . \]
Finally,
\[ r^{- 1} U_1 \left(\begin{array}{ccc}
     0 & 0 & 0\\
     0 & 0 & 0\\
     1 & 0 & 0
   \end{array}\right) U_2^{- 1} = \Phi \]
Let

(3.4)
\[ N (s) = r e^{- \frac{R_2}{2} s} \hat{N} - \frac{1}{m} (\frac{1}{a_2 c_2} -
   \frac{1}{a_1 c_1}) \left(\begin{array}{c}
     c_2\\
     - b_2\\
     a_2
   \end{array}\right) \left(\begin{array}{c}
     c_1, m, a_1
   \end{array}\right) = \]
\[ r \hat{N} e^{- \frac{R_1}{2} s} - \frac{1}{m} (\frac{1}{a_2 c_2} -
   \frac{1}{a_1 c_1}) \left(\begin{array}{c}
     c_2\\
     - b_2\\
     a_2
   \end{array}\right) \left(\begin{array}{c}
     c_1, m, a_1
   \end{array}\right) \]
Then we have the following crucial representation of all ``rational
isomorphs''.

\begin{tabular}{l}
  
\end{tabular}

{\tmstrong{3.1 Proposition}} If $Q \begin{array}{l}
  \epsilon
\end{array} \tmop{SL} (3, \mathbb{Q})$, then

(3.5)
\[ Q^t M_2 Q = M_1, \]
if and only if there exists a rational number $s$ such that $Q = N (s)$.

\begin{tabular}{l}
  
\end{tabular}

{\tmstrong{Proof}} First, by our discussion above, we know that if (3.5) holds
true, then there exist rational numbers $s \tmop{and} t$, such that
\[ Q = r (E + s  S_2 + t S_2^2) \hat{N} . \]
Now given this representation, $Q$ satisfies (3.5) if and only if
\[ (3.6) r (T_1^t)^{- 1} \left(\begin{array}{ccc}
     1 & s & t\\
     0 & 1 & s\\
     0 & 0 & 1
   \end{array}\right) T_2^t - r^{- 1} U_1 \left(\begin{array}{cccccc}
     1 &  & 0 &  &  & 0\\
     - s &  & 1 &  &  & 0\\
     s^2 - t &  & - s &  &  & 1
   \end{array}\right) U_2^{- 1} = 0 . \]
Employing the above decompositions, the left hand side of (3.6) turns into
\[ \frac{1}{m^2} \Gamma_0^t + \frac{1}{m} \Gamma_1^t + \Gamma_2^t +
   \frac{s}{m} \Omega_0^t + s \Omega_1^t + t \Phi  - \frac{1}{m^2} \Theta_0 -
   \frac{1}{m} \Theta_1 - \Theta_2 + \frac{s}{m} \Lambda_0 + s \Lambda_1 -
   (s^2 - t) \Phi  . \]
Since
\[ \Gamma_0^t = \Theta_0, \Gamma_1^t = - \Theta_1 = (\frac{1}{a_2 c_2} -
   \frac{1}{a_1 c_1}) \Phi, \Gamma_2^t = \Theta_2, \]
the left hand side of (3.6) simplifies to
\[ \frac{s}{m} (\Omega_0^t + \Lambda_0) + s (\Omega_1^t + \Lambda_1) + (2
   (\frac{1}{m} (\frac{1}{a_2 c_2} - \frac{1}{a_1 c_1}) + t) - s^2) \Phi . \]
But
\[ \Omega_1^t + \Lambda_1 = \Phi + \left(\begin{array}{ccc}
     0 & c_1 b_2 & 0\\
     0 & 0 & 0\\
     0 & a_1 b_2 & 0
   \end{array}\right), \]
while
\[ \Omega_0^t + \Lambda_0 = - m \left(\begin{array}{ccc}
     0 & c_1 b_2 & 0\\
     0 & 0 & 0\\
     0 & a_1 b_2 & 0
   \end{array}\right), \]
so that the left hand side of (3.6) finally takes the form
\[ (2 (\frac{1}{m} (\frac{1}{a_2 c_2} - \frac{1}{a_1 c_1}) + t) + s - s^2)
   \Phi . \]
This expression is equal to zero if and only if
\[ t = \frac{1}{2} (s^2 - s) - \frac{1}{m} (\frac{1}{a_2 c_2} - \frac{1}{a_1
   c_1}), \]
which is equivalent to $Q = N (s)$. \ $\Box$

\begin{tabular}{l}
  
\end{tabular}

{\tmstrong{Remarks}} 1) If $a_1 = a_2$, $c_1 = c_2$, then the proof of
Proposition 3.1 shows that all ``automorphs'' of an MT-matrix are of the form
$e^{\frac{R }{6} s}$ for some integer $s$.

2) All integral ``isomorphs'' are actually contained in a proper congruence
subgroup of $\tmop{SL} (3, \mathbb{Z})$, namely the matrices which are
orthogonal modulo 3.

3) Notice that due to cancellation the matrix $N (s)$ can be written more
compactly as follows,
\[ N (s) = \frac{1}{m^2} \Gamma_0 + \Gamma_2 + \frac{s}{m} (\Omega_0 + m
   \Omega_1) + \frac{s^2 - s}{2} \Phi^t . \]
\begin{tabular}{l}
  
\end{tabular}

\begin{tabular}{l}
  
\end{tabular}

{\tmem{{\tmstrong{4 Markoff triples and quadratic residues}}}}

\begin{tabular}{l}
  
\end{tabular}

The point of departure in this section is the following matrix identity within
the settings of section 3. Let
\[ W = W (a_i, b_i, c_i) = \left(\begin{array}{ccc}
     c_i & 0 & 2 a_i - m c_i\\
     - b_i & 1 & c_i^2 - a_i^2\\
     a_i & 0 & m a_i - 2 c_i
   \end{array}\right), i = 1, 2 \]
then

(4.1)
\[ N (0) = W (a_2, b_2, c_2) W (a_1, b_1, c_1)^{- 1} = \frac{1}{2 m^2} W (a_2,
   b_2, c_2) W (a_1, b_1, c_1)^{\tmop{adj}} \]
\[ = \frac{1}{2 m^2} \left(\begin{array}{ccc}
     c_2 & 0 & 2 a_2 - m c_2\\
     - b_2 & 1 & c_2^2 - a_2^2\\
     a_2 & 0 & m a_2 - 2 c_2
   \end{array}\right) \left(\begin{array}{ccc}
     m a_1 - 2 c_1 & 0 & m c_1 - 2 a_1\\
     2 m c_1 & 2 m^2 & 2 m a_1\\
     - a_1 & 0 & c_1
   \end{array}\right) . \]
This identity separates the two Markoff triples with the property $m = a_1 c_1
- b_1 = a_2 c_2 - b_2 $. For a single Markoff triple $(a, b, c)$ the matrix $W
(a, b, c)$ in turn gives rise to the following identity

(4.2)
\[ W (a, b, c)^t M (a, b, c) W (a, b, c) = \left(\begin{array}{ccc}
     0 & m & 0\\
     m & 1 & m^2\\
     0 & - m^2 & - 4 m^2
   \end{array}\right) . \]

Significantly, the matrix on the right hand side depends on $m$ only. Also
notice that an application of the matrix \ $R_{}$ to the second column of $W 
$ yields the third column, while an application of $R_{}$to the third column
of $W  $ yields $4 m$ times the first column of $W  $. The following identity
exhibits the intrinsic symmetry of this matrix,
\[ W (c, b, a) =\mathcal{J}W (a, b, c) \left(\begin{array}{ccc}
     1 & 0 & 0\\
     0 & 1 & 0\\
     0 & 0 & - 1
   \end{array}\right) . \]
To get a better understanding of the architecture of the matrix on the right
hand side of (4.2) we observe first that
\[ \left(\begin{array}{ccc}
     0 & m & 0\\
     m & 1 & m^2\\
     0 & - m^2 & - 4 m^2
   \end{array}\right) = \left(\begin{array}{ccc}
     0 & m & 0\\
     m & 1 & m^2\\
     0 & m^2 & 4 m^2
   \end{array}\right)  \left(\begin{array}{ccc}
     1 & 0 & 0\\
     0 & 1 & 0\\
     0 & 0 & - 1
   \end{array}\right) . \]
Both of these factors are associated with the nilpotent matrix on the right
hand side of the conjugation
\[ W (a, b, c)^{- 1} R (a, b, c) W (a, b, c) = \left(\begin{array}{ccc}
     0 & 0 & 4 m\\
     0 & 0 & 0\\
     0 & 1 & 0
   \end{array}\right) \]
as follows: The first factor, which is \ self-adjoint, conjugates the matrix
$\left(\begin{array}{ccc}
  0 & 0 & 4 m\\
  0 & 0 & 0\\
  0 & 1 & 0
\end{array}\right)$ to its adjoint $\left(\begin{array}{ccc}
  0 & 0 & 0\\
  0 & 0 & 1\\
  4 m & 0 & 0
\end{array}\right)$, while the second factor, which is a self-adjoint
involution, conjugates $\left(\begin{array}{ccc}
  0 & 0 & 4 m\\
  0 & 0 & 0\\
  0 & 1 & 0
\end{array}\right)$ to $- \left(\begin{array}{ccc}
  0 & 0 & 4 m\\
  0 & 0 & 0\\
  0 & 1 & 0
\end{array}\right)$. Any matrix $Y$ of this design has the following
``automorph'' property,
\[ \exp \left( s \left(\begin{array}{c}
     \begin{array}{l}
       0
     \end{array} \begin{array}{l}
       0
     \end{array} \begin{array}{l}
       0
     \end{array}\\
     \begin{array}{l}
       0
     \end{array} \begin{array}{l}
       0
     \end{array} \begin{array}{l}
       1
     \end{array}\\
     4 m \begin{array}{l}
       0
     \end{array} \begin{array}{l}
       0
     \end{array}
   \end{array}\right) \right) Y \exp \left( s \left(\begin{array}{c}
     \begin{array}{l}
       0
     \end{array} \begin{array}{l}
       0
     \end{array} \begin{array}{l}
       4 m
     \end{array}\\
     \begin{array}{l}
       0
     \end{array} \begin{array}{l}
       0 \begin{array}{l}
         
       \end{array}
     \end{array} \begin{array}{l}
       0
     \end{array}\\
     \begin{array}{l}
       0
     \end{array} \begin{array}{l}
       1 \begin{array}{l}
         
       \end{array}
     \end{array} \begin{array}{l}
       0
     \end{array}
   \end{array}\right) \right) = Y, \]
where exp(.) denotes the exponential function, and s is a rational number.
This construction works essentially for any nilpotent 3x3 matrix. Conversely,
any matrix $Y$ with the indicated ``automorph'' property must have the form
\[ Y = \left(\begin{array}{ccc}
     0 & \alpha & 0\\
     \alpha & \gamma & \beta\\
     0 & - \beta & - 2 m \alpha
   \end{array}\right), \]
for arbitrary values $\alpha, \beta \tmop{and} \gamma$. In the context of
(1.2) we have $\alpha = m, \beta = m^2 \tmop{and} \gamma = 1$.

\begin{tabular}{l}
  
\end{tabular}

{\tmstrong{Remark}} The following observation, which will not be used in the
sequel, is of some interest. Let
\[ \mathcal{U}= W \left(\begin{array}{ccc}
     1 & 0 & 0\\
     - m & 1 & 0\\
     0 & 0 & 1
   \end{array}\right) = \left(\begin{array}{ccc}
     c & 0 & 2 a - m c\\
     - a c & 1 & c^2 - a^2\\
     a & 0 & m a - 2 c
   \end{array}\right), \]
and let
\[ \mathcal{Q}= \left(\begin{array}{ccc}
     1 \begin{array}{l}
       
     \end{array} & \frac{a}{2} \begin{array}{l}
       
     \end{array} & \frac{b}{2}\\
     \frac{a}{2} \begin{array}{l}
       
     \end{array} & 1 \begin{array}{l}
       
     \end{array} & \frac{c}{2}\\
     \frac{b}{2} \begin{array}{l}
       
     \end{array} & \frac{c}{2} \begin{array}{l}
       
     \end{array} & 1
   \end{array}\right) = \frac{1}{2} (M (a, b, c) + M (a, b, c)^t) \]
Then
\[ \mathcal{U}^t \mathcal{Q}\mathcal{U}= \left(\begin{array}{ccc}
     - m^2 & 0 & 0\\
     0 & 1 & 0\\
     0 & 0 & - 4 m^2
   \end{array}\right) . \]
This identity shows that the column vectors of the matrix $\mathcal{U}$ form
an orthogonal basis with respect to the (indefinite) ternary quadratic form
determined by the symmetric matrix $\mathcal{Q}$, which has a determinant
equal to 1 if and only if $(a, b, c)$ is a Markoff triple.

\begin{tabular}{l}
  
\end{tabular}

We are now going to state a property that exhibits the intrinsic rigidity of
the identity (4.2).

\begin{tabular}{l}
  
\end{tabular}

{\tmstrong{4.1 Proposition}} For any four positive integers $a, b, c, q$, the
following two conditions are equivalent

a) The triple $(a, b, c)$ is Markoff, and $q = a c - b$.

b) There exists an integral 3x3 matrix $W = (w_{\tmop{ij}})$ with the
properties
\[ W^t M (a, b, c) W = \left(\begin{array}{ccc}
     0 & q & 0\\
     q & 1 & q^2\\
     0 & - q^2 & - 4 q^2
   \end{array}\right), \]
and $w_{12} = w_{32} = 0, w_{22} = 1$.

\begin{tabular}{l}
  
\end{tabular}

{\tmstrong{Proof}} By (4.2), a) implies b). To show that b) implies a), we
first observe that

$\det \left(\begin{array}{ccc}
  0 & m & 0\\
  m & 1 & m^2\\
  0 & - m^2 & - 4 m^2
\end{array}\right) = 4 m^4$, and therefore $\det (W) = \pm 2 m^2$. Replacing
$W$ by $- W$ if necessary, and letting $X = (x_{\tmop{ij}}) = W^{\tmop{adj}}$,
we can restate the matrix identity in b) as follows,
\[ \left(\begin{array}{ccc}
     x_{11} & x_{21} & x_{31}\\
     0 & 2 m^2 & 0\\
     x_{13} & x_{23} & x_{33}
   \end{array}\right) \left(\begin{array}{ccc}
     0 & m & 0\\
     m & 1 & m^2\\
     0 & - m^2 & - 4 m^2
   \end{array}\right) \left(\begin{array}{ccc}
     x_{11} & 0 & x_{13}\\
     x_{21} & 2 m^2 & x_{23}\\
     x_{31} & 0 & x_{33}
   \end{array}\right) = 4 m^4 \left(\begin{array}{ccc}
     1 & a & b\\
     0 & 1 & c\\
     0 & 0 & 1
   \end{array}\right) . \]
Reading off the identities for those entries only which are located on or
below the diagonal, with the exception of entry (2,2) which is trivial,
yields,

Entry (1,1): \ $2 m x_{11} x_{21} + x_{21}^2 - 4 m^2 x_{31}^2 = 4 m^4$

Entry (2,1): \ $m x_{11} + x_{21} + m^2 x_{31} = 0$

Entry (3,1): \ $m x_{11} x_{23} + m x_{13} x_{21} + x_{21} x_{23} - m^2 x_{21}
x_{33} + m^2 x_{23} x_{31} - 4 m^2 x_{31} x_{33} = 0$

Entry (3,2): \ $m x_{13} + x_{23} - m^2 x_{33} = 0$

Entry (3,3): \ $2 m x_{13} x_{23} + x_{23}^2 - 4 m^2 x_{33}^2 = 4 m^4$.

Combining the identities from entries (1,1) and (2,1) yields

(4.3)
\begin{eqnarray*}
  x^2_{21} + 2 m^2 x_{21} x_{31} + 4 m^2 x^2_{31} + 4 m^4 = 0 & . & 
\end{eqnarray*}
Combining the identities from entries (3,2) and (3,3) yields

(4.4)
\[ x_{23}^2 - 2 m^2 x_{23} x_{33} + 4 m^2 x_{33}^2 + 4 m^4 = 0 . \]
Finally, substituting the identities from entries (2,1)and (3,2) into the
identity for entry (3,1) yields

(4.5)
\[ x_{21} x_{23} + 4 m^2 x_{31} x_{33} = 0 . \]
It follows from (4.3) and (4.4), respectively, that $x_{21}$ and $x_{23}$ are
divisible by $2 m$. Thus, letting
\[ x^{\ast}_{21} = \frac{x_{21}}{2 m}, x^{\ast}_{23} = \frac{x_{23}}{2 m}, \]
we obtain

$(4.3)^{\ast}$
\[ (x^{\ast}_{21})^2 + m x^{\ast}_{21} x_{31} + x^2_{31} + m^2 = 0, \]
$(4.4)^{\ast}$
\[ (x^{\ast}_{23})^2 - m x^{\ast}_{23} x_{33} + x^2_{33} + m^2 = 0, \]
$(4.5)^{\ast}$
\[ x^{\ast}_{21} x^{\ast}_{23} + x_{31} x_{33} = 0 . \]
It follows from $(4.3)^{\ast}$ through $(4.5)^{\ast}$ that
\[ | x^{\ast}_{21} | = | x _{33} |, | x^{\ast}_{23} | = | x _{31} |,
   x^{\ast}_{21} x _{31} < 0 . \]
These last four conditions show that, up to a minus sign, the matrix $X$ has
exactly the form of the adjugate matrix of $W$ in (4.1). It is now
straightforward to check that the numbers $a, b, c, m$ have the properties
claimed in a). $\Box$

\begin{tabular}{l}
  
\end{tabular}

{\tmstrong{Remarks}} 1) More generally, the following can be shown.

There exists a matrix $X = (x^{(1)}, x^{(2)}, x^{(3)}) \begin{array}{l}
  \epsilon
\end{array} M_3 (\mathbb{Z})$ solving the matrix equation
\[ X^t M (a, b, c) X = \left(\begin{array}{ccc}
     0 & q & 0\\
     q & 1 & q^2\\
     0 & - q^2 & - 4 q^2
   \end{array}\right), \]
such that the vector $x^{(2)}$ has length 1, if and only if $(a, b, c)$ is a
Markoff triple. Moreover
\[ q = c \tmop{if} x^{(2)} = \left(\begin{array}{c}
     1\\
     0\\
     0
   \end{array}\right), \]
\[ q = a c - b \tmop{if} x^{(2)} = \left(\begin{array}{c}
     0\\
     1\\
     0
   \end{array}\right), \]
\[ q = a \tmop{if} x^{(2)} = \left(\begin{array}{c}
     0\\
     0\\
     1
   \end{array}\right) \]

2) It is quite natural to wonder to what degree the matrix on the right hand
side of (4.2) is uniquely determined by the discussion so far. The answer is,
not as much as one is led to suspect. As a matter of fact, essentially
everything that has been said so far works with minor adjustments just as well
if the matrix $W$=$W (a, b, c)$ is being replaced by
\[ Z^{\tmop{adj}} = \left(\begin{array}{ccc}
     - 2 c & 2 a - m c & m c\\
     2 a c & c^2 - a^2 & - m b\\
     - 2 a & m a - 2 c & m a
   \end{array}\right), \]
where
\[ Z = Z (a, b, c)  = \left(\begin{array}{ccc}
     c & m & a\\
     - a & 0 & c\\
     a & 2 & c
   \end{array}\right), \]
satisfying \ $\det (Z) = 2 m^2$. All one has to do is to replace the second
condition in Proposition 4.1 part b) by the following,
\[ z_{12} = m, z_{22} = 0, z_{32} = 2 . \]
The identity taking the place of (4.2) then becomes

(4.6)
\[ (Z^{\tmop{adj}})^t M (a, b, c) Z^{\tmop{adj}} = \left(\begin{array}{ccc}
     - 4 m^2 & 2 m^3 & 2 m^3\\
     - 2 m^3 & - 4 m^2 & 0\\
     2 m^3 & 0 & 0
   \end{array}\right) = 2 m^2 \left(\begin{array}{ccc}
     - 2 & m  & m \\
     - m  & - 2 & 0\\
     m  & 0 & 0
   \end{array}\right) . \]
One also has the identity,
\[ Z W = \left(\begin{array}{ccc}
     0 & m & 0\\
     0 & 0 & 2 m^2\\
     2 m & 2 & 0
   \end{array}\right) . \]
In a way the matrix $W$ is related to the nilpotent matrix $R$, while the
matrix $Z$ is related to $R^t$. In fact, $R^t$ applied to the last column
vector of $Z^t$ yields four times the second column vector of $Z^t$, while an
application of $R^t$ to the second column vector of $Z^t$ yields $2 m$ times
the first row vector of $Z^t$. \ As we shall see shortly, however, the
relationship between these two matrices is more intimate than appears to be
the case at first sight.

\begin{tabular}{l}
  
\end{tabular}

The next step is to consider the general diophantine matrix equation

(4.7)
\[ X^t M (a, b, c) X = \left(\begin{array}{ccc}
     0 & q & 0\\
     q & 1 & q^2\\
     0 & - q^2 & - 4 q^2
   \end{array}\right) ; \begin{array}{l}
     
   \end{array} X \begin{array}{l}
     \epsilon
   \end{array} \tmmathbf{M} _3 \text{(} \mathbb{Z}), \begin{array}{l}
     
   \end{array} q \begin{array}{l}
     \epsilon
   \end{array} \mathbb{N}, \]
where $(a, b, c)$ is a Markoff triple. We shall see shortly that, after
imposing a slightly technical restriction, this equation has a solution if and
only if $q$ is divisible by 3 and $- 1$ is a quadratic residue modulo
$\frac{q}{3}$. Before we go into that we give a brief summary of some
pertinent number theoretic facts, which can be readily gleaned from the
standard literature (see [L], for part a) and b); [Pn], Satz 2.4, or more
generally, [M1], [Ni] for part c)).

\begin{tabular}{l}
  
\end{tabular}

{\tmstrong{4.2 Lemma}} a) For any integer $n$ there exists an element
$\varepsilon$ in the residue class ring $\mathbb{Z}_n$, such that
$\varepsilon^2 = - 1$ (in other words, $- 1$ is a quadratic residue modulo
$n$) if and only if any odd prime factor $p$ of $n$ has the property $p = 1$
modulo 4.

b) If $n$ is an integer which is not divisible by 4 such that $- 1$ is a
quadratic residue modulo $n$, and $l$ is the number of distinct odd prime
factors dividing $n$, then there are exactly $2^l$ elements in $\mathbb{Z}_n$
whose square is equal to $- 1$.

c) For any solution of the diophantine equation $k^2 + 1 = n l ;
\begin{array}{l}
  
\end{array} n, l > 0$, there exists a matrix

$\left(\begin{array}{c}
  \begin{array}{l}
    p\\
    r
  \end{array} \begin{array}{l}
    q\\
    s
  \end{array}
\end{array}\right) \begin{array}{l}
  \epsilon
\end{array} \tmop{Sl} (2, \mathbb{Z})$ with such that
\[ \left(\begin{array}{c}
     \begin{array}{l}
       p\\
       r
     \end{array} \begin{array}{l}
       q\\
       s
     \end{array}
   \end{array}\right) \left(\begin{array}{c}
     \begin{array}{l}
       p\\
       r
     \end{array} \begin{array}{l}
       q\\
       s
     \end{array}
   \end{array}\right) \begin{array}{l}
     t\\
     
   \end{array} = \left(\begin{array}{c}
     \begin{array}{l}
       p\\
       r
     \end{array} \begin{array}{l}
       q\\
       s
     \end{array}
   \end{array}\right) \left(\begin{array}{c}
     \begin{array}{l}
       p\\
       q
     \end{array} \begin{array}{l}
       r\\
       s
     \end{array}
   \end{array}\right) = \left(\begin{array}{c}
     \begin{array}{l}
       n\\
       k
     \end{array} \begin{array}{l}
       k\\
       l
     \end{array}
   \end{array}\right) . \]
Given such a matrix $\left(\begin{array}{c}
  \begin{array}{l}
    p\\
    r
  \end{array} \begin{array}{l}
    q\\
    s
  \end{array}
\end{array}\right)$, all the other matrices in $\tmop{Sl} (2, \mathbb{Z})$
having this property are $\left(\begin{array}{c}
  \begin{array}{l}
    - p\\
    - r
  \end{array} \begin{array}{l}
    - q\\
    - s
  \end{array}
\end{array}\right)$, $\left(\begin{array}{c}
  \begin{array}{l}
    - q\\
    - s
  \end{array} \begin{array}{l}
    p\\
    r
  \end{array}
\end{array}\right)$, and $\left(\begin{array}{c}
  \begin{array}{l}
    q\\
    s
  \end{array} \begin{array}{l}
    - p\\
    - r
  \end{array}
\end{array}\right)$.

\begin{tabular}{l}
  
\end{tabular}

If $(a, b, c)$ is a Markoff triple, we let $\mathfrak{m}= \frac{a c - b}{3}$
if $a c - b$ is odd, and $\mathfrak{m}= \frac{a c - b}{6}$ if $a c - b$ is
even. Since $a^2 + c^2 \equiv 0$ (mod $\mathfrak{m}$), we have $\alpha^2 = -
1$, where $\alpha$ is the element in $\mathbb{Z}_{\mathfrak{m}} $
corresponding to $\frac{a}{c}$.

We are now going to tackle the system (4.7). First we take the transpose
matrices on both sides of (4.7) and subtract the result from (4.7), yielding

(4.8)
\[ X^t \left(\begin{array}{ccc}
     0 & a & b\\
     - a & 0 & c\\
     - b & - c & 0
   \end{array}\right) X = \left(\begin{array}{ccc}
     0 & 0 & 0\\
     0 & 0 & 2 q^2\\
     0 & - 2 q^2 & 0
   \end{array}\right) . \]
The matrix in the middle on the left hand side has rank 2, and the vector $(c,
- b, a)^t$ is in the kernel of this matrix. Writing
\[ X = (x^{(1)}, x^{(2)}, x^{(3)}), \]
we infer from (4.8) that
\[ (x^{(1)})^t \left(\begin{array}{ccc}
     0 & a & b\\
     - a & 0 & c\\
     - b & - c & 0
   \end{array}\right) x^{(i)} = 0, \begin{array}{l}
     
   \end{array} \tmop{for} i = 1, 2, 3. \]
Since $X$ is invertible, its column vectors are linearly independent, and this
entails
\[ (x^{(1)})^t \left(\begin{array}{ccc}
     0 & a & b\\
     - a & 0 & c\\
     - b & - c & 0
   \end{array}\right) = 0, \]
which in turn implies that $x^{(1)}$ and $(c, - b, a)^t$ are linearly
dependent. We now impose the technical restriction mentioned above.

(4.9)
\[ (x^{(1)})^t = (c, - b, a)^t . \]
Under this assumption, we are first going to deal with the Markoff triple
(3,3,3). In other words, we want to solve the system

(4.10)
\[ \left(\begin{array}{ccc}
     x_{11} & x_{21} & x_{31}\\
     x_{12} & x_{22} & x_{32}\\
     x_{13} & x_{23} & x_{33}
   \end{array}\right) \left(\begin{array}{ccc}
     1 & 3 & 3\\
     0 & 1 & 3\\
     0 & 0 & 1
   \end{array}\right) \left(\begin{array}{ccc}
     x_{11} & x_{12} & x_{13}\\
     x_{21} & x_{22} & x_{23}\\
     x_{31} & x_{32} & x_{33}
   \end{array}\right) = \]
\[ \left(\begin{array}{ccc}
     x^2_{11} + 3 x_{11} x_{21} + x_{21}^2 & x_{11} x_{12} + 3 x_{11} x_{22} +
     x_{21} x_{22} \begin{array}{l}
       
     \end{array} & x_{11} x_{13} + 3 x_{11} x_{23} + x_{21} x_{23}\\
     + 3 x_{11} x_{31} + 3 x_{21} x_{31} + x^2_{31} \begin{array}{l}
       
     \end{array} & + 3 x_{11} x_{32} + 3 x_{21} x_{32} + x_{31} x_{32}
     \begin{array}{l}
       
     \end{array} \begin{array}{l}
       
     \end{array} & + 3 x_{11} x_{33} + 3 x_{21} x_{33} + x_{31} x_{33}\\
     &  & \\
     x_{11} x_{12} + 3 x_{12} x_{21} + x_{21} x_{22} \begin{array}{l}
       
     \end{array} & x^2_{12} + 3 x_{12} x_{22} + x^2_{22} & x_{12} x_{13} + 3
     x_{12} x_{23} + x_{22} x_{23}\\
     + 3 x_{12} x_{31} + 3 x_{22} x_{31} + x_{31} x_{32} \begin{array}{l}
       
     \end{array} & + 3 x_{12} x_{32} + 3 x_{22} x_{32} + x^2_{32} & + 3 x_{12}
     x_{33} + 3 x_{22} x_{33} + x_{32} x_{33}\\
     &  & \\
     x_{11} x_{13} + 3 x_{13} x_{21} + x_{21} x_{23} \begin{array}{l}
       
     \end{array} \begin{array}{l}
       
     \end{array} & x_{12} x_{13} + 3 x_{13} x_{22} + x_{22} x_{23} & x^2_{13}
     + 3 x_{13} x_{23} + x^2_{23}\\
     + 3 x_{13} x_{31} + 3 x_{23} x_{31} + x_{31} x_{33} \begin{array}{l}
       
     \end{array} & + 3 x_{13} x_{32} + 3 x_{23} x_{32} + x_{32} x_{33} & + 3
     x_{13} x_{33} + 3 x_{23} x_{33} + x^2_{33}
   \end{array}\right) \]
\[ = \left(\begin{array}{ccc}
     0 & q & 0\\
     q & 1 & q^2\\
     0 & - q^2 & - 4 q^2
   \end{array}\right), \]
where at this point we shall assume that $q$ is an odd integer divisible by 3.
By assumption we have,

(4.11)
\[ x_{11} = x_{31} = - x_{21} = 3 . \]
The entry (1,2) or (2,1) yields,

(4.12)
\[ x_{12} + 2 x_{22} + x_{32} = \frac{q}{3} . \]
The entry (1.3) or (3.1) yields,

(4.13)
\[ x_{13} + 2 x_{23} + x_{33} = 0 . \]
While entry (3,3) yields only $x_{23} + x_{33} = \pm 2 q$, subtracting entry
(3,2) from entry (2,3) yields,

(4.14)
\[ \text{ } x_{23} + x_{33} = 2 q . \]
Combination of \ (4.12) with entry (2,2) yields,

(4.15)
\[ (\frac{q}{3})^2 + \frac{q}{3} (x_{32} - x_{22}) - (x_{32} + x_{22})^2 = 1
\]
Finally, combining (4.12), (4.14) and entry (2,3) yields,

(4.16)
\[ x_{33} - 6 (x_{32} + x_{22}) = q. \]
Up to this point we have only extracted necessary conditions for the
solvability of (4.9) and (4.10). We turn now to their sufficiency.

Letting

(4.17)
\[ \tmmathbf{\alpha}= x_{23}, \tmmathbf{\beta}= x_{33}, \tmmathbf{\gamma}=
   x_{32} - x_{22}, \varepsilon = x_{32} + x_{22}, \]
we can recast the above identities as follows. Instead of (4.14) we write

(4.18)
\[ \tmmathbf{\alpha}+\tmmathbf{\beta}= 2 q . \]
Instead of (4.15) we write

(4.19)
\[ (\frac{q}{3}) (\frac{q}{3} +\tmmathbf{\gamma}) - \varepsilon^2 = 1 . \]
Instead of (4.16) we write

(4.20)
\[ \tmmathbf{\beta}= 6 \varepsilon + q . \]
The identity (4.19) is telling us that
\[ \varepsilon^2 = - 1 (\tmop{mod} \frac{q}{3}) . \]
So if we let $\varepsilon$ be any number with this property,
\[ \varepsilon^2 = - 1 + \frac{q}{3} j, \]
and if we let
\[ \tmmathbf{\gamma}= j - \frac{q}{3}, \tmmathbf{\beta}= 6 \varepsilon + q,
   \tmmathbf{\alpha}= 2 q -\tmmathbf{\beta}, \]
then we can use (4.12) and the last two identities in (4.17) to solve for all
three entries in the second column of the matrix $X$. Note that, since by our
assumption $q$ is an odd integer, (4.19) shows that $\varepsilon$ is odd
(even) if and only if $\tmmathbf{\gamma}$ is odd (even). This ensures that by
virtue of the last two identities in (4.17) the numbers $x_{22}$ and $x_{32}$
are integers. Hence, $x_{12}$ is an integer as well. By (4.13), (4.18) and
(4.20) we can now solve for the three entries in the last column of the matrix
$X$ in terms of $\varepsilon \tmop{and} j$ as well, all numbers being
integers. To summarize, we have shown that all integral solutions of the
system (4.9) and (4.10) have the form

(4.21)
\[ X = \left(\begin{array}{ccc}
     3 \begin{array}{l}
       
     \end{array} & \frac{q}{6} - \frac{3 \varepsilon}{2} + \frac{j}{2} & - 3 q
     + 6 \varepsilon\\
     - 3 \begin{array}{l}
       
     \end{array} & \frac{1}{2} (\varepsilon - j + \frac{q}{3}) & q - 6
     \varepsilon\\
     3 \begin{array}{l}
       
     \end{array} & \frac{1}{2} (\varepsilon + j + \frac{q}{3}) & q + 6
     \varepsilon
   \end{array}\right), \]
provided $q$ is an odd integer which is divisible by 3, and the integers
$\varepsilon \tmop{and} j$ solve the diophantine equation $\epsilon^2 = - 1 +
\frac{q}{3} j$. In order to see that the same conclusion holds for even
integers $q$ as well, we observe that if $q$ is even, then $\varepsilon$ has
to be odd. But this means that $\varepsilon^2 + 1 = 4 n + 2$ for some integer
$n$, and therefore $j$ has to be odd as well. So it follows that all three
entries in the second column of the matrix $X$ are integers in case $q$ is
even. In conclusion, what we have shown is the first part of the following
statement.

\begin{tabular}{l}
  
\end{tabular}

{\tmstrong{4.3 Proposition}} a) The system (4.9) and (4.10) has an integral
solution $X$ if and only if $- 1$ is a quadratic residue modulo the integer
$\frac{q}{3}$.

b) The integral solutions of (4.9) and (4.10) are completely parametrized by
all the square roots of $- 1$ in the residue class ring associated with
$\frac{q}{3}$ in the following sense: Given two pairs of numbers,
$(\varepsilon_1, j_1)$ and $(\varepsilon_2, j_2)$ satisfying (4.21), such that
$\varepsilon_1 \equiv \varepsilon_2  (\tmop{mod} \frac{q}{3})$, each giving
rise to the solutions $X_1 \tmop{and} X_2$ of (4.9) and (4.10), respectively,
there exists an integer $i$ such that
\[ e^{\frac{i}{2} \mathcal{R}} X_1 = X_2, \]
where
\[ \mathcal{R}= \left(\begin{array}{ccc}
     - 3 & - 4 & - 1\\
     1 & 0 & - 1\\
     1 & 4 & 3
   \end{array}\right) . \]

\

The proof of the second part of this proposition will be given below where we
deal with general Markoff triples.

\begin{tabular}{l}
  
\end{tabular}

{\tmstrong{Remark}} It follows from (4.16), (4.14) and (4.13) that all entries
in the last column of $X$ are divisible by 3.

\begin{tabular}{l}
  
\end{tabular}

The following corollary, which is of some independent interest, will not be
used in the sequel.

\begin{tabular}{l}
  
\end{tabular}

{\tmstrong{4.4 Corollary}} If $X$ is an integral solution of (4.9) and (4.10),
and
\[ X^{\star} = \left(\begin{array}{ccc}
     0 & 0 & 1\\
     0 & 1 & 0\\
     1 & 0 & 0
   \end{array}\right) X \left(\begin{array}{ccc}
     1 & 0 & 0\\
     0 & 1 & 0\\
     0 & 0 & - 1
   \end{array}\right), \]
then
\[ e^{(E - \frac{\tmmathbf{\beta}}{2 q} \mathcal{R})} X = X^{\star} . \]
\begin{tabular}{l}
  
\end{tabular}

{\tmstrong{4.5 Corollary}} a) For any Markoff triple the system (4.7) and
(4.9) has an integral solution $X$ if and only if $- 1$ is a quadratic residue
modulo the integer $\frac{q}{3}$.

b) The integral solutions of (4.7) and (4.9) are completely parametrized by
all the square roots of $- 1$ in the residue class ring associated with
$\frac{q}{3}$ in the following sense: Given two pairs of numbers,
$(\varepsilon_1, j_1)$ and $(\varepsilon_2, j_2)$ satisfying (4.21), such that
$\varepsilon_1 \equiv \varepsilon_2  (\tmop{mod} \frac{q}{3})$, each giving
rise to the solutions $X_1 \tmop{and} X_2$ of (4.7) and (4.9), respectively,
there exists an integer $i$ such that
\[ e^{\frac{i}{2} \mathcal{R}} X_1 = X_2, \]
where
\[ \mathcal{R}= \frac{1}{3} \left(\begin{array}{ccccc}
     a^2 + b^2 - a b c & 2 a + b c - a c^2 &  &  & 2 b - a c\\
     b c - 2 a & c^2 - a^2  &  &  & 2 c - a b\\
     a c - 2 b & - 2 c - a b + a^2 c &  &  & a b c - b^2 - c^2
   \end{array}\right) . \]
\begin{tabular}{l}
  
\end{tabular}

{\tmstrong{Proof}} By Proposition 1.2 a) any integral solution of the system
(4.7) for the Markoff triple (3,3,3) can be transformed into an integral
solution of the system (4.7) for an arbitrary Markoff triple, and vice versa.
The first identity in Proposition 1.2 b) ensures that property (4.9) is
preserved under such a transformation. $\Box$

\begin{tabular}{l}
  
\end{tabular}

Our next task is to characterize the solutions of the system (4.7) and (4.9),
whose existence has been established in Proposition 4.5, for arbitrary Markoff
triples more specifically.

\begin{tabular}{l}
  
\end{tabular}

{\tmstrong{4.6 Proposition}} Suppose that $q \begin{array}{l}
  \epsilon
\end{array} \mathbb{N}$ is divisible by 3 and $- 1$ is a quadratic residue
modulo $\frac{q}{3}$. Then given an integral solution $X$ of the system (4.7)
and (4.9), there exist two integers $\varepsilon$ and $j$ such that
$\varepsilon^2 + 1 = \frac{q}{3} j$, and there exist three integers $\alpha,
k, l$, such that $c \alpha - a = m k$ and $\alpha^2 + 1 = \frac{m}{3} l$,
having the following property,

(4.22)
\[ Z X = Q =\mathcal{A}\mathcal{B}, \]
where
\[ Z = \left(\begin{array}{ccc}
     c & m & a\\
     - a & 0 & c\\
     a & 2 & c
   \end{array}\right), Q = \left(\begin{array}{ccc}
     0 \begin{array}{l}
       
     \end{array} & q & 0\\
     0 \begin{array}{l}
       
     \end{array} & m \varepsilon - q \alpha & 2 m q\\
     2 m \begin{array}{l}
       
     \end{array} & \frac{1}{3} (q l + m j) - 2 \alpha \varepsilon
     \begin{array}{l}
       
     \end{array} & 4 (m \varepsilon - q \alpha)
   \end{array}\right) \]
\[ \mathcal{A}= \left(\begin{array}{ccc}
     1 & 0 & 0\\
     - \alpha & m & 0\\
     \frac{l}{3} & - 2 \alpha & m
   \end{array}\right), \mathcal{B}= \left(\begin{array}{ccc}
     0 & q & 0\\
     0 & \varepsilon & 2 q\\
     2 & \frac{j}{3} & 4 \varepsilon
   \end{array}\right) . \]
Conversely, any integral solution of the form (4.22) is also a solution of the
system (4.7) and (4.9).

\begin{tabular}{l}
  
\end{tabular}

{\tmstrong{Remarks}} 1) Notice the separation of data pertaining to the
quadratic residues for the numbers $m \tmop{and} q$, respectively, which
results from the factorization of the matrix $Q$ into $\mathcal{A} \tmop{and}
\mathcal{B}$.

2) If $(a, b, c)$ is an arbitrary triple of positive integers admitting a
solution that can be factored as in (4.22), then $(a, b, c)$ is a Markoff
triple. This is a consequence of the identity
\[ X^t \left(\begin{array}{c}
     c \\
     m\\
     a 
   \end{array}\right) = \left(\begin{array}{c}
     0\\
     q\\
     0
   \end{array}\right), \]
as well as the observation that the first entry in this vector identity is
equivalent to the Markoff property.

3) Denoting the $i$-th column vector of the matrix $X$ in (4.22) by $x^{(i)}$,
the following two identities hold
\[ R (a, b, c) x^{(2)} = x^{(3)}, R (a, b, c) x^{(3)} = 4 q x^{(1)} . \]

\

Before we turn to the proof of Proposition 4.6 we state some consequences and
identities. We call two solutions$X_1 \tmop{and} X_2$ of the system (4.7) and
(4.9) equivalent if and only if
\[ e^{\frac{i}{2} \mathcal{R}} X_1 = X_2, \tmop{where} \mathcal{R}=
   \frac{1}{3} \left(\begin{array}{ccccc}
     - c^2  & 2 a - c m &  &  & 2 b - a c\\
     b c - 2 a & c^2 - a^2  &  &  & 2 c - a b\\
     a c - 2 b & a m - 2 c &  &  & a^2
   \end{array}\right), \]
for some integer $i$. The following restates part b) of Proposition 4.3 and
Part b) of Corollary 4.5.

\begin{tabular}{l}
  
\end{tabular}

{\tmstrong{4.7 Corollary}} Two solutions of the system (4.7) and (4.9) are
equivalent if and only if they are associated with two numbers $\varepsilon_1
\tmop{and} \varepsilon_2$ which are congruent modulo $\frac{q}{3}$. In
particular the number of inequivalent integral solutions of the system (4.7)
and (4.9) is the same for all Markoff triples.

\begin{tabular}{l}
  
\end{tabular}

{\tmstrong{Proof}} First we observe that for $R = 3\mathcal{R}$
\[ R (Z^{- 1} \mathcal{A}\mathcal{B}) = (Z^{- 1} \mathcal{A}\mathcal{B})
   \left(\begin{array}{ccc}
     0 & 0 & 4 q\\
     0 & 0 & 0\\
     0 & 1 & 0
   \end{array}\right) . \]
This follows from the following sequence of basic identities.
\[ R Z^{- 1} = Z^{- 1} \left(\begin{array}{ccc}
     0 & 0 & 0\\
     2 m & 0 & 0\\
     0 & 4 & 0
   \end{array}\right), \]
\[ \left(\begin{array}{ccc}
     0 & 0 & 0\\
     m & 0 & 0\\
     0 & 2 & 0
   \end{array}\right) \mathcal{A}=\mathcal{A} \left(\begin{array}{ccc}
     0 & 0 & 0\\
     1 & 0 & 0\\
     0 & 2 & 0
   \end{array}\right), \]
\[ \left(\begin{array}{ccc}
     0 & 0 & 0\\
     1 & 0 & 0\\
     0 & 2 & 0
   \end{array}\right) \mathcal{B}= \frac{1}{2} \mathcal{B}
   \left(\begin{array}{ccc}
     0 & 0 & 4 q\\
     0 & 0 & 0\\
     0 & 1 & 0
   \end{array}\right) . \]
Next we observe, writing temporarily $\mathcal{B}=\mathcal{B}_{\varepsilon}$,
\[ \mathcal{B}_{\varepsilon} e x p (\frac{x}{2} \left(\begin{array}{ccc}
     0 & 0 & 4 q\\
     0 & 0 & 0\\
     0 & 1 & 0
   \end{array}\right)) =\mathcal{B}_{\varepsilon + x q}, \]
where $e x p (x) = e^x $. Putting all this together, yields
\[ e x p (\frac{1}{2} x\mathcal{R}) (Z^{- 1}
   \mathcal{A}\mathcal{B}_{\varepsilon}) = Z^{- 1}
   \mathcal{A}\mathcal{B}_{\varepsilon + x \frac{q}{3}}, \]
from which the claim follows. $\Box$

\

\begin{tabular}{l}
  
\end{tabular}

{\tmstrong{4.8 Corollary}} Given any integers $\alpha, l$, such that \
$\alpha^2 + 1 = \frac{m}{3} l$ the following identities hold.

(4.23)
\[ (Z^{- 1} \mathcal{A})^t M (a, b, c) (Z^{- 1} \mathcal{A}) = \frac{1}{2}
   \left(\begin{array}{ccc}
     0 & 1 & 1\\
     - 1 & - 2 & 0\\
     1 & 0 & 0
   \end{array}\right) . \]
Moreover, the numbers \ $\alpha, k, l$ can be chosen such that

(4.24)
\[ \mathcal{A}^{- 1} Z \begin{array}{l}
     \epsilon
   \end{array}  \frac{1}{3} \tmmathbf{M}_3 \text{(} \mathbb{Z})
   \begin{array}{l}
     
   \end{array} \tmop{if} \begin{array}{l}
     
   \end{array} m \begin{array}{l}
     
   \end{array} \tmop{is} \begin{array}{l}
     
   \end{array} \tmop{odd}, \begin{array}{l}
     
   \end{array} \mathcal{A}^{- 1} Z \begin{array}{l}
     \epsilon
   \end{array}  \frac{1}{6} \tmmathbf{M}_3 \text{(} \mathbb{Z})
   \begin{array}{l}
     
   \end{array} \tmop{if} \begin{array}{l}
     
   \end{array} m \begin{array}{l}
     
   \end{array} \tmop{is} \begin{array}{l}
     
   \end{array} \tmop{even} . \]
\begin{tabular}{l}
  
\end{tabular}

{\tmstrong{Proof}} The identity (4.23) follows from the identity (4.6), as
well as the identity

(4.25)
\[ \mathcal{A}^t \left(\begin{array}{ccc}
     - 2 & m & m\\
     - m & - 2 & 0\\
     m & 0 & 0
   \end{array}\right) \mathcal{A}= m^2 \left(\begin{array}{ccc}
     0 & 1 & 1\\
     - 1 & - 2 & 0\\
     1 & 0 & 0
   \end{array}\right) . \]
In order to show the validity of (4.24) we choose $q$ in Proposition 4.6 such
that $\frac{q}{3} = 1$ modulo 4 and $\frac{q}{3}$ is a prime number which does
not divide $m$. Then, by (4.22),
\[ \mathcal{A}^{- 1} Z = \frac{1}{\det (\mathcal{A})} \mathcal{A}^{\tmop{adj}}
   Z = \frac{1}{m^2} \mathcal{A}^{\tmop{adj}} Z = \frac{1}{2 q^2}
   \mathcal{B}X^{\tmop{adj}} = \frac{1}{\det (X)} \mathcal{B}X^{\tmop{adj}}
   =\mathcal{B}X^{- 1} . \]
Since $\frac{q}{3}$ does not divide $\frac{m}{3}$ by assumption, it follows
that none of the denominators of the reduced fractions in the entries of the
matrix $\mathcal{A}^{- 1} Z$ is divisible by a prime factor of $m$ distinct
from 2 or 3. Furthermore, choosing $\alpha$ and $k$ as in Proposition 4.6,
\[ c \alpha - a = m k \text{ ,} \]
which, by virtue of the Markoff property implies that there exists an integer
$k^{\ast}$ such that
\[ a \alpha + c = m k^{\ast}, \]
we obtain
\[ \mathcal{A}^{- 1} Z = \frac{1}{m^2} \left(\begin{array}{ccc}
     m^2 & 0 & 0\\
     \alpha m & m & 0\\
     \frac{m l}{3} - 2 & 2 \alpha & m
   \end{array}\right) \left(\begin{array}{ccc}
     c & m & a\\
     - a & 0 & c\\
     a & 2 & c
   \end{array}\right) = \]
\[ \frac{1}{m^2} \left(\begin{array}{ccc}
     m^2 c & m^3 & a m^2\\
     m (c \alpha - a) & \alpha m^2 & m (a \alpha + c)\\
     m (\frac{c l}{3} + a) - 2 (a \alpha + c) \begin{array}{l}
       
     \end{array} & \begin{array}{l}
       
     \end{array} \frac{m^2 l}{3} \begin{array}{l}
       
     \end{array} & m (\frac{a l}{3} + c) + 2 (c \alpha - a)
   \end{array}\right) = \]
\[ \frac{1}{m } \left(\begin{array}{ccc}
     m  c & m^2 & a m \\
     c \alpha - a & \alpha m  & a \alpha + c\\
     \frac{c l}{3} + a - 2 k^{\ast} \begin{array}{l}
       
     \end{array} & \begin{array}{l}
       
     \end{array} \frac{m  l}{3} \begin{array}{l}
       
     \end{array} & \frac{a l}{3} + c + 2 k
   \end{array}\right), \]
and so the denominators of the reduced fractions in the entries of
$\mathcal{A}^{- 1} Z$ divide 3 if $m$ is odd, and they divide 6 if $m$ is
even. $\Box$

\begin{tabular}{l}
  
\end{tabular}

{\tmstrong{Remark }}Swapping the roles of the matrices $W \tmop{and} Z$ and,
accordingly $\mathcal{A} \tmop{and} \mathcal{B}$, one can obtain an identity
akin to (4.23). Combining (4.2) with the identity,

(4.26)
\[ (\mathcal{B}^{- 1})^t \left(\begin{array}{ccc}
     0 & q & 0\\
     q & 1 & q^2\\
     0 & - q^2 & - 4 q^2
   \end{array}\right) \mathcal{B}^{- 1} = \frac{1}{2} \left(\begin{array}{ccc}
     0 & 1 & 1\\
     - 1 & - 2 & 0\\
     1 & 0 & 0
   \end{array}\right), \]
for $q = m$leads to

(4.27)
\[ (W\mathcal{B}^{- 1})^t M (a, b, c) (W\mathcal{B}^{- 1}) = \frac{1}{2}
   \left(\begin{array}{ccc}
     0 & 1 & 1\\
     - 1 & - 2 & 0\\
     1 & 0 & 0
   \end{array}\right) . \]
The following two identities shed some more light on the nature of the
matrices $\mathcal{A} \tmop{and} \mathcal{B}$.

\begin{tabular}{l}
  
\end{tabular}

{\tmstrong{4.9 Lemma}} If $\alpha_i, l_i, \varepsilon_i, j_i ; i = 1, 2$ are
two sets of data as in Proposition 4.6, and $\mathcal{A}_i, \mathcal{B}_i$ are
the corresponding matrices associated with them, then

(4.28)
\[ \mathcal{A}^{- 1}_1 \mathcal{A}_2 = e x p ((\frac{\alpha_1 - \alpha_2}{m})
   \left(\begin{array}{ccc}
     0 & 0 & 0\\
     1 & 0 & 0\\
     0 & 2 & 0
   \end{array}\right)), \mathcal{A} _2 \mathcal{A}^{- 1}_1 = e x p
   ((\frac{\alpha_1 - \alpha_2}{m}) \left(\begin{array}{ccc}
     0 & 0 & 0\\
     m & 0 & 0\\
     0 & 2 & 0
   \end{array}\right)), \]
(4.29)
\[ \mathcal{B} _2 \mathcal{B}^{- 1}_1 = e x p ((\frac{\varepsilon_2 -
   \varepsilon_1}{q}) \left(\begin{array}{ccc}
     0 & 0 & 0\\
     1 & 0 & 0\\
     0 & 2 & 0
   \end{array}\right)), \mathcal{B}^{- 1}_1 \mathcal{B}_2 = e x p
   ((\frac{\varepsilon_2 - \varepsilon_1}{q}) \left(\begin{array}{ccc}
     0 & 0 & 4 q\\
     0 & 0 & 0\\
     0 & 1 & 0
   \end{array}\right)) . \]

\

The proof of Lemma 4.7 is obtained through straightforward manipulations.

\begin{tabular}{l}
  
\end{tabular}

{\tmstrong{Proof of Proposition 4.6}} That any matrix $X$ having the
properties stipulated in \ (4.22) is a solution of the system (4.7) and (4.9)
can be seen through an application of (4.6), (4.25) and (4.26). In order to
show that any integral solution of (4.7) and (4.9) has the claimed form, we
are going to proceed as in the line of reasoning leading up to Proposition
4.3. This means that we need to solve the system (4.7) and (4.9) in such a way
as to exhibit the dependence of the solutions on the data pertaining to the
corresponding quadratic residues. Extending (4.10) to the case of a general
Markoff triple we consider,

(4.30)
\[ \left(\begin{array}{ccc}
     x_{11} & x_{21} & x_{31}\\
     x_{12} & x_{22} & x_{32}\\
     x_{13} & x_{23} & x_{33}
   \end{array}\right) \left(\begin{array}{ccc}
     1 & a & b\\
     0 & 1 & c\\
     0 & 0 & 1
   \end{array}\right) \left(\begin{array}{ccc}
     x_{11} & x_{12} & x_{13}\\
     x_{21} & x_{22} & x_{23}\\
     x_{31} & x_{32} & x_{33}
   \end{array}\right) = \]
\[ \left(\begin{array}{ccc}
     x^2_{11} + a x_{11} x_{21} + x_{21}^2 & x_{11} x_{12} + a x_{11} x_{22} +
     x_{21} x_{22} \begin{array}{l}
       
     \end{array} & x_{11} x_{13} + a x_{11} x_{23} + x_{21} x_{23}\\
     + b x_{11} x_{31} + c x_{21} x_{31} + x^2_{31} \begin{array}{l}
       
     \end{array} & + b x_{11} x_{32} + c x_{21} x_{32} + x_{31} x_{32}
     \begin{array}{l}
       
     \end{array} \begin{array}{l}
       
     \end{array} & + b x_{11} x_{33} + c x_{21} x_{33} + x_{31} x_{33}\\
     &  & \\
     x_{11} x_{12} + a x_{12} x_{21} + x_{21} x_{22} \begin{array}{l}
       
     \end{array} & x^2_{12} + a x_{12} x_{22} + x^2_{22} & x_{12} x_{13} + a
     x_{12} x_{23} + x_{22} x_{23}\\
     + b x_{12} x_{31} + c x_{22} x_{31} + x_{31} x_{32} \begin{array}{l}
       
     \end{array} & + b x_{12} x_{32} + c x_{22} x_{32} + x^2_{32} & + b x_{12}
     x_{33} + c x_{22} x_{33} + x_{32} x_{33}\\
     &  & \\
     x_{11} x_{13} + a x_{13} x_{21} + x_{21} x_{23} \begin{array}{l}
       
     \end{array} \begin{array}{l}
       
     \end{array} & x_{12} x_{13} + a x_{13} x_{22} + x_{22} x_{23} & x^2_{13}
     + a x_{13} x_{23} + x^2_{23}\\
     + b x_{13} x_{31} + c x_{23} x_{31} + x_{31} x_{33} \begin{array}{l}
       
     \end{array} & + b x_{13} x_{32} + c x_{23} x_{32} + x_{32} x_{33} & + b
     x_{13} x_{33} + c x_{23} x_{33} + x^2_{33}
   \end{array}\right) \]
\[ = \left(\begin{array}{ccc}
     0 & q & 0\\
     q & 1 & q^2\\
     0 & - q^2 & - 4 q^2
   \end{array}\right) . \]
By (4.9) we have

(4.31)
\[ x_{11} = c, x_{21} = - b, x_{31} = a . \]
Entry (1,2) yields

(4.32)
\[ c x_{12} + m x_{22} + a x_{32} = q . \]
Entry (1,3) yields

(4.33)
\[ c x_{13} + m x_{23} + a x_{33} = 0 . \]
Multiplying the (3,3) entries by $m^2$ and eliminating $x_{23}$ yields $c
x_{33} - a x_{13} = \pm 2 m q$. In order to determine the proper sign we
multiply the (2,3) entries by $m^2$, eliminate $x_{23}$ by means of (4.33),
then we proceed in exactly the same fashion with the (3,2) entries, and obtain
after subtracting the latter from the former,

(4.34)
\[ c x_{33} - a x_{13} = 2 m q . \]
Multiplying the (2,2) entries by $m^2$ and combining the result with (4.32)
yields,

(4.35)
\[ (a x_{12} - c x_{32})^2 - (a x_{12} + 2 x_{22} + c x_{32}) m q = - m^2 -
   q^2 . \]
We consider the following factorization,

(4.36)
\[ m = p d, q = p e ; d, e \tmop{and} p \tmop{pairwise} \tmop{relatively}
   \tmop{prime} . \]
Then (4.35) implies that $a x_{12} - c x_{32}$ is divisible by $p$. Since $d
\tmop{and} e$ are relatively prime, there exist integers $\varepsilon_0$ and
$\alpha_0$ such that
\[ e \alpha_0 - d \varepsilon_0 = \frac{a x_{12} - c x_{32}}{p} . \]
or

(4.37)
\[ q \alpha_0 - m \varepsilon_0 = a x_{12} - c x_{32} \]
Combining (4.35) and (4.37) yields

(4.38)
\[ m^2 (\varepsilon^2_0 + 1) + q^2 (\alpha^2_0 + 1) = (a x_{12} + 2 x_{22} + c
   x_{32} + 2 \varepsilon_0 \alpha_0) m q. \]
It follows from this that there exist integers $j_0$ and $k_0$ such that

(4.39)
\[ \varepsilon^2_0 + 1 = e j_0, \alpha^2_0 + 1 = d k_0 . \]
Now (4.37) implies

(4.40)
\[ \alpha_0 \equiv \frac{1}{q} (a x_{12} - c x_{32}) \quad (\tmop{mod} d), \]
while (4.32) implies

(4.41)
\[ x_{32} \equiv \frac{1}{a} (q - c x_{12}) \begin{array}{l}
     
   \end{array}  (\tmop{mod} d) . \]
Combining (4.40) and (4.41) we obtain by virtue of the Markoff property
\[ \alpha_0 \equiv \frac{1}{q} (a x_{12} - \frac{c}{a} (q - c x_{12})) \equiv
   (a + \frac{c^2}{a}) \frac{x_{12}}{q} - \frac{c}{a} \equiv \frac{m b}{a q}
   x_{12} - \frac{c}{a} \equiv - \frac{c}{a} \begin{array}{l}
     
   \end{array} (\tmop{mod} d) . \]
Hence,
\[ a \alpha_0 + a \equiv 0 \begin{array}{l}
     
   \end{array} (\tmop{mod} d), \]
or, again by the Markoff property,

(4.42)
\[ c \alpha_0 - a \equiv 0 \begin{array}{l}
     
   \end{array} (\tmop{mod} d) . \]
Let $\alpha$ and $k$ be integers such that

(4.43)
\[ c \alpha - a = m k, \]
and let $l$ be an integer such that

(4.44)
\[ \alpha^2 + 1 = \frac{m}{3} l . \]
Now (4.42) and (4.43) imply that there exists an integer $u$ such that

(4.45)
\[ \alpha - \alpha_0 = d u . \]
Let

(4.46)
\[ \varepsilon = \varepsilon_0 + e u . \]
Then (4.37), (4.45) and (4.46) yield
\[ q \alpha - m \varepsilon = q \alpha - m \varepsilon_0 - p d e u = q \alpha
   - m \varepsilon_0 - q d u = q \alpha - m \varepsilon_0 - q (\alpha -
   \alpha_0) = q \alpha_0 - m \varepsilon_0 \]
\[ = a x_{12} - c x_{32} . \]
In conclusion

(4.47)
\[ q \alpha - m \varepsilon = a x_{12} - c x_{32} . \]
Moreover,
\[ \varepsilon^2 + 1 = \varepsilon^2_0 + 2 e \varepsilon_0 u + e^2 u^2 + 1 = e
   (j_0 + 2 \varepsilon_0 u + e u^2), \]
or, letting $j^{\ast} = j_0 + 2 \varepsilon_0 u + e u^2$,

(4.48)
\[ \varepsilon^2 + 1 = e j^{\ast} . \]
Since (4.35) and (4.47) yield
\[ m^2 (\varepsilon^2 + 1) + q^2 (\alpha^2 + 1) = (a x_{12} + 2 x_{22} + c
   x_{32} + 2 \varepsilon \alpha) m q, \]
we finally obtain by virtue of (4.44) and (4.48),

(4.49)
\[ a x_{12} + 2 x_{22} + c x_{32} = \frac{q}{3} l + d j^{\ast} - 2 \alpha
   \varepsilon \]
Putting it all together, (4.32) provides the first, (4.47) the second, and
(4.49) the third linear identity, respectively, of the following system

(4.50)
\[ \left(\begin{array}{ccc}
     c & m & a\\
     - a & 0 & c\\
     a & 2 & c
   \end{array}\right) \left(\begin{array}{c}
     x_{12}\\
     x_{22}\\
     x_{32}
   \end{array}\right) = \left(\begin{array}{c}
     q\\
     m \varepsilon - q \alpha\\
     \frac{q}{3} l + d j^{\ast} - 2 \alpha \varepsilon
   \end{array}\right) \]
We turn now to the third column of the matrix $X$. We have found two linear
identities already, namely (4.33) and (4.34). To determine the third, we
multiply the (2,3) entries by $m^2$, which simplifies to
\[ (a x_{12} - c x_{32}) (c x_{33} - a x_{13}) + (x_{23} + c x_{33}) m q = m^2
   q^2 . \]
Combining this with (4.34) yields
\[ 2 (a x_{12} - c x_{32}) + (x_{23} + c x_{33}) = m q, \]
or by (4.47),
\[ x_{23} + c x_{33} = m q + 2 (m \varepsilon - q \alpha) . \]
Multiplying this by 2 and subtracting (4.34) from it, yields
\[ a x_{13} + 2 x_{23} + c x_{33} = 4 (m \varepsilon - q \alpha), \]
which is the missing third identity. Putting it all together again, we have

(4.51)
\[ \left(\begin{array}{ccc}
     c & m & a\\
     - a & 0 & c\\
     a & 2 & c
   \end{array}\right) \left(\begin{array}{c}
     x_{13}\\
     x_{23}\\
     x_{33}
   \end{array}\right) = \left(\begin{array}{c}
     0\\
     2 m q\\
     4 (m \varepsilon - q \alpha)
   \end{array}\right) . \]
Finally, Markoff's property ensures that (4.31) implies

(4.52)
\[ \left(\begin{array}{ccc}
     c & m & a\\
     - a & 0 & c\\
     a & 2 & c
   \end{array}\right) \left(\begin{array}{c}
     x_{11}\\
     x_{21}\\
     x_{31}
   \end{array}\right) = \left(\begin{array}{c}
     0\\
     0\\
     2 m
   \end{array}\right) . \]
Next we shall first deal with a special case, namely the situation where the
numbers $\frac{m}{3}$ and $\frac{q}{3}$ are relatively prime. In this case $p
= 3$ and \ $d = \frac{m}{3}$, and so, letting $j = j^{\ast}$, (4.50) reads
\[ \left(\begin{array}{ccc}
     c & m & a\\
     - a & 0 & c\\
     a & 2 & c
   \end{array}\right) \left(\begin{array}{c}
     x_{12}\\
     x_{22}\\
     x_{32}
   \end{array}\right) = \left(\begin{array}{c}
     q\\
     m \varepsilon - q \alpha\\
     \frac{1}{3} (q l + m j) - 2 \alpha \varepsilon
   \end{array}\right) . \]

This, combined with (4.51) and (4.52) yields the claimed identity $Z X = Q$,
as well as the claimed factorization of $Q$ into the matrices $\mathcal{A}$
and $\mathcal{B}$, thus settling the claim in this particular case. In order
to deal with the general case we observe that the combination of (4.50)
through (4.52) yields the factorization
\[ Z X =\mathcal{A} \left(\begin{array}{ccc}
     0 & q & 0\\
     0 & \varepsilon & 2 q\\
     2 & \frac{j^{\ast}}{p} & 4 \varepsilon
   \end{array}\right) . \]
Solving for the second factor on the right hand side yields

(4.53)
\[ \left(\begin{array}{ccc}
     0 & q & 0\\
     0 & \varepsilon & 2 q\\
     2 & \frac{j^{\ast}}{p} & 4 \varepsilon
   \end{array}\right) =\mathcal{A}^{- 1} Z X . \]
Next we are going to invoke (4.24) in Corollary 4.8, whose proof was based on
the special case \ we have just settled. By (4.24) the reduced fractions in
the entries of the matrix $\mathcal{A}^{- 1} Z$are either integers or rational
numbers whose denominator divides 6. Since $X$ is integral, the same is true
for the entries of the matrix $\mathcal{A}^{- 1} Z X$ on the right hand side
of (4.53). Since the denominator of the reduced fraction in the entry (3,2) of
the matrix on the left hand side of (4.53) has exactly one factor 3, which,
due to the fact that $j^{\ast}$ as a product of prime factors which are either
equal to 2 or equal to 1 modulo 4, does not cancel, that denominator can only
be equal to 3 or 6. If at least one of the integers $m$ or $q$ is odd, then
$p$ has to be odd as well, in which case the said denominator is equal to 3.
If both, $m$ and $q$ are even, then (4.47) implies that the integers $a
x_{12}$ and $c x_{32}$ are both either odd or even. Since $m$, being a Markoff
number, can have at most one even prime factor, and since this implies that
$d$ has to be odd, this together with (4.49) implies that $j^{\ast}$ has to be
even. In conclusion, the said denominator can not be 6, and therefore it has
to be equal to 3. So, letting $j = \frac{3 j^{\ast}}{p}$, our claim follows in
the general case as well. $\Box$

\begin{tabular}{l}
  
\end{tabular}

{\tmstrong{Remarks}} 1) Note that all the arguments in the proof of
Proposition 4.6 are necessary for the existence of an integral solution of the
system (4.7) and (4.9), they are not sufficient. Sufficiency rests entirely
upon Corollary 4.5. (See also the discussion at the end of this section.)
However, since the factorization obtained in Proposition 4.6 is valid for all
Markoff triples, in particular for the triple $(3, 3, 3)$, the parametrization
of the equivalence classes of all integral solutions of the system (4.7) and
(4.9) through quadratic residues of $- 1$ is unique.

2) The case $m = q$ in Proposition 4.6 calls for some special attention. In
this case we have the following particular form for the matrix product on the
right hand side of (4.22),
\[ \mathcal{A}\mathcal{B}= \left(\begin{array}{ccc}
     0 & m & 0\\
     0 & m (\varepsilon - \alpha) & 2 m^2\\
     2 m & (\varepsilon - \alpha)^2 + 2 & 4 m (\varepsilon - \alpha)
   \end{array}\right) . \]
Significantly, this matrix depends on $m$ and the integer $\varepsilon -
\alpha$ only. Since $\alpha^2 = - 1, \varepsilon^2 = - 1$ modulo $m$, it
follows that there exists a factorization $\frac{m}{3} = p q$ if $m$ is odd,
$\frac{m}{6} = p q$ if $m$ is even, with $p$ and $q$ being relatively prime
integers. This observation reflects a rather general pattern. Given any
integer $n$ with a quadratic residue of $- 1$ in the residue class ring
$\mathbb{Z}_n$, and given a particular choice of such a quadratic residue, $k$
say, the differences \ $k - k^{\ast}$, where $k^{\ast}$ ranges over all the
other quadratic residues of $- 1$ in $\mathbb{Z}_n$, correspond in a
one-to-one way to all the ordered pairs of positive integers whose product
equals $n$ divided by the product of its even prime factors. Since all initial
choices are equivalent in this regard, we can not hope to characterize a
specific one, in our case $\alpha$, and settle the uniqueness question within
this framework. The characterization of $\alpha$ among all the other quadratic
residues of $- 1$ modulo $\frac{m}{3}$ is more implicit. It is expressed
through the ``almost'' integrality of the matrix $\mathcal{A}^{- 1} Z$ and its
inverse. The first part of this statement was established in (4.24). Being of
such an implicit nature however, this property is too elusive in order to be
useful for the establishment of the uniqueness of the set $\nobracket \{
\alpha, - \alpha \}$ modulo $\frac{m}{3}$, which is equivalent to the
uniqueness claim of the Theorem.

3) The question arises how much of the formalism in this section is particular
to the setting of Markoff triples. In order to obtain a partial answer to this
question we introduce the following concept:

If $\mathfrak{p} \geqslant 3$ is an integer, then we call $(\mathfrak{a},
\mathfrak{b}, \mathfrak{c}) \begin{array}{l}
  \epsilon
\end{array} \mathbb{N}^3$ a $\mathfrak{p}-$triple if

(*)
\[ \mathfrak{a}^2 +\mathfrak{b}^2 +\mathfrak{c}^2
   =\mathfrak{p}\mathfrak{a}\mathfrak{b}\mathfrak{c}+ 3 -\mathfrak{p}. \]
The following statements hold true for $\mathfrak{p}-$triples:

(I) $(1, 1, 1)$ is a $\mathfrak{p}-$triple.

(II) If $(\mathfrak{a}, \mathfrak{b}, \mathfrak{c})$ is a
$\mathfrak{p}-$triple then $(\mathfrak{a}, \mathfrak{b},
\mathfrak{p}\mathfrak{a}\mathfrak{c}-\mathfrak{b})$ and
$(\mathfrak{p}\mathfrak{a}\mathfrak{c}-\mathfrak{b}, \mathfrak{a},
\mathfrak{c})$ are $\mathfrak{p}-$triples as well.

(III) Up to permutations of the components, any $\mathfrak{p}-$triple can be
obtained through finitely many transitions as stipulated in (II).

(IV) All solutions of the diophantine matrix equation
\[ X^t M (\mathfrak{p}\mathfrak{a}, \mathfrak{p}\mathfrak{b},
   \mathfrak{p}\mathfrak{c}) X = \left(\begin{array}{ccc}
     0 & \mathfrak{p}\mathfrak{q} & 0\\
     \mathfrak{p}\mathfrak{q} & 1 & \mathfrak{p}^2 \mathfrak{q}^2\\
     0 & -\mathfrak{p}^2 \mathfrak{q}^2 & - 4\mathfrak{p}^2 \mathfrak{q}^2
   \end{array}\right) \]
for a given integer $\mathfrak{q}$ for which $- 1$ is a quadratic residue
modulo $\mathfrak{q}$ can be parametrized by the solutions of the equation
$\varepsilon^2 \equiv - 1$\begin{tabular}{l}
  
\end{tabular}(mod $\mathfrak{q}$).

Properties (I) through (III) tell us that a tree of $\mathfrak{p}-$triples can
be built in parallel to the developments in Section 1. Property (IV) is a
reflection of the fact that Proposition 4.3a) and its proof carry over to the
settings of $\mathfrak{p}-$triples.

The case $\mathfrak{p}= 3$ yields Markoff triples of course, and the
factorization obtained in Proposition 4.6 is particular to this case. If we
choose $\mathfrak{p}= 0$ in (*), then $(\mathfrak{a}, \mathfrak{b},
\mathfrak{c})$ is a solution of the equation if and only if $| \mathfrak{a} |
= | \mathfrak{b} | = | \mathfrak{c} | = 1$. But if we choose $\mathfrak{p}= -
1$, then we recover (essentially, i. e. up to a minus sign) the case of the
Tchebycheff polynomials briefly discussed in Section 2.

\

We conclude the present section by extending some of the considerations from
the quadratic residues of $- 1$ modulo a given Markoff number $\mathfrak{m}=
\frac{m}{3}$ which are determined by a Markoff triple dominated by $m$ to
arbitrary quadratic residues of $- 1$ \ modulo a given Markoff number. We will
particularly focus on the divisibility property (4.24) in Corollary 4.8 in
this larger context. The first objective is to obtain a characterization of
the equivalence classes of these quadratic residues in terms of factorizations
of the corresponding Markoff number. Throughout we shall be using the
following notation. Let $\mathbb{m}= \frac{m}{3}$ if $m$ is odd, and
$\mathbb{m}= \frac{m}{6}$ if $m$ is even.

\begin{tabular}{l}
  
\end{tabular}

{\tmstrong{4.10 Lemma}} a) Let $n$ be an integer such that $n^2 + 1
=\mathbb{m}l$. Then there exists a (unique) ordered pair $(p, q)$ of
relatively prime positive integers $p$ and $q$ such that $\mathbb{m}= p q$ and

(4.54)
\[ c n + a = p u, \begin{array}{l}
     
   \end{array} c n - a = q v ; \begin{array}{l}
     
   \end{array} u, v \begin{array}{l}
     \epsilon
   \end{array} \mathbb{Z}. \]
b) For any ordered pair $(p, q)$ of relatively prime positive integers $p$ and
$q$ such that $\mathbb{m}= p q$ there is exactly one equivalence class of
numbers $n$ modulo $\mathfrak{m}$ such that $n^2 + 1 =\mathbb{m}l$ and (5.1)
holds.

\begin{tabular}{l}
  
\end{tabular}

{\tmstrong{Proof}} a) Since
\[ (c n + a) (c n - a) = c^2 n^2 - a^2 = c^2 (\mathbb{m}l - 1) - a^2 = c^2
   l\mathbb{m}- (a^2 + c^2) = c^2 l\mathbb{m}- b m, \]
there exists a factorization $\mathbb{m}= p q$ such that
\[ \frac{c n + a}{p}, \frac{c n - a}{q} \begin{array}{l}
     \epsilon
   \end{array} \mathbb{Z}. \]
Since both, $p$ and $q$ are relatively prime to $a c$ as well as to $n$, the
numbers $p$ and $q$ have to be relatively prime as well, and therefore the
ordered pair $(p, q)$ is uniquely determined.

b) By Lemma 4.2 part (a) and (b) the number of equivalence classes of integers
for which $- 1$ is a quadratic residue modulo $\mathbb{m}$, and the number of
ordered pairs $(p, q)$ of relatively prime integers $p$ and $q$ such that
$\mathbb{m}= p q$, where $p$ and $q$ are products of odd primes only, is the
same. Since part a) ensures already that we have an injective map from one set
into the other, the claim follows. $\Box$

\

{\tmstrong{Remarks}} 1) Note that the integers $p \tmop{and} v$ occurring in
(4.54) are relatively prime. To see this, suppose the opposite were true,
namely that both $p \tmop{and} v$ are divisible by a prime factor $r$. This
prime factor has to be larger or equal to 5, because $\mathbb{m}$ has only
prime factors larger or equal to 5 and hence $p$, being a factor of
$\mathbb{m}$, has to be larger or equal to 5 as well. \ Subtracting the \
second identity in (4.54) from the first shows that $r$ divides $2 a$, which
in turn implies that $r$ divides $\frac{a}{3}$. Adding the two identities in
(4.54) together shows that $r$ divides $c n$. But since \ $n^2 + 1
=\mathbb{m}l$, $r$ cannot be a factor of $n$, and therefore $r$ has to divide
$\frac{c}{3}$. Since $\frac{a}{3}$ and $\frac{c}{3}$ are relatively prime, we
have reached a contradiction. By the same line of reasoning one can show that
the integers $q \tmop{and} u$ in (4.54) are relatively prime.

2) Note that, in the context of Lemma 4.10, if the ordered pair $(p, q)$
corresponds to the integer $n$, then the ordered pair $(q, p)$ corresponds to
the integer $- n$.

3) The notion of ``relatively prime integers'' employed in Lemma 4.10 does
include the case $\begin{array}{l}
  p = 1
\end{array} \tmop{or} q = 1$. The quadratic residues that correspond to this
choice are exactly the two with opposite signs determined by the Markoff
numbers $a \tmop{and} c$.

\

Next we consider the quadratic polynomial
\[ F (n) =\mathfrak{c}n^2 + (3\mathfrak{m}\mathfrak{c}- 2\mathfrak{a}) n
   -\mathfrak{c} \]
for an integer parameter $n$.

\

{\tmstrong{4.11 Lemma}} The following identity holds true
\[ F (n) F (- n) =\mathfrak{c}^2 (n^2 + 1)^2 + (9\mathfrak{c}^2 - 4)
   \mathfrak{m}^2 n^2 . \]
In particular $F (n) F (- n)$ is divisible by $\mathfrak{m}^2$ if and only if
$n$ is a quadratic residue of $- 1$ modulo $\mathfrak{m}$

\

{\tmstrong{Proof}} \ Invoking once the Markoff property we obtain,
\[ F (n) F (- n) = [\mathfrak{c}n^2 + (3\mathfrak{m}\mathfrak{c}-
   2\mathfrak{a}) n -\mathfrak{c}] [\mathfrak{c}n^2 -
   (3\mathfrak{m}\mathfrak{c}- 2\mathfrak{a}) n -\mathfrak{c}] = \]
\[ \mathfrak{c}^2 n^4 - [(3\mathfrak{m}\mathfrak{c}- 2\mathfrak{a})^2 +
   2\mathfrak{c}^2] n^2 +\mathfrak{c}^2 = \]
\[ \mathfrak{c}^2 n^4 - (9\mathfrak{m}^2 \mathfrak{c}^2 -
   12\mathfrak{a}\mathfrak{c}\mathfrak{m}+ 4\mathfrak{a}^2 + 2\mathfrak{c}^2)
   n^2 +\mathfrak{c}^2 = \]
\[ \mathfrak{c}^2 n^4 - [9\mathfrak{m}^2 \mathfrak{c}^2 - 4 (\mathfrak{a}^2
   +\mathfrak{c}^2 +\mathfrak{m}^2) + 4\mathfrak{a}^2 + 2\mathfrak{c}^2] n^2
   +\mathfrak{c}^2 = \]
\[ (n^4 + 2 n^2 + 1) \mathfrak{c}^2 - (9\mathfrak{c}^2 - 4) \mathfrak{m}^2 n^2
   = \]
\[ (n^2 + 1)^2 \mathfrak{c}^2 - (9\mathfrak{c}^2 - 4) \mathfrak{m}^2 n^2 . \]
$\Box$

\

If $n$ is a quadratic residue of $- 1$ modulo $\mathfrak{m}$, more
specifically $n^2 + 1 =\mathfrak{m}l$, then $3 F (n)$ is equal to the (1,1)
entry of the matrix
\[ \mathbb{A}= \left(\begin{array}{ccc}
     - 2 (c + a n) + m c (\frac{l}{3} + n) \begin{array}{l}
       
     \end{array} & 2 m (a - c n) - m^2 c \begin{array}{l}
       
     \end{array} & m^2 c\\
     2 a c - (c^2 - a^2) n - m b \frac{l}{3} \begin{array}{l}
       
     \end{array} & m (c^2 - a^2 + 2 b n) \begin{array}{l}
       
     \end{array} & - m^2 b\\
     2 (c n - a) + m a (\frac{l}{3} - n) \begin{array}{l}
       
     \end{array} & m^2 a - 2 m (c + a n) \begin{array}{l}
       
     \end{array} & m^2 a
   \end{array}\right), \]
which is obtained from the matrix $Z^{\tmop{adj}} \mathcal{A}$ by replacing
$\alpha$ by $n$. If $(p, q)$ is affiliated with $n$ as in Lemma 4.10, then it
follows from Lemma 4.11, after invoking the Markoff property, as well as from
the first remark following Lemma 4.10, that $F (n)$ is divisible by $q^2$,
while $F (- n)$ is divisible by $p^2$. In particular the (1,1) entry in the
matrix $\mathbb{A}$ is divisible by $q^2$. The following statement shows then
that all entries in the first column of $\mathbb{A}$ are divisible by $q^2$.

\

{\tmstrong{4.12 Lemma}} If $q$ is a factor of $\mathbb{m}$ then $q^2$ divides
any of the three entries in the first column of the matrix $\mathbb{A}$, if
and only if $q^2$ divides the other two entries in the same column as well.

\

{\tmstrong{Proof}} This is a consequence of the following two observations.
Subtracting $a$ times the first entry in the first column of $\mathbb{A}$ from
$c$ times the third entry in the same column yields,
\[ 2 (a^2 + c^2 - m a c) \alpha = 2 (m b - m a c) \alpha = - 2 m^2 \alpha . \]
Likewise, adding $a$ times the first entry and the second entry in the first
column of $\mathbb{A}$ together yields,
\[ (- a^2 - c^2 + m a c) \alpha + m^2 \frac{l}{3} = m^2 (\alpha + \frac{l}{3})
   . \]
$\Box$

\

Since all the entries in the second and the third column of $\mathbb{A}$ are
obviously divisible by $q^2$, we conclude that all entries of $\mathbb{A}$ are
divisible by $q^2$. In case the first number of the pair $(p, q)$ is equal to
1, we obtain a short proof of the divisibility property (4.24), without the
restriction that $m$ be dominant.

\

{\tmstrong{Remark}} If $m = q$ in Proposition 4.6, then a more specific
description of the solutions of the system (4.7) and (4.9) can be obtained. If
$\mathbb{m}= p q$, then the corresponding solutions are of the form
\[ X = e^{\frac{n}{18 q} R} W, \]
where $n$ solves the diophantine equation
\[ q\mathfrak{c}n + p x = 2\mathfrak{a}. \]
One can see this by exploiting the fact that the second column vector of the
matrix $W$ is a unit vector. Each equivalence class of solutions contains
solutions for which the first and the third entry in the second column is
divisible by $m$. For these solutions the parameter $n$ is characterized by
the following diophantine equation,
\[ q^2 \mathfrak{c}n + p^2 x = 2\mathfrak{a}- 3\mathfrak{m}\mathfrak{c}. \]
By the Markoff property, this diophantine equation is equivalent to the
following,
\[ q^2 \mathfrak{a}n + p^2 y = - 2\mathfrak{c}+ 3\mathfrak{m}\mathfrak{a}. \]
\begin{tabular}{l}
  
\end{tabular}

\begin{tabular}{l}
  
\end{tabular}

{\tmem{{\tmstrong{5 Determination of the matrix
$Z^{\tmop{adj}}$}}$\mathcal{A}$}}

\begin{tabular}{l}
  
\end{tabular}

In this section we embark on a refined analysis of the matrix

(5.1)
\[ Z^{\tmop{adj}} \mathcal{A}= \left(\begin{array}{ccc}
     - 2 c & 2 a - m c & m c\\
     2 a c & c^2 - a^2 & - m b\\
     - 2 a & m a - 2 c & m a
   \end{array}\right) \left(\begin{array}{ccc}
     1 & 0 & 0\\
     - \alpha & m & 0\\
     \frac{l}{3} & - 2 \alpha & m
   \end{array}\right) = \]
\[ \left(\begin{array}{ccc}
     - 2 (c + a \alpha) + m c (\frac{l}{3} + \alpha) \begin{array}{l}
       
     \end{array} & 2 m (a - c \alpha) - m^2 c \begin{array}{l}
       
     \end{array} & m^2 c\\
     2 a c - (c^2 - a^2) \alpha - m b \frac{l}{3} \begin{array}{l}
       
     \end{array} & m (c^2 - a^2 + 2 b \alpha) \begin{array}{l}
       
     \end{array} & - m^2 b\\
     2 (c \alpha - a) + m a (\frac{l}{3} - \alpha) \begin{array}{l}
       
     \end{array} & m^2 a - 2 m (c + a \alpha) \begin{array}{l}
       
     \end{array} & m^2 a
   \end{array}\right), \]
which appeared in Proposition 4.6. Let $\mathfrak{m}= \frac{m}{3}$. Since
$\det (Z^{- 1} \mathcal{A}) = \frac{1}{2}$ it follows from (4.24) (or
alternatively from Lemma 4.11 and Lemma 4.12) that all entries in this matrix
are divisible by $\mathfrak{m}^2$. This is trivially true for the entries in
last column. By the specification of $\alpha$ in Proposition 4.6 and by the
Markoff property this is also clear for the entries in the second column. By
Lemma 4.12, the divisibility by $ \mathfrak{m}^2$ of any of the three entries
in the first column implies the divisibility by $\mathfrak{m}^2$ of the other
two entries. As we shall see, the divisibility by $\mathfrak{m}^2$ of any (and
hence all) of the three entries in the first column will be crucial for our
objectives. First, however, for the purpose of transparency, we are going to
change the notations in Proposition 4.6 for the quadratic residues, to the
effect that they reflect the respective Markoff numbers they are affiliated
with. Letting $k_m = \alpha, k_c = k$, by the Markoff property there exists an
integer $k_a$ such that

(5.2)
\[ c k_m - m k_c = a, m k_a - a k_m = c . \]
Moreover, there exist positive integers $l_a, l_m, l_c$ such that

(5.3)
\[ k^2_a + 1 =\mathfrak{a}l_a, k^2_m + 1 =\mathfrak{m}l_m, k^2_c + 1
   =\mathfrak{c}l_c . \]
Below we are going to restrict the values of these parameters. With this
notation in place we shall prove,

\begin{tabular}{l}
  
\end{tabular}

{\tmstrong{5.1 Proposition}} The following identities hold true,

(5.4)
\[ k_m l_a - k_a l_m = l_c + 3 k_c, k_c l_m - k_m l_c = l_a - 3 k_a, \]
\[ \mathfrak{m}l_a -\mathfrak{a}l_m = 2 k_c + 3\mathfrak{c}, \mathfrak{c}l_m
   -\mathfrak{m}l_c = 2 k_a - 3\mathfrak{a}. \]
\begin{tabular}{l}
  
\end{tabular}

{\tmstrong{Remarks}} 1) Rewriting the first and the last entry in the first
column of the matrix in (5.1) in terms of the notation introduced in (5.3) and
(5.4) one can see that the first identity in (5.4) implies the divisibility by
$m\mathfrak{m}$ of the first entry, while the second identity in (5.4) implies
the divisibility by $m\mathfrak{m}$ of the third entry. As the proof of
Proposition 5.1 will show, however, the identities in (5.4) are equivalent to
these two respective divisibility properties.

2) The first two identities in (5.4) appear in [F], Gesammelte Abhandlungen
Band III, p. 604 (18.) without proof.

\begin{tabular}{l}
  
\end{tabular}

The proof of Proposition 5.1 will be based on several lemmas which we are now
going to tackle. For any integer $x$ let

(5.5)
\[ k_a (x) = k_a +\mathfrak{a}x, k_m (x) = k_m +\mathfrak{m}x, k_c (x) = k_c
   +\mathfrak{c}x ; \]
(5.6)
\[ l_a (x) = l_a + 2 k_a x +\mathfrak{a}x^2, l_m (x) = l_m + 2 k_m x
   +\mathfrak{m}x^2, l_c (x) = l_c + 2 k_c x +\mathfrak{c}x^2 . \]
Abusing the terminology somewhat, we use in each instant the same symbol for
the polynomial as well as its constant term. This should not give rise to any
confusion because we will throughout denote the polynomial by the respective
symbol followed by $(x)$. If the constants $\alpha = k_m$ and $l = l_m$ in the
matrix $Z^{\tmop{adj}} \mathcal{A}$ in (5.1) are replaced by the polynomials
$k_m (x)$ and $l_m (x)$, respectively, then the resulting polynomial matrix

(5.7)
\[ \left(\begin{array}{ccc}
     - 2 (c + a k_m (x)) + m c (\frac{l_m (x)}{3} + k_m (x)) \begin{array}{l}
       
     \end{array} & 2 m (a - c k_m (x)) - m^2 c \begin{array}{l}
       
     \end{array} & m^2 c\\
     2 a c - (c^2 - a^2) k_m (x) - m b \frac{l_m (x)}{3} \begin{array}{l}
       
     \end{array} & m (c^2 - a^2 + 2 b k_m (x)) \begin{array}{l}
       
     \end{array} & - m^2 b\\
     2 (c k_m (x) - a) + m a (\frac{l_m (x)}{3} - k_m (x)) \begin{array}{l}
       
     \end{array} & m^2 a - 2 m (c + a k_m (x)) \begin{array}{l}
       
     \end{array} & m^2 a
   \end{array}\right) = \]
\[ \left(\begin{array}{ccc}
     {}[- 2 k_a (x) + c (k_m + \frac{l_m (x)}{3})] m \begin{array}{l}
       
     \end{array} & 2 m (a - c k_m (x)) - m^2 c \begin{array}{l}
       
     \end{array} & m^2 c\\
     2 a c - (c^2 - a^2) k_m (x) - m b \frac{l_m (x)}{3} \begin{array}{l}
       
     \end{array} & m (c^2 - a^2 + 2 b k_m (x)) \begin{array}{l}
       
     \end{array} & - m^2 b\\
     {}[2 k_c (x) + a (\frac{l_m (x)}{3} - k_m (x))] m \begin{array}{l}
       
     \end{array} & m^2 a - 2 m (c + a k_m (x)) \begin{array}{l}
       
     \end{array} & m^2 a
   \end{array}\right) \]
has still the property that every entry is divisible by $\mathfrak{m}^2$ for
all $x \epsilon \mathbb{Z}$. Considering the entries (1,1) and (3,1),
respectively, it follows that there exist integer valued functions $u_a (x)$
and $u_c (x)$ such that

(5.8)
\[ \mathfrak{m}^2 u_c (x) = - 2 (\mathfrak{c}+\mathfrak{a}k_m (x))
   +\mathfrak{m}\mathfrak{c} (l_m (x) + 3 k_m (x)), \]
\[ {\mathfrak{m}^2}  u_a (x) = 2 (\mathfrak{c}k_m (x) -\mathfrak{a})
   +\mathfrak{m}\mathfrak{a} (l_m (x) - 3 k_m (x)) \begin{array}{l}
     
   \end{array} \]
or alternatively

(5.9)
\[ \mathfrak{m}u_c (x) = - 2 k_a (x) +\mathfrak{c} (l_m (x) + 3 k_m (x)),
   \mathfrak{m}u_a (x) = 2 k_c (x) +\mathfrak{a} (l_m (x) - 3 k_m (x)) \]
\begin{tabular}{l}
  
\end{tabular}

{\tmstrong{5.2 Lemma}} The functions $u_a (x)$ and $u_c (x)$ are quadratic
polynomials with integral coefficients which have the form

(5.10)
\[ u_a (x) =\mathfrak{a}x^2 + (2 k_a - 3\mathfrak{a}) x + u_a, u_c (x)
   =\mathfrak{c}x^2 + (2 k_c + 3\mathfrak{c}) x + u_c, \]
where $u_a$ and $u_c$ are integers.

\begin{tabular}{l}
  
\end{tabular}

{\tmstrong{Proof}} Obviously, the expressions on the right hand side of the
two identities in (5.9) are quadratic polynomials, and therefore the \
functions $u_a (x)$ and $u_c (x)$ are quadratic polynomials as well, which, by
virtue of the divisibility properties of the entries in the matrix (5.7), must
have integral coefficients. So, if $u_a (x) = u^{(2)}_a x^2 + u^{(1)}_a x +
u_a$, then (5.5), (5.6) and the second identity in (5.9) yield
\[ \mathfrak{m} (u^{(2)}_a x^2 + u^{(1)}_a x + u_a)
   =\mathfrak{m}\mathfrak{a}x^2 + (2\mathfrak{c}+ 2\mathfrak{a}k_m -
   3\mathfrak{a}\mathfrak{m}) x + 2 k_c +\mathfrak{a}l_m - 3\mathfrak{a}k_m,
\]
which in turn, by the second identity in (5.2), leads to

(5.11)
\[ \mathfrak{m} (u^{(2)}_a x^2 + u^{(1)}_a x + u_a)
   =\mathfrak{m}\mathfrak{a}x^2 +\mathfrak{m} (2 k_a - 3\mathfrak{a}) x + 2
   k_c +\mathfrak{a}l_m - 3\mathfrak{a}k_m . \]
Comparing the coefficients for the quadratic and the linear terms in (5.11) on
both sides, respectively, we get,
\[ u^{(2)}_a =\mathfrak{a}, u^{(1)}_a = 2 k_a - 3\mathfrak{a}, \]
as claimed. The second identity in (5.10) can be settled in a similar way.
$\Box$

\begin{tabular}{l}
  
\end{tabular}

Let

(5.12)
\[ v_a (x) = u_a (x) + 3 k_a (x), v_c (x) = u_c (x) - 3 k_c (x) . \]
\begin{tabular}{l}
  
\end{tabular}

{\tmstrong{5.3 Lemma}} The following identity holds true,

(5.13)
\[ \left(\begin{array}{c}
     k_m (x) \begin{array}{l}
       
     \end{array} \begin{array}{l}
       
     \end{array} \begin{array}{l}
       
     \end{array} \begin{array}{l}
       
     \end{array} 1\\
     - 1 \begin{array}{l}
       
     \end{array} \begin{array}{l}
       
     \end{array} \begin{array}{l}
       \begin{array}{l}
         
       \end{array}
     \end{array} k_m (x)
   \end{array}\right) \left(\begin{array}{c}
     v_a (x)\\
     - v_c (x)
   \end{array}\right) = \left(\begin{array}{c}
     \begin{array}{l}
       
     \end{array} \begin{array}{l}
       3
     \end{array} \begin{array}{l}
       
     \end{array} \begin{array}{l}
       
     \end{array} \begin{array}{l}
       
     \end{array} \begin{array}{l}
       
     \end{array} l_m (x)\\
     - l_m (x) \begin{array}{l}
       
     \end{array} \begin{array}{l}
       
     \end{array} \begin{array}{l}
       
     \end{array} \begin{array}{l}
       - 3
     \end{array}
   \end{array}\right) \left(\begin{array}{c}
     k_c (x)\\
     k_a (x)
   \end{array}\right) \]

{\tmstrong{Proof}} First, the identities in (5.8) can be written as a linear
system in $\mathfrak{a}$ and $\mathfrak{c}$,
\[ \left(\begin{array}{c}
     -\mathfrak{m} (l_m (x) - 3 k_m (x)) + 2 \begin{array}{l}
       
     \end{array} \begin{array}{l}
       
     \end{array} \begin{array}{l}
       
     \end{array} \begin{array}{l}
       
     \end{array} \begin{array}{l}
       
     \end{array} \begin{array}{l}
       
     \end{array} \begin{array}{l}
       - 2 k_m (x)
     \end{array} \begin{array}{l}
       
     \end{array} \begin{array}{l}
       
     \end{array} \begin{array}{l}
       
     \end{array} \begin{array}{l}
       
     \end{array}\\
     \begin{array}{l}
       
     \end{array} \begin{array}{l}
       
     \end{array} \begin{array}{l}
       - 2 k_m (x)
     \end{array} \begin{array}{l}
       
     \end{array} \begin{array}{l}
       
     \end{array} \begin{array}{l}
       
     \end{array} \begin{array}{l}
       
     \end{array} \begin{array}{l}
       
     \end{array} \begin{array}{l}
       
     \end{array} \begin{array}{l}
       
     \end{array} \begin{array}{l}
       
     \end{array} \begin{array}{l}
       
     \end{array} \begin{array}{l}
       
     \end{array} \begin{array}{l}
       
     \end{array} \mathfrak{m} (l_m (x) + 3 k_m (x)) - 2
   \end{array}\right) \left(\begin{array}{c}
     \mathfrak{a}\\
     \mathfrak{c}
   \end{array}\right) =\mathfrak{m}^2 \left(\begin{array}{c}
     - u_a (x)\\
     u_c (x)
   \end{array}\right) . \]
The determinant of the matrix on the left hand side is equal to
\[ \mathfrak{m}^2 (9 k_m (x)^2 - l_m (x)^2) . \]
Hence solving the above linear system for $\left(\begin{array}{c}
  \mathfrak{a}\\
  \mathfrak{c}
\end{array}\right)$ yields
\[ [2 \left(\begin{array}{c}
     \begin{array}{l}
       - 1
     \end{array} \begin{array}{l}
       
     \end{array} \begin{array}{l}
       
     \end{array} k_m (x)\\
     k_m (x) \begin{array}{l}
       
     \end{array} \begin{array}{l}
       
     \end{array} \begin{array}{l}
       
     \end{array} \begin{array}{l}
       1
     \end{array} \begin{array}{l}
       
     \end{array}
   \end{array}\right) +\mathfrak{m} \left(\begin{array}{c}
     l_m (x) + 3 k_m (x) \begin{array}{l}
       
     \end{array} \begin{array}{l}
       
     \end{array} \begin{array}{l}
       
     \end{array} \begin{array}{l}
       
     \end{array} \begin{array}{l}
       
     \end{array} \begin{array}{l}
       
     \end{array} \begin{array}{l}
       
     \end{array} \begin{array}{l}
       0
     \end{array} \begin{array}{l}
       
     \end{array} \begin{array}{l}
       
     \end{array} \begin{array}{l}
       
     \end{array} \begin{array}{l}
       
     \end{array}\\
     \begin{array}{l}
       
     \end{array} \begin{array}{l}
       
     \end{array} \begin{array}{l}
       
     \end{array} \begin{array}{l}
       
     \end{array} \begin{array}{l}
       0
     \end{array} \begin{array}{l}
       
     \end{array} \begin{array}{l}
       
     \end{array} \begin{array}{l}
       
     \end{array} \begin{array}{l}
       
     \end{array} \begin{array}{l}
       
     \end{array} \begin{array}{l}
       
     \end{array} \begin{array}{l}
       
     \end{array} \begin{array}{l}
       
     \end{array} - l_m (x) + 3 k_m (x)
   \end{array}\right) ] \left(\begin{array}{c}
     - u_a (x)\\
     u_c (x)
   \end{array}\right) = \]
\[ (9 k_m (x)^2 - l_m (x)^2) \left(\begin{array}{c}
     \mathfrak{a}\\
     \mathfrak{c}
   \end{array}\right) . \]
By (5.2),
\[ \left(\begin{array}{c}
     \begin{array}{l}
       - 1
     \end{array} \begin{array}{l}
       
     \end{array} \begin{array}{l}
       
     \end{array} k_m (x)\\
     k_m (x) \begin{array}{l}
       
     \end{array} \begin{array}{l}
       
     \end{array} \begin{array}{l}
       
     \end{array} \begin{array}{l}
       1
     \end{array} \begin{array}{l}
       
     \end{array}
   \end{array}\right) \left(\begin{array}{c}
     \mathfrak{a}\\
     \mathfrak{c}
   \end{array}\right) =\mathfrak{m} \left(\begin{array}{c}
     k_{\mathfrak{c}} (x)\\
     k_{\mathfrak{a}} (x)
   \end{array}\right), \]
or equivalently,
\[ \left(\begin{array}{c}
     \mathfrak{a}\\
     \mathfrak{c}
   \end{array}\right) = \frac{1}{l_m (x)} \left(\begin{array}{c}
     \begin{array}{l}
       
     \end{array} \begin{array}{l}
       - 1
     \end{array} \begin{array}{l}
       
     \end{array} \begin{array}{l}
       
     \end{array} \begin{array}{l}
       
     \end{array} k_m (x)\\
     k_m (x) \begin{array}{l}
       
     \end{array} \begin{array}{l}
       
     \end{array} \begin{array}{l}
       
     \end{array} \begin{array}{l}
       1
     \end{array} \begin{array}{l}
       
     \end{array}
   \end{array}\right) \left(\begin{array}{c}
     k_{\mathfrak{c}} (x)\\
     k_{\mathfrak{a}} (x)
   \end{array}\right) . \]
Substituting this into the identity above, and multiplying the result from the
left by the inverse of the matrix on the right hand side yields,
\[ \left[ 2 l_m (x) \left(\begin{array}{c}
     1 \begin{array}{l}
       
     \end{array} \begin{array}{l}
       
     \end{array} 0\\
     0 \begin{array}{l}
       
     \end{array} \begin{array}{l}
       
     \end{array} 1
   \end{array}\right) + \left(\begin{array}{c}
     - l_m (x) - 3 k_m (x) \begin{array}{l}
       
     \end{array} \begin{array}{l}
       
     \end{array} \begin{array}{l}
       
     \end{array} \begin{array}{l}
       
     \end{array} \begin{array}{l}
       - k_m (x) l_m (x) + 3 k^2_m (x)
     \end{array}\\
     k_m (x) l_m (x) + 3 k^2_m (x) \begin{array}{l}
       
     \end{array} \begin{array}{l}
       
     \end{array} \begin{array}{l}
       
     \end{array} \begin{array}{l}
       
     \end{array} \begin{array}{l}
       - l_m (x) + 3 k_m (x)
     \end{array} \begin{array}{l}
       
     \end{array}
   \end{array}\right) ]  \left(\begin{array}{c}
     - u_a (x)\\
     u_c (x)
   \end{array}\right) = \right. \]
\[ (9 k_m (x)^2 - l_m (x)^2) \left(\begin{array}{c}
     - k_{\mathfrak{c}} (x)\\
     - k_{\mathfrak{a}} (x)
   \end{array}\right), \]
or
\[ \left(\begin{array}{c}
     l_m (x) - 3 k_m (x) \begin{array}{l}
       
     \end{array} \begin{array}{l}
       
     \end{array} \begin{array}{l}
       
     \end{array} \begin{array}{l}
       
     \end{array} \begin{array}{l}
       - k_m (x) l_m (x) + 3 k^2_m (x)
     \end{array}\\
     k_m (x) l_m (x) + 3 k^2_m (x) \begin{array}{l}
       
     \end{array} \begin{array}{l}
       
     \end{array} \begin{array}{l}
       
     \end{array} \begin{array}{l}
       
     \end{array} \begin{array}{l}
       l_m (x) + 3 k_m (x)
     \end{array} \begin{array}{l}
       
     \end{array}
   \end{array}\right) \left(\begin{array}{c}
     u_a (x)\\
     - u_c (x)
   \end{array}\right) = (l_m (x)^2 - 9 k_m (x)^2) \left(\begin{array}{c}
     k_{\mathfrak{c}} (x)\\
     k_{\mathfrak{a}} (x)
   \end{array}\right), \]
which, since the matrix on the left hand side is equal to the product
\[ \left(\begin{array}{c}
     l_m (x) - 3 k_m (x) \begin{array}{l}
       
     \end{array} \begin{array}{l}
       
     \end{array} \begin{array}{l}
       
     \end{array} \begin{array}{l}
       
     \end{array} \begin{array}{l}
       
     \end{array} \begin{array}{l}
       
     \end{array} \begin{array}{l}
       0
     \end{array} \begin{array}{l}
       
     \end{array} \begin{array}{l}
       
     \end{array} \begin{array}{l}
       
     \end{array} \begin{array}{l}
       
     \end{array}\\
     \begin{array}{l}
       
     \end{array} \begin{array}{l}
       
     \end{array} \begin{array}{l}
       
     \end{array} \begin{array}{l}
       
     \end{array} \begin{array}{l}
       0
     \end{array} \begin{array}{l}
       
     \end{array} \begin{array}{l}
       
     \end{array} \begin{array}{l}
       
     \end{array} \begin{array}{l}
       
     \end{array} \begin{array}{l}
       
     \end{array} \begin{array}{l}
       
     \end{array} \begin{array}{l}
       
     \end{array} l_m (x) + 3 k_m (x)
   \end{array}\right) \left(\begin{array}{c}
     \begin{array}{l}
       1
     \end{array} \begin{array}{l}
       
     \end{array} \begin{array}{l}
       
     \end{array} - k_m (x)\\
     k_m (x) \begin{array}{l}
       
     \end{array} \begin{array}{l}
       
     \end{array} \begin{array}{l}
       
     \end{array} \begin{array}{l}
       1
     \end{array} \begin{array}{l}
       
     \end{array}
   \end{array}\right), \]
is equivalent to
\[ \left(\begin{array}{c}
     \begin{array}{l}
       1
     \end{array} \begin{array}{l}
       
     \end{array} \begin{array}{l}
       
     \end{array} - k_m (x)\\
     k_m (x) \begin{array}{l}
       
     \end{array} \begin{array}{l}
       
     \end{array} \begin{array}{l}
       
     \end{array} \begin{array}{l}
       1
     \end{array} \begin{array}{l}
       
     \end{array}
   \end{array}\right) \left(\begin{array}{c}
     u_a (x)\\
     - u_c (x)
   \end{array}\right) = \left(\begin{array}{c}
     l_m (x) + 3 k_m (x) \begin{array}{l}
       
     \end{array} \begin{array}{l}
       
     \end{array} \begin{array}{l}
       
     \end{array} \begin{array}{l}
       
     \end{array} \begin{array}{l}
       
     \end{array} \begin{array}{l}
       
     \end{array} \begin{array}{l}
       0
     \end{array} \begin{array}{l}
       
     \end{array} \begin{array}{l}
       
     \end{array} \begin{array}{l}
       
     \end{array} \begin{array}{l}
       
     \end{array}\\
     \begin{array}{l}
       
     \end{array} \begin{array}{l}
       
     \end{array} \begin{array}{l}
       
     \end{array} \begin{array}{l}
       
     \end{array} \begin{array}{l}
       0
     \end{array} \begin{array}{l}
       
     \end{array} \begin{array}{l}
       
     \end{array} \begin{array}{l}
       
     \end{array} \begin{array}{l}
       
     \end{array} \begin{array}{l}
       
     \end{array} \begin{array}{l}
       
     \end{array} \begin{array}{l}
       
     \end{array} l_m (x) - 3 k_m (x)
   \end{array}\right) \left(\begin{array}{c}
     k_{\mathfrak{c}} (x)\\
     k_{\mathfrak{a}} (x)
   \end{array}\right) . \]
This in turn leads to
\[ (k^2_m (x) + 1) \left(\begin{array}{c}
     u_a (x)\\
     - u_c (x)
   \end{array}\right) = \]
\[ \left(\begin{array}{c}
     \begin{array}{l}
       1
     \end{array} \begin{array}{l}
       
     \end{array} \begin{array}{l}
       
     \end{array} k_m (x)\\
     - k_m (x) \begin{array}{l}
       
     \end{array} \begin{array}{l}
       
     \end{array} \begin{array}{l}
       
     \end{array} \begin{array}{l}
       
     \end{array} \begin{array}{l}
       1
     \end{array} \begin{array}{l}
       
     \end{array}
   \end{array}\right) \left(\begin{array}{c}
     l_m (x) + 3 k_m (x) \begin{array}{l}
       
     \end{array} \begin{array}{l}
       
     \end{array} \begin{array}{l}
       
     \end{array} \begin{array}{l}
       
     \end{array} \begin{array}{l}
       
     \end{array} \begin{array}{l}
       
     \end{array} \begin{array}{l}
       0
     \end{array} \begin{array}{l}
       
     \end{array} \begin{array}{l}
       
     \end{array} \begin{array}{l}
       
     \end{array} \begin{array}{l}
       
     \end{array}\\
     \begin{array}{l}
       
     \end{array} \begin{array}{l}
       
     \end{array} \begin{array}{l}
       
     \end{array} \begin{array}{l}
       
     \end{array} \begin{array}{l}
       0
     \end{array} \begin{array}{l}
       
     \end{array} \begin{array}{l}
       
     \end{array} \begin{array}{l}
       
     \end{array} \begin{array}{l}
       
     \end{array} \begin{array}{l}
       
     \end{array} \begin{array}{l}
       
     \end{array} \begin{array}{l}
       
     \end{array} l_m (x) - 3 k_m (x)
   \end{array}\right) \left(\begin{array}{c}
     k_{\mathfrak{c}} (x)\\
     k_{\mathfrak{a}} (x)
   \end{array}\right) = \]
\[ \  \]
\[ \left(\begin{array}{c}
     l_m (x) + 3 k_m (x) \begin{array}{l}
       
     \end{array} \begin{array}{l}
       
     \end{array} \begin{array}{l}
       
     \end{array} \begin{array}{l}
       
     \end{array} \begin{array}{l}
       k_m (x) l_m (x) - 3 k^2_m (x)
     \end{array}\\
     - k_m (x) l_m (x) - 3 k^2_m (x) \begin{array}{l}
       
     \end{array} \begin{array}{l}
       
     \end{array} \begin{array}{l}
       
     \end{array} \begin{array}{l}
       
     \end{array} \begin{array}{l}
       l_m (x) - 3 k_m (x)
     \end{array} \begin{array}{l}
       
     \end{array} \begin{array}{l}
       
     \end{array} \begin{array}{l}
       
     \end{array}
   \end{array}\right) \left(\begin{array}{c}
     k_{\mathfrak{c}} (x)\\
     k_{\mathfrak{a}} (x)
   \end{array}\right) = \]
\[ \left(\begin{array}{c}
     l_m (x) + 3 k_m (x) \begin{array}{l}
       
     \end{array} \begin{array}{l}
       
     \end{array} \begin{array}{l}
       
     \end{array} \begin{array}{l}
       
     \end{array} \begin{array}{l}
       
     \end{array} \begin{array}{l}
       
     \end{array} \begin{array}{l}
       
     \end{array} \begin{array}{l}
       
     \end{array} \begin{array}{l}
       
     \end{array} \begin{array}{l}
       
     \end{array} k_m (x) l_m (x) + 3 - 3 (k^2_m (x) + 1)\\
     - k_m (x) l_m (x) + 3 - 3 (k^2_m (x) + 1) \begin{array}{l}
       
     \end{array} \begin{array}{l}
       
     \end{array} \begin{array}{l}
       
     \end{array} \begin{array}{l}
       
     \end{array} \begin{array}{l}
       
     \end{array} \begin{array}{l}
       
     \end{array} \begin{array}{l}
       
     \end{array} \begin{array}{l}
       
     \end{array} \begin{array}{l}
       
     \end{array} \begin{array}{l}
       l_m (x) - 3 k_m (x)
     \end{array} \begin{array}{l}
       
     \end{array} \begin{array}{l}
       
     \end{array} \begin{array}{l}
       
     \end{array}
   \end{array}\right) \left(\begin{array}{c}
     k_{\mathfrak{c}} (x)\\
     k_{\mathfrak{a}} (x)
   \end{array}\right) = \]
\[ [ \left(\begin{array}{c}
     l_m (x) + 3 k_m (x) \begin{array}{l}
       
     \end{array} \begin{array}{l}
       
     \end{array} k_m (x) l_m (x) + 3\\
     - k_m (x) l_m (x) + 3 \begin{array}{l}
       
     \end{array} \begin{array}{l}
       l_m (x) - 3 k_m (x)
     \end{array} \begin{array}{l}
       
     \end{array} \begin{array}{l}
       
     \end{array}
   \end{array}\right) - 3 (k^2_m (x) + 1) \left(\begin{array}{c}
     0 \begin{array}{l}
       
     \end{array} \begin{array}{l}
       
     \end{array} 1\\
     1 \begin{array}{l}
       
     \end{array} \begin{array}{l}
       
     \end{array} 0
   \end{array}\right) ]  \left(\begin{array}{c}
     k_{\mathfrak{c}} (x)\\
     k_{\mathfrak{a}} (x)
   \end{array}\right) . \]
This entails,
\[ (k^2_m (x) + 1)  [ \left(\begin{array}{c}
     u_a (x)\\
     - u_c (x)
   \end{array}\right) + 3 \left(\begin{array}{c}
     k_{\mathfrak{a}} (x)\\
     k_{\mathfrak{c}} (x)
   \end{array}\right) ] = \left(\begin{array}{c}
     k_m (x) \begin{array}{l}
       
     \end{array} \begin{array}{l}
       
     \end{array} \begin{array}{l}
       - 1
     \end{array} \begin{array}{l}
       
     \end{array}\\
     1 \begin{array}{l}
       
     \end{array} \begin{array}{l}
       
     \end{array} \begin{array}{l}
       \begin{array}{l}
         
       \end{array}
     \end{array} k_m (x)
   \end{array}\right) \left(\begin{array}{c}
     \begin{array}{l}
       
     \end{array} \begin{array}{l}
       3
     \end{array} \begin{array}{l}
       
     \end{array} \begin{array}{l}
       
     \end{array} \begin{array}{l}
       
     \end{array} l_m (x)\\
     - l_m (x) \begin{array}{l}
       
     \end{array} \begin{array}{l}
       
     \end{array} \begin{array}{l}
       - 3
     \end{array}
   \end{array}\right) \left(\begin{array}{c}
     k_{\mathfrak{c}} (x)\\
     k_{\mathfrak{a}} (x)
   \end{array}\right), \]
and finally, after multiplying this from the left by the inverse of the matrix
$\left(\begin{array}{c}
  k_m (x) \begin{array}{l}
    
  \end{array} \begin{array}{l}
    
  \end{array} \begin{array}{l}
    - 1
  \end{array} \begin{array}{l}
    
  \end{array}\\
  1 \begin{array}{l}
    
  \end{array} \begin{array}{l}
    
  \end{array} \begin{array}{l}
    \begin{array}{l}
      
    \end{array}
  \end{array} k_m (x)
\end{array}\right)$,
\[ \left(\begin{array}{c}
     k_m (x) \begin{array}{l}
       
     \end{array} \begin{array}{l}
       
     \end{array} \begin{array}{l}
       
     \end{array} \begin{array}{l}
       
     \end{array} 1\\
     - 1 \begin{array}{l}
       
     \end{array} \begin{array}{l}
       
     \end{array} \begin{array}{l}
       \begin{array}{l}
         
       \end{array}
     \end{array} k_m (x)
   \end{array}\right) [ \left(\begin{array}{c}
     u_a (x)\\
     - u_c (x)
   \end{array}\right) + 3 \left(\begin{array}{c}
     k_{\mathfrak{a}} (x)\\
     k_{\mathfrak{c}} (x)
   \end{array}\right) ] = \left(\begin{array}{c}
     \begin{array}{l}
       
     \end{array} \begin{array}{l}
       3
     \end{array} \begin{array}{l}
       
     \end{array} \begin{array}{l}
       
     \end{array} \begin{array}{l}
       
     \end{array} \begin{array}{l}
       
     \end{array} l_m (x)\\
     - l_m (x) \begin{array}{l}
       
     \end{array} \begin{array}{l}
       
     \end{array} \begin{array}{l}
       
     \end{array} \begin{array}{l}
       - 3
     \end{array}
   \end{array}\right) \left(\begin{array}{c}
     k_{\mathfrak{c}} (x)\\
     k_{\mathfrak{a}} (x)
   \end{array}\right), \]
which is the claimed identity (5.12) $\Box$

\begin{tabular}{l}
  
\end{tabular}

\

{\tmstrong{5.4 Lemma}} The following identity holds true

(5.14)
\[ \left(\begin{array}{ccc}
     \mathfrak{c} & \mathfrak{m} & \mathfrak{a}\\
     k_c & k_m & k_a\\
     l_c & l_m & l_a
   \end{array}\right) \left(\begin{array}{c}
     \mathfrak{c}\\
     -\mathfrak{b}\\
     \mathfrak{a} 
   \end{array}\right) = \left(\begin{array}{c}
     0 \\
     0\\
     2
   \end{array}\right) \]
\begin{tabular}{l}
  
\end{tabular}

{\tmstrong{Proof}} It follows from (5.2),
\[ c k_m (x) - a = m k_c (x), a k_m (x) + c = m k_a (x) . \]
Squaring these two identities, and adding the results together yields, after
invoking the second identity in (5.3)
\[ (a^2 + c^2) \mathfrak{m}l_m (x) = m^2 (k_a (x)^2 + k_c (x)^2), \]
or equivalently

(5.15)
\[ \mathfrak{b}l_m (x) = k_a (x)^2 + k_c (x)^2 . \]
The following is a consequence of the first and the third identity in (5.3)

(5.16)
\[ k_a (x)^2 + k_c (x)^2 =\mathfrak{a}l_a (x) +\mathfrak{c}l_c (x) - 2. \]
Combining (5.15) and (5.16) yields,

\[ \mathfrak{a} (l_a + 2 k_a x +\mathfrak{a}x^2) -\mathfrak{b} (l_m + 2 k_m x
   +\mathfrak{m}x^2) +\mathfrak{c} (l_c + 2 k_c x +\mathfrak{c}x^2) = 2 . \]
Comparing the coefficients for the quadratic, the linear, and the constant
terms, yields the first, the second, and the third entry, respectively, in the
vector identity (5.14). $\Box$

\

{\tmstrong{Remark}} Since $\mathfrak{a}, \mathfrak{b}, \mathfrak{c}$ are
relatively prime, it follows from the identity (5.14) that

$\det \left(\begin{array}{ccc}
  \mathfrak{c} & \mathfrak{m} & \mathfrak{a}\\
  k_c & k_m & k_a\\
  l_c & l_m & l_a
\end{array}\right) = \pm 2$. This conclusion has already been reached by
Frobenius in [F], Gesam-

melte Abhandlungen, Band III, p.604 (13.).

\begin{tabular}{l}
  
\end{tabular}

{\tmstrong{5.5 Lemma}} The following identity holds true

(5.17)
\[ \mathfrak{a}v_a +\mathfrak{c}v_{\mathfrak{c}} =\mathfrak{a}l_a
   +\mathfrak{c}l_c \]
\begin{tabular}{l}
  
\end{tabular}

{\tmstrong{Proof}} Squaring the two entries of the vector on the left hand
side of (5.13), and adding the results together yields,
\[ (k_m (x) v_a (x) - v_c (x))^2 + (v_a (x) + k_m (x) v_c (x))^2 = (k_m (x)^2
   + 1) (v_a (x)^2 + v_c (x)^2) \]
\[ =\mathfrak{m}l_m (x) (v_a (x)^2 + v_c (x)^2) . \]
Doing the same thing for the vector on the right hand side of (5.13) yields,
\[ (k_a (x) l_m (x) - 3 k_c (x))^2 + (k_c (x) l_m (x) - 3 k_a (x))^2 \]
\[ 9 (k_a (x)^2 + k_c (x)^2) + (k^2_a (x) l_m (x) + k^2_a (x) l_m (x) - 12 k_a
   (x) k_c (x)) l_m (x) . \]
Therefore, by (5.13),
\[ \mathfrak{m}l_m (x) (v_a (x)^2 + v_c (x)^2) = 9 (k_a (x)^2 + k_c (x)^2) +
   (k^2_a (x) l_m (x) + k^2_a (x) l_m (x) - 12 k_a (x) k_c (x)) l_m (x) . \]
Combining this with (5.15) leads to
\[ \mathfrak{m} (v_a (x)^2 + v_c (x)^2) = (9 + l_m (x)^2) b + 12 k_a k_c . \]
Evaluating this identity for the quadratic terms leads to
\[ \mathfrak{m} (2\mathfrak{a}v_a + 4 k^2_a + 2\mathfrak{c}v_c + 4 k^2_c) =
   (2\mathfrak{m}l_m + 4 k^2_m) \mathfrak{b}+ 12\mathfrak{a}\mathfrak{c}, \]
which simplifies to
\[ \mathfrak{a}v_a + 2\mathfrak{a}l_a +\mathfrak{c}v_{\mathfrak{c}} +
   2\mathfrak{c}l_c = 6 + 3\mathfrak{b}l_m . \]
But the last row in (5.14) implies
\[ \mathfrak{a}l_a +\mathfrak{c}l_c = 2 +\mathfrak{b}l_m, \]
and so the claimed identity follows. $\Box$

\begin{tabular}{l}
  
\end{tabular}

The following statement will allow us to make specific choices for the matrix
in (5.14). It has been known for a long time and can be found for instance in
[R], p.163.

\begin{tabular}{l}
  
\end{tabular}

{\tmstrong{5.6 Lemma}} The numbers $k_a, k_m, k_c$ in (5.2) can be chosen so
that they all have the same sign, and

(5.18)
\[ \begin{array}{l}
     | k_a |
   \end{array} \leqslant \frac{\mathfrak{a}}{2}, \begin{array}{l}
     | k_m |
   \end{array} \leqslant \frac{\mathfrak{m}}{2}, \begin{array}{l}
     | k_c |
   \end{array} \leqslant \frac{\mathfrak{c}}{2} . \]

{\tmstrong{Remark}} All the arguments in this section up to Lemma 5.6 are
valid without the requirement that $m = \max \{ a, m, c \}$. Departing from
the arrangement in (5.2) it is through Lemma 5.6, and only through this lemma,
that the maximality of $m$ is being exploited in the present section.

\

The inequalities in (5.18) obviously imply the following,

(5.19)
\[ l_a \leqslant \frac{\mathfrak{a}}{4} + \frac{1}{\mathfrak{a}}, l_m
   \leqslant \frac{\mathfrak{m}}{4} + \frac{1}{\mathfrak{m}}, l_c \leqslant
   \frac{\mathfrak{c}}{4} + \frac{1}{\mathfrak{c}} . \]
Henceforth we shall restrict the parameters in question to those satisfying
(5.18).

\begin{tabular}{l}
  
\end{tabular}

{\tmstrong{5.7 Lemma}} The following identities hold true,

(5.20)
\[ v_a = l_a, v_c = l_c . \]
\begin{tabular}{l}
  
\end{tabular}

{\tmstrong{Proof}} Since, by (5.17),
\[ \mathfrak{a} (v_a - l_a) =\mathfrak{c} (l_c - v_c), \]
and since $\mathfrak{a} \tmop{and} \mathfrak{c}$ are relatively prime, we
conclude that

(5.21)
\[ \frac{v_a - l_a}{\mathfrak{c}} \begin{array}{l}
     \epsilon
   \end{array} \mathbb{Z}, \frac{v_c - l_c}{\mathfrak{a}} \begin{array}{l}
     \epsilon
   \end{array} \mathbb{Z}. \]
By (5.13),
\[ \left(\begin{array}{c}
     v_a\\
     - v_c
   \end{array}\right) = \frac{1}{\mathfrak{m}l_m} \left(\begin{array}{c}
     k_m \begin{array}{l}
       
     \end{array} \begin{array}{l}
       
     \end{array} - 1\\
     1 \begin{array}{l}
       
     \end{array} \begin{array}{l}
       \begin{array}{l}
         
       \end{array}
     \end{array} k_m
   \end{array}\right) \left(\begin{array}{c}
     \begin{array}{l}
       
     \end{array} \begin{array}{l}
       3
     \end{array} \begin{array}{l}
       
     \end{array} \begin{array}{l}
       
     \end{array} l_m\\
     - l_m \begin{array}{l}
       
     \end{array} \begin{array}{l}
       - 3
     \end{array}
   \end{array}\right) \left(\begin{array}{c}
     - k_c\\
     k_a
   \end{array}\right) \]
\[ = \frac{1}{\mathfrak{m}l_m} \left(\begin{array}{c}
     l_m + 3 k_m \begin{array}{l}
       
     \end{array} \begin{array}{l}
       
     \end{array} k_m l_m + 3\\
     - k_m l_m + 3 \begin{array}{l}
       
     \end{array} \begin{array}{l}
       l_m - 3 k_m
     \end{array} \begin{array}{l}
       
     \end{array} \begin{array}{l}
       
     \end{array}
   \end{array}\right) \left(\begin{array}{c}
     - k_c\\
     k_a
   \end{array}\right), \]
in particular

(5.22)
\[ v_a = \frac{1}{\mathfrak{m}l_m} (- (l_m + 3 k_m) k_c + (k_m l_m + 3) k_a \]
Suppose that $\mathfrak{a} \leqslant \mathfrak{c} \leqslant \mathfrak{m}$. We
are now going to use \ (5.18), (5.19) and (5.22) to obtain an upper bound for
$| \frac{v_a - l_a}{\mathfrak{c}} |$. First,
\[ | \frac{l_m + 3 k_m}{\mathfrak{m}l_m} | \leqslant \frac{1}{\mathfrak{m}} +
   3 \frac{1}{\sqrt{\mathfrak{m}l_m}}, \]
\[ | \frac{k_m l_m + 3}{\mathfrak{m}l_m} | \leqslant \frac{1}{2} +
   \frac{3}{\mathfrak{m}l_m}, \]
and therefore, by (5.22) and (5.18)
\[ \begin{array}{l}
     | v_a |
   \end{array} \leqslant (\frac{1}{\mathfrak{m}} + 3
   \frac{1}{\sqrt{\mathfrak{m}l_m}}) \frac{\mathfrak{c}}{2} + (\frac{1}{2} +
   \frac{3}{\mathfrak{m}l_m}) \frac{\mathfrak{a}}{2} \leqslant \frac{1}{2}
   (\frac{1}{\mathfrak{m}} + 3 \frac{1}{\sqrt{\mathfrak{m}l_m}} + \frac{1}{2}
   + \frac{3}{\mathfrak{m}l_m}) \mathfrak{c}. \]
Hence, by (5.19),
\[ | \frac{v_a - l_a}{\mathfrak{c}} | \leqslant \frac{| v_a |}{\mathfrak{c}} +
   \frac{l_a}{\mathfrak{c}} \leqslant \frac{1}{2} (\frac{1}{\mathfrak{m}} + 3
   \frac{1}{\sqrt{\mathfrak{m}l_m}} + \frac{1}{2} + \frac{3}{\mathfrak{m}l_m})
   + \frac{1}{4} + \frac{1}{\mathfrak{a}\mathfrak{c}} . \]
If $\mathfrak{m} \geqslant 29$, then $\mathfrak{a}\mathfrak{c} \geqslant 10$,
and in this case the right hand side of this inequality is a number which is
less than 1. For those Markoff triples which meet this condition, it follows
from (5.21) that $v_a$ is equal to $l_a$, and hence by (5.17), $v_c$ is equal
to $l_c$, settling the claim of the lemma in case $\mathfrak{a} \leqslant
\mathfrak{c} \leqslant \mathfrak{m}$. If $\mathfrak{c} \leqslant \mathfrak{a}
\leqslant \mathfrak{m}$, then the same type of estimates for $| \frac{v_c -
l_c}{\mathfrak{a}} |$ in place of $| \frac{v_a - l_a}{\mathfrak{c}} |$ lead to
the same conclusion. For the Markoff triples for which the largest member is
less than 29, namely the triples (1,1,1); (1,1,2); (1,2,5) and (1,5,13), the
validity of the claim can be checked through inspection. $\Box$

\begin{tabular}{l}
  
\end{tabular}

{\tmstrong{Proof of Proposition 5.1}} The first two identities follow from
Lemma 5.3 and Lemma 5.7. In order to establish the third identity we
reintroduce the parameter $x$ into the first identity, \ write out the result
in terms of $x$,
\[ (k_m +\mathfrak{m}x) (l_a + 2 k_a x +\mathfrak{a}x^2) - (k_a
   +\mathfrak{a}x) (l_m + 2 k_m x +\mathfrak{m}x^2) = l_c + 2 k_c
   +\mathfrak{c}x^2 + 3 k_c + 3\mathfrak{c}, \]
and compare the coefficients of the linear terms. This yields the third
identity. The fourth identity can be shown in exactly the same way by
employing the second identity. $\Box$

\begin{tabular}{l}
  
\end{tabular}

Putting together what has been established so far, we can summarize the
situation through the following (incomplete) matrix identities,

(5.23)
\[ \left(\begin{array}{ccc}
     \mathfrak{c} & \mathfrak{m} & \mathfrak{a}\\
     k_c & k_m & k_a\\
     l_c & l_m & l_a
   \end{array}\right) \begin{array}{l}
     - 1\\
     \\
     
   \end{array} = \frac{1}{2} \left(\begin{array}{ccc}
     l_c + 3 k_c & - (2 k_c + 3\mathfrak{c}) & \mathfrak{c}\\
     ? & ? & -\mathfrak{b}\\
     l_a - 3 k_a & - (2 k_a - 3\mathfrak{a}) & \mathfrak{a}
   \end{array}\right), \]
(5.24)
\[ Z^{- 1} \mathcal{A}= \frac{1}{2} \left(\begin{array}{ccc}
     l_c + 3 k_c & - (2 k_c + 3\mathfrak{c}) & \mathfrak{c}\\
     ? & ? & -\mathfrak{b}\\
     l_a - 3 k_a & - (2 k_a - 3\mathfrak{a}) & \mathfrak{a}
   \end{array}\right) \left(\begin{array}{ccc}
     \frac{1}{3} & 0 & 0\\
     0 & 1 & 0\\
     0 & 0 & 3
   \end{array}\right) . \]
The factor $\frac{1}{2}$ on the right hand side of (5.23) is a consequence of
(5.14). What follows are comments on the seven enunciated entries of the
matrix on the right hand side of (5.23). If we consider the $x$-parameter
version of the first identity in (5.4), namely $k_m (x) l_a (x) - k_a (x) l_m
(x) = l_c (x) + 3 k_c (x)$, then comparing the coefficients of the constant
terms yields the entry (1,1), comparing the coefficients of the linear terms
yields the entry (1,2), and comparing the coefficients of the quadratic terms
yields the entry (1,3) of the matrix on the right hand side in (5.23).
Likewise, if we consider the $x$-parameter version of the second identity in
(5.4), namely $k_c (x) l_m (x) - k_m (x) l_c (x) = l_a (x) - 3 k_a (x)$, then
comparing the coefficients of the constant terms yields the entry (3,1),
comparing the coefficients of the linear terms yields the entry (3,2), and
comparing the coefficients of the quadratic terms yields the entry (3,3) of
the matrix on the right hand side in (5.23). The entry (2,3) is a consequence
of (5.14). Turning to the matrix identity (5.24), it suffices to note that
this is a consequence of (5.23) and Lemma 5.7.

Our next objective is to obtain more information about the entries (2,1) and
(2,2) of the matrix on the right hand side in (5.23). Let

(5.25)
\[ \mathfrak{A} (c, m, a) = \left(\begin{array}{ccc}
     \mathfrak{c} & \mathfrak{m} & \mathfrak{a}\\
     k_c & k_m & k_a\\
     l_c & l_m & l_a
   \end{array}\right), \]
\[ \mathfrak{B} (c, m, a) = \frac{1}{2 m^2} \left(\begin{array}{ccc}
     - 2 (c + a k_m) + m c (\frac{l_m}{3} + k_m) \begin{array}{l}
       
     \end{array} & 2 m (a - c k_m) - m^2 c \begin{array}{l}
       
     \end{array} & m^2 c\\
     2 a c - (c^2 - a^2) k_m - m b \frac{l_m}{3} \begin{array}{l}
       
     \end{array} & m (c^2 - a^2 + 2 b k_m) \begin{array}{l}
       
     \end{array} & - m^2 b\\
     2 (c k_m - a) + m a (\frac{l_m}{3} - k_m) \begin{array}{l}
       
     \end{array} & m^2 a - 2 m (c + a k_m) \begin{array}{l}
       
     \end{array} & m^2 a
   \end{array}\right) . \]
Note that, due to the specification of the parameters in (5.18), and since $b
= a c - m,$both of these matrices are uniquely determined by the Markoff
triple $(a, m, c)$ up to the common sign chosen for $k_a, k_m, k_c$ (see Lemma
5.6). Now, instead of $(a, m, c)$ we consider the Markoff triple $(a, b, c)$
in this context. Since $m \geqslant \max (a, c)$ and $m = \tmop{ac} - b$, we
must have $b \leqslant$max$(a, c)$. This means that max$(a, b, c) \epsilon \{
a, c \}$. Suppose max$(a, b, c) = c$. Then, considering the matrices
$\mathfrak{A} \tmop{and} \mathfrak{B}$ in this context, we have two distinct
choices to arrange the members of the triple $(a, b, c)$, so that the
resulting situation is consistent with our settings for $(a, m, c)$, namely
\[ (a, c, b) \tmop{or} (b, c, a) . \]
In the first case (5.23) and (5.24) turn into, respectively,
\[ \mathfrak{A} (a, c, b)^{- 1} = \frac{1}{2} \left(\begin{array}{ccc}
     l_a + 3 k_a & - (2 k_a + 3\mathfrak{a}) & \mathfrak{a}\\
     ? & ? & - (3\mathfrak{a}\mathfrak{b}-\mathfrak{c})\\
     l_{\mathfrak{b}} - 3 k_b & - (2 k_b - 3\mathfrak{b}) & \mathfrak{b}
   \end{array}\right), \]
\[ \mathfrak{B} (a, c, b) = \frac{1}{2} \left(\begin{array}{ccc}
     l_a + 3 k_a & - (2 k_a + 3\mathfrak{a}) & \mathfrak{a}\\
     ? & ? & - (3\mathfrak{a}\mathfrak{b}-\mathfrak{c})\\
     l_{\mathfrak{b}} - 3 k_b & - (2 k_b - 3\mathfrak{b}) & \mathfrak{b}
   \end{array}\right) \left(\begin{array}{ccc}
     \frac{1}{3} & 0 & 0\\
     0 & 1 & 0\\
     0 & 0 & 3
   \end{array}\right), \]
while we get in the second case,
\[ \mathfrak{A} (b, c, a)^{- 1} = \frac{1}{2} \left(\begin{array}{ccc}
     l_b + 3 k_b & - (2 k_b + 3\mathfrak{b}) & \mathfrak{b}\\
     ? & ? & - (3\mathfrak{a}\mathfrak{b}-\mathfrak{c})\\
     l_a - 3 k_a & - (2 k_a - 3\mathfrak{a}) & \mathfrak{a}
   \end{array}\right) \]
\[ \mathfrak{B} (b, c, a) = \frac{1}{2} \left(\begin{array}{ccc}
     l_b + 3 k_b & - (2 k_b + 3\mathfrak{b}) & \mathfrak{b}\\
     ? & ? & - (3\mathfrak{a}\mathfrak{b}-\mathfrak{c})\\
     l_a - 3 k_a & - (2 k_a - 3\mathfrak{a}) & \mathfrak{a}
   \end{array}\right) \left(\begin{array}{ccc}
     \frac{1}{3} & 0 & 0\\
     0 & 1 & 0\\
     0 & 0 & 3
   \end{array}\right) . \]
Either case leads to the following further specification of (5.23),

(5.26)
\[ \left(\begin{array}{ccc}
     \mathfrak{c} & \mathfrak{m} & \mathfrak{a}\\
     k_c & k_m & k_a\\
     l_c & l_m & l_a
   \end{array}\right) \begin{array}{l}
     - 1\\
     \\
     
   \end{array} = \frac{1}{2} \left(\begin{array}{ccc}
     l_c + 3 k_c & - (2 k_c + 3\mathfrak{c}) & \mathfrak{c}\\
     - (l_b + 3 \nu k_b) & 2 k_b + 3 \nu \mathfrak{b} & -\mathfrak{b}\\
     l_a - 3 k_a & - (2 k_a - 3\mathfrak{a}) & \mathfrak{a}
   \end{array}\right) ; \begin{array}{l}
     
   \end{array} \nu \epsilon \{ - 1, 1 \} . \]

Before proceeding to give a more specific determination of the matrix $Z^{- 1}
\mathcal{A}$ which involves a certain quadratic equation, we need to look at
that equation first.

\begin{tabular}{l}
  
\end{tabular}

{\tmstrong{5.7 Lemma}} The quadratic equation

(5.27)
\[ \mathfrak{m}y^2 - 4\mathfrak{b}k_m y - (9\mathfrak{b}^2 - 4) \mathfrak{m}+
   4\mathfrak{b}^2 l_m - 8\mathfrak{b} \sigma = 0, \]
has always two rational solutions in case $\sigma = 1$, and it has no rational
solutions in case $\sigma = - 1$.

\begin{tabular}{l}
  
\end{tabular}

{\tmstrong{Proof}} First we show that the discriminant $D$ of this equation is
a perfect square in case $\sigma = 1$.
\[ D = 16\mathfrak{b}^2 k^2_m + 4\mathfrak{m}^2 (9\mathfrak{b}^2 - 4) -
   16\mathfrak{b}^2 \mathfrak{m}l_m - 8\mathfrak{b}\mathfrak{m}= -
   16\mathfrak{b}^2 + 4\mathfrak{m}^2 (9\mathfrak{b}^2 - 4) -
   8\mathfrak{b}\mathfrak{m}  \]
\[ \frac{D}{4} = 9\mathfrak{m}^2 \mathfrak{b}^2 - 4 (\mathfrak{m}^2
   {+\mathfrak{b}^2} ) - 8\mathfrak{b}\mathfrak{m}  \frac{}{} = 9
   (\mathfrak{a}^2 +\mathfrak{c}^2)^2 - 4
   (3\mathfrak{a}\mathfrak{c}-\mathfrak{b})^2 - 4\mathfrak{b}^2 - 8
   (\mathfrak{a}^2 +\mathfrak{c}^2) \]
\[ = 9 (\mathfrak{a}^4 + 2\mathfrak{a}^2 \mathfrak{c}^2 +\mathfrak{c}^4) - 4
   (9\mathfrak{a}^2 \mathfrak{c}^2 -
   6\mathfrak{a}\mathfrak{b}\mathfrak{c}+\mathfrak{b}^2) - 4\mathfrak{b}^2 - 8
   (\mathfrak{a}^2 +\mathfrak{c}^2) \]
\[ = 9\mathfrak{a}^4 + 18\mathfrak{a}^2 \mathfrak{c}^2 + 9\mathfrak{c}^4 -
   36\mathfrak{a}^2 \mathfrak{c}^2 + 8 (\mathfrak{a}^2 +\mathfrak{b}^2
   +\mathfrak{c}^2) - 8\mathfrak{b}^2 - 8 (\mathfrak{a}^2 +\mathfrak{c}^2) \]
\[ = 9\mathfrak{a}^4 - 18\mathfrak{a}^2 \mathfrak{c}^2 + 9\mathfrak{c}^4 = 9
   (\mathfrak{a}^2 -\mathfrak{c}^2)^2 . \]
It follows that (5.27) has two rational solutions in case $\sigma = 1$. In
order to show that (5.27) does not have a rational solution in case $\sigma =
- 1$, we are going to show that its discriminant $D +
64\mathfrak{m}\mathfrak{b}$ is not a perfect square, or rather that
\[ \frac{D}{4} + 16\mathfrak{m}\mathfrak{b}= 9 (\mathfrak{a}+\mathfrak{c})^2
   (\mathfrak{a}-\mathfrak{c})^2 + 16\mathfrak{m}\mathfrak{b}= 9
   (\mathfrak{m}\mathfrak{b}+ 2\mathfrak{a}\mathfrak{c})
   (\mathfrak{m}\mathfrak{b}- 2\mathfrak{a}\mathfrak{c}) +
   16\mathfrak{m}\mathfrak{b} \]
\[ = 9 (\mathfrak{m}^2 \mathfrak{b}^2 - 4\mathfrak{a}^2 \mathfrak{c}^2) +
   16\mathfrak{m}\mathfrak{b}= 9\mathfrak{m}^2 \mathfrak{b}^2 - 4
   (\mathfrak{m}+\mathfrak{b})^2 + 16\mathfrak{m}\mathfrak{b} \]
\[ = 9\mathfrak{m}^2 \mathfrak{b}^2 - 4\mathfrak{m}^2 -
   8\mathfrak{m}\mathfrak{b}- 4\mathfrak{b}^2 + 16\mathfrak{m}\mathfrak{b}=
   9\mathfrak{m}^2 \mathfrak{b}^2 - 4 (\mathfrak{m}-\mathfrak{b})^2 \]
is not a perfect square. Suppose this were false, which means that there
exists an integer $w$ such that,
\[ 4 (\mathfrak{m}-\mathfrak{b})^2 + w^2 = 9\mathfrak{m}^2 \mathfrak{b}^2 . \]
Since $\mathfrak{m}-\mathfrak{b}$ and $\mathfrak{m}^2 \mathfrak{b}^2$ are
relatively prime (because $\mathfrak{m}$ and $\mathfrak{b}$ are relatively
prime), and since $\mathfrak{m}-\mathfrak{b}$ is not divisible by 3 (because
$\mathfrak{m}+\mathfrak{b}$ is divisible by 3, but neither $\mathfrak{m}$ nor
$\mathfrak{b}$ are divisible by 3), by the standard parametrization of
primitive Pythagorean triples, there exist integers $u \tmop{and} v$ such
that,
\[ \mathfrak{m}-\mathfrak{b}= u v, 3\mathfrak{m}\mathfrak{b}= u^2 + v^2 . \]
This implies that $u^2 + v^2$ is divisible by 3. However, since $u^2 + v^2$
divided be the square of the greatest common divisor of $u \tmop{and} v$
cannot be divisible by 3, due to the fact that such a number can only have odd
prime factors which are equal to 1 modulo 4, both $u \tmop{and} v$ have to be
divisible by 3. this in turn implies that $\mathfrak{m}\mathfrak{b}$ has to be
divisible by 3, which means that $\mathfrak{m} \tmop{or} \mathfrak{b}$ is
divisible by 3. This not the case, because $\mathfrak{m} \tmop{and}
\mathfrak{b}$ are Markoff numbers. That contradiction settles our claim.
$\Box$

\begin{tabular}{l}
  
\end{tabular}

{\tmstrong{5.8 Proposition}} The following identity holds true,

(5.28)
\[ Z^{- 1} \mathcal{A}= \frac{1}{2} \left(\begin{array}{ccc}
     l_c + 3 k_c & - (2 k_c + 3\mathfrak{c}) & \mathfrak{c}\\
     - (l_b + 3 \nu k_b) & 2 k_b + 3 \nu \mathfrak{b} & -\mathfrak{b}\\
     l_a - 3 k_a & - (2 k_a - 3\mathfrak{a}) & \mathfrak{a}
   \end{array}\right) \left(\begin{array}{ccc}
     \frac{1}{3} & 0 & 0\\
     0 & 1 & 0\\
     0 & 0 & 3
   \end{array}\right), \]
where

(5.29)
\[ \nu = \{ \begin{array}{l}
     1\\
     - 1
   \end{array} \begin{array}{l}
     
   \end{array} \begin{array}{l}
     \tmop{if} \mathfrak{a}<\mathfrak{c}\\
     \tmop{if} \mathfrak{a}>\mathfrak{c}
   \end{array} \]
Moreover,

(5.30)
\[ \det \left(\begin{array}{ccc}
     \mathfrak{c} & \mathfrak{m} & \mathfrak{a}\\
     k_c & k_m & k_a\\
     l_c & l_m & l_a
   \end{array}\right) = 2 \]
\begin{tabular}{l}
  
\end{tabular}

{\tmstrong{Proof}} Note that in each row of the matrix
$\left(\begin{array}{ccc}
  l_c + 3 k_c & - (2 k_c + 3\mathfrak{c}) & \mathfrak{c}\\
  - (l_b + 3 \nu k_b) & 2 k_b + 3 \nu \mathfrak{b} & -\mathfrak{b}\\
  l_a - 3 k_a & - (2 k_a - 3\mathfrak{a}) & \mathfrak{a}
\end{array}\right)$ the entries represent the coefficients of a binary
quadratic form which is equivalent to a (reduced) Markoff form. For
convenience we shall henceforth address the discriminant of a quadratic form
whose coefficients agree with the entries of a row vector as the discriminant
of that row. Note that inverting the entries in such a row vector leads to the
same discriminant. For instance the second row in the matrix
$\left(\begin{array}{ccc}
  l_c + 3 k_c & - (2 k_c + 3\mathfrak{c}) & \mathfrak{c}\\
  - (l_b + 3 \sigma k_b) & 2 k_b + 3 \sigma \mathfrak{b} & -\mathfrak{b}\\
  l_a - 3 k_a & - (2 k_a - 3\mathfrak{a}) & \mathfrak{a}
\end{array}\right)$ has the discriminant $9\mathfrak{b}^2 - 4$. We are now
going to show that the second row vector in the matrix
\[ Z^{- 1} \mathcal{A} \left(\begin{array}{ccc}
     3 & 0 & 0\\
     0 & 1 & 0\\
     0 & 0 & \frac{1}{3}
   \end{array}\right) \]
has the discriminant $9\mathfrak{b}^2 - 4$, or equivalently, that the row
vector

(5.31)
\[ \frac{1}{\mathfrak{m}} \left( 2\mathfrak{a}\mathfrak{c}- (\mathfrak{c}^2
   -\mathfrak{a}^2) k_m -\mathfrak{m}\mathfrak{b} \frac{l_m}{3},
   \begin{array}{l}
     
   \end{array} \mathfrak{m}(3 (\mathfrak{c}^2 -\mathfrak{a}^2) +
   2\mathfrak{b}k_m), \begin{array}{l}
     
   \end{array} - 3\mathfrak{m}^2 \mathfrak{b} \right) \]
has the discriminant $\mathfrak{m}^2 (9\mathfrak{b}^2 - 4)$,
\[ \frac{1}{\mathfrak{m}^2} [(\mathfrak{m}(3 (\mathfrak{c}^2 -\mathfrak{a}^2)
   + 2\mathfrak{b}k_m))^2 - 4 (2\mathfrak{a}\mathfrak{c}- (\mathfrak{c}^2
   -\mathfrak{a}^2) k_m -\mathfrak{m}\mathfrak{b} \frac{l_m}{3}) (-
   3\mathfrak{m}^2 \mathfrak{b})] = \]
\[ 3 (\mathfrak{c}^2 -\mathfrak{a}^2) + 2\mathfrak{b}k_m)^2 + 4\mathfrak{b}
   (6\mathfrak{a}\mathfrak{c}- 3 (\mathfrak{c}^2 -\mathfrak{a}^2) k_m
   -\mathfrak{m}\mathfrak{b}l_m) = \]
\[ 9 (\mathfrak{c}^2 -\mathfrak{a}^2)^2 + 12 (\mathfrak{c}^2 -\mathfrak{a}^2)
   \mathfrak{b}k_m + 4\mathfrak{b}^2 k^2_m +
   24\mathfrak{a}\mathfrak{b}\mathfrak{c}- 12\mathfrak{b} (\mathfrak{c}^2
   -\mathfrak{a}^2) k_m - 4\mathfrak{b}^2 \mathfrak{m}l_m = \]
\[ 9 (\mathfrak{c}^2 -\mathfrak{a}^2)^2 - 4\mathfrak{b}^2 +
   24\mathfrak{a}\mathfrak{b}\mathfrak{c}= 9 (\mathfrak{c}+\mathfrak{a})^2
   (\mathfrak{c}-\mathfrak{a})^2 - 4\mathfrak{b}^2 +
   24\mathfrak{a}\mathfrak{b}\mathfrak{c}= \]
\[ 9 (\mathfrak{m}\mathfrak{b}+ 2\mathfrak{a}\mathfrak{c})
   (\mathfrak{m}\mathfrak{b}- 2\mathfrak{a}\mathfrak{c}) - 4\mathfrak{b}^2 +
   24\mathfrak{a}\mathfrak{b}\mathfrak{c}= 9 (\mathfrak{m}^2 \mathfrak{b}^2 -
   4\mathfrak{a}^2 \mathfrak{c}^2) - 4\mathfrak{b}^2 +
   24\mathfrak{a}\mathfrak{b}\mathfrak{c}= \]
\[ 9\mathfrak{m}^2 \mathfrak{b}^2 - 12\mathfrak{a}\mathfrak{c}
   (3\mathfrak{a}\mathfrak{c}-\mathfrak{b}) - 4\mathfrak{b}^2 +
   12\mathfrak{a}\mathfrak{b}\mathfrak{c}= 9\mathfrak{m}^2 \mathfrak{b}^2 -
   12\mathfrak{a}\mathfrak{c}\mathfrak{m}- 4\mathfrak{b}^2 + 4 (\mathfrak{a}^2
   +\mathfrak{b}^2 +\mathfrak{c}^2) = \]
\[ 9\mathfrak{m}^2 \mathfrak{b}^2 - 12\mathfrak{a}\mathfrak{c}\mathfrak{m}+
   4\mathfrak{b}\mathfrak{m}= 9\mathfrak{m}^2 \mathfrak{b}^2 - 4\mathfrak{m}
   (3\mathfrak{a}\mathfrak{c}-\mathfrak{b}) =\mathfrak{m}^2 (9\mathfrak{b}^2 -
   4) . \]
Let $x$ be the (2,1) entry, and let $y$ be the (2,2) entry in the matrix
$\left(\begin{array}{ccc}
  l_c + 3 k_c & - (2 k_c + 3\mathfrak{c}) & \mathfrak{c}\\
  ? & ? & -\mathfrak{b}\\
  l_a - 3 k_a & - (2 k_a - 3\mathfrak{a}) & \mathfrak{a}
\end{array}\right)$ in (5.24). Since det$(\mathcal{A}^{- 1} Z) = 2$, it
follows from (5.23) and (5.24), as well as the determination of the
discriminant of the vector in (5.31), that $x$ and $y$ solve the following two
diophantine equations,
\[ \mathfrak{m}x + k_m y -\mathfrak{b}l_m = 2, \begin{array}{l}
     
   \end{array} y^2 + 4\mathfrak{b}x = 9\mathfrak{b}^2 - 4 . \]
This in turn leads to the quadratic equation (5.27) for the case $\sigma = 1$.
The same line of

reasoning applied to the matrix $\left(\begin{array}{ccc}
  l_c + 3 k_c & - (2 k_c + 3\mathfrak{c}) & \mathfrak{c}\\
  - (l_b + 3 \nu k_b) & 2 k_b + 3 \nu \mathfrak{b} & -\mathfrak{b}\\
  l_a - 3 k_a & - (2 k_a - 3\mathfrak{a}) & \mathfrak{a}
\end{array}\right)$ in (5.26) leads us to a similar conclusion, namely that
the second entry in the second row of this matrix solves the quadratic
equation (5.27) for $\sigma = 1 \tmop{or} \sigma = - 1$. But since Lemma 5.7
states that there are no rational solutions to (5.27) in case $\sigma = - 1$,
it follows once again that $\sigma = 1$, which settles (5.30). Comparison of
the outcome of these two lines of reasoning leads to the conclusion that
(5.28) holds true for some $\nu \epsilon \begin{array}{l}
  \{ - 1, 1 \} \nobracket
\end{array}$. In order to establish (5.29) we observe that (5.28) evaluated
for the entry (2,2) of that matrix yields the following identity,
\[ 3 (\mathfrak{c}^2 -\mathfrak{a}^2) + 2\mathfrak{b}k_m =\mathfrak{m} (2
   k_{\mathfrak{b}} + 3 \nu \mathfrak{b}) = 2\mathfrak{m}k_b + 3 \nu
   (\mathfrak{c}^2 +\mathfrak{a}^2) . \]
This in turn leads to,
\[ \mathfrak{b}k_m -\mathfrak{m}k_{\mathfrak{b}} = \{ \begin{array}{l}
     3\mathfrak{a}^2\\
     - 3\mathfrak{c}^2
   \end{array} \begin{array}{l}
     
   \end{array} \begin{array}{l}
     \tmop{if} \nu = 1\\
     \tmop{if} \nu = - 1
   \end{array} . \]
But since
\[ \begin{array}{l}
     | \mathfrak{b}k_m -\mathfrak{m}k_{\mathfrak{b}} |
   \end{array} \leqslant \begin{array}{l}
     \mathfrak{b} | k_m | +\mathfrak{m} | k_{\mathfrak{b}} |
   \end{array} \leqslant \frac{\mathfrak{m}\mathfrak{b}}{2} +
   \frac{\mathfrak{m}\mathfrak{b}}{2} =\mathfrak{m}\mathfrak{b}=\mathfrak{c}^2
   +\mathfrak{a}^2 \leqslant 2 (\max (\mathfrak{a}, \mathfrak{c}))^2 < 3 (\max
   (\mathfrak{a}, \mathfrak{c}))^2, \]
(5.29) follows. $\Box$

\begin{tabular}{l}
  
\end{tabular}

Returning to the settings in the third remark following Proposition 1.2 in
Section 1 we can now give a conclusive description of the parameter $\nu$ in
(5.26) in relation to the tree of Markoff triples. First in (5.18) all the
parameters are positive (cf. [Zh1], Lemma 2). If $(A, A B, B)$ is an
admissible triple of 2x2 matrices, then the matrix $N$ constructed in
Proposition 1.2 such that
\[ N^t M (3, 3, 3) N = M (\tmop{tr} (A), \tmop{tr} (\tmop{AB}), \tmop{tr}
   (B)), \]
has the property
\[ N = \frac{1}{2} \left(\begin{array}{ccc}
     1 & - 3 & 1\\
     1 & 1 & - 1\\
     - 1 & 1 & 1
   \end{array}\right) \mathfrak{A} (\tmop{tr} (A), \tmop{tr} (A B), \tmop{tr}
   (B)) . \]
Moreover,
\[ \nu = - 1 \tmop{for} \mathfrak{A} (\tmop{tr} (A), \tmop{tr} (A^2 B),
   \tmop{tr} (A B))^{- 1} \tmop{and} \nu = 1 \tmop{for} \mathfrak{A}
   (\tmop{tr} (A B), \tmop{tr} (A  B^2), \tmop{tr} (B))^{- 1} . \]
In other words, replacing $A$ results in a positive value for $\nu$, while
replacing $B$ results in a negative value for $\nu$. Since $\tmop{tr}
(\tmop{AB}) \geqslant \max (\tmop{tr} (A), \tmop{tr} (B))$ we can finally
conclude that (5.29) implies a complete determination of the matrix $Z^{- 1}
\mathcal{A}$ in terms of Markoff triples and their affiliated quadratic
residues subject to the specification (5.18).

\begin{tabular}{l}
  
\end{tabular}

{\tmstrong{5.9 Corollary}} The following identity holds true,
\[ Z^{- 1} \mathcal{A}= \left(\begin{array}{ccc}
     \mathfrak{c} & \mathfrak{m} & \mathfrak{a}\\
     k_c & k_m & k_a\\
     l_c & l_m & l_a
   \end{array}\right) \begin{array}{l}
     - 1\\
     \\
     
   \end{array} \left(\begin{array}{ccc}
     \frac{1}{3} & 0 & 0\\
     0 & 1 & 0\\
     0 & 0 & 3
   \end{array}\right) . \]
\begin{tabular}{l}
  
\end{tabular}

\begin{tabular}{l}
  
\end{tabular}

{\tmstrong{{\tmem{6 Some number theoretic conclusions}}}}

\begin{tabular}{l}
  
\end{tabular}

The point of departure in the present section is the observation that,
disregarding the factor $\frac{1}{2}$, on the one hand, the entries in the
rows of the matrix on the right hand side of the identity (5.26) are the
coefficients of indefinite binary quadratic forms which are equivalent to
Markoff forms associated with the corresponding Markoff numbers in the last
column, and which in the case of the last row corresponds to a reduced form,
i.e. it actually is a Markoff form. On the other hand, the entries of the
columns of the matrix on the left hand side are representations of the number
1 by the ternary quadratic form

(6.1)
\[ Q (x, y, z) = x z - y^2 \]
All elements in the group of automorphs of this form are given by the
matrices

(6.2)
\[ \left(\begin{array}{ccc}
     p^2 & 2 p q & q^2\\
     p r \begin{array}{l}
       
     \end{array} & p s + q r \begin{array}{l}
       
     \end{array} & q s\\
     r^2 & 2 r s & s^2
   \end{array}\right), \tmop{where} \left(\begin{array}{c}
     \begin{array}{l}
       p\\
       r
     \end{array} \begin{array}{l}
       q\\
       s
     \end{array}
   \end{array}\right) \begin{array}{l}
     \epsilon
   \end{array} \tmop{SL} (2, \mathbb{Z}) . \]
This observation goes all the way back to Gauss (cf. [Ba], Kapitel I, pp.
22-23). If $p, q, r, s$ are elements in an arbitrary commutative ring, then we
always have the formula
\[ \det \left(\begin{array}{ccc}
     p^2 & 2 p q & q^2\\
     p r \begin{array}{l}
       
     \end{array} & p s + q r \begin{array}{l}
       
     \end{array} & q s\\
     r^2 & 2 r s & s^2
   \end{array}\right) = (\det \left(\begin{array}{c}
     \begin{array}{l}
       p\\
       r
     \end{array} \begin{array}{l}
       q\\
       s
     \end{array}
   \end{array}\right))^3 . \]
For a given matrix $A \begin{array}{l}
  \epsilon
\end{array} \tmop{SL} (2, \mathbb{Z})$ let $\Psi (A)$ be the corresponding 3x3
matrix in (6.2). Then $\Psi$ determines an isomorphism from the group
PSL(2,$\mathbb{Z}$) onto the group of automorphs of the form (6.1) with
determinant 1. For any integral solution of the equation

(6.3)
\[ x z - y^2 = 1, \]
we define the number $| y |$ as the height of the triple $(x, y, z)$.
Implementing a procedure akin to the continued fraction algorithm, it is an
elementary task to show that, by employing a finite sequence of matrices of
the form
\[ \Psi (\left(\begin{array}{c}
     \begin{array}{l}
       1\\
       r
     \end{array} \begin{array}{l}
       0\\
       1
     \end{array}
   \end{array}\right)) \tmop{or} \Psi (\left(\begin{array}{c}
     \begin{array}{l}
       1\\
       0
     \end{array} \begin{array}{l}
       q\\
       1
     \end{array}
   \end{array}\right)) \]
to a triple $(x, y, z)$ solving (6.3), one can reduce the height of such a
triple to the smallest possible value, which is 0. In the sequel we shall need
the following by-product of this procedure.

\begin{tabular}{l}
  
\end{tabular}

{\tmstrong{6.1 Lemma}} If $(x, y, z)$ is a solution of (6.3), then the
application of a matrix of the form

$\Psi (\left(\begin{array}{c}
  \begin{array}{l}
    1\\
    r
  \end{array} \begin{array}{l}
    0\\
    1
  \end{array}
\end{array}\right)) \tmop{or} \Psi (\left(\begin{array}{c}
  \begin{array}{l}
    1\\
    0
  \end{array} \begin{array}{l}
    q\\
    1
  \end{array}
\end{array}\right))$ to the vector $(x, y, z)^t$ does not change the sign of
$x \tmop{and} z$.

\begin{tabular}{l}
  
\end{tabular}

{\tmstrong{Proof}} It suffices to note that the extreme value of the quadratic
polynomial $x r^2 + 2 y r + z$ is equal to $\frac{1}{x}$, while the extreme
value of the quadratic polynomial $z q^2 + 2 y q + x$ is equal to
$\frac{1}{z}$. $\Box$

\begin{tabular}{l}
  
\end{tabular}

An alternative way of looking at this situation is as follows. For any triple
$(x, y, z)$ consider the binary quadratic form $Q (s, t) = x s^2 + 2 y s t + z
t^2$. Then the triple $(x, y, z)$ solves the equation (6.3) if and only if the
corresponding quadratic form $Q$ is positive definite and has a discriminant
which is equal to $- 4$. And so the argument just made turns out to be
equivalent to the longstanding wisdom that all quadratic forms with this
property are equivalent. Now let $A \epsilon \tmop{SL} (2, \mathbb{Z})$ be
such that

\[ \Psi (A) \left(\begin{array}{c}
     \begin{array}{l}
       \mathfrak{m}\\
       k_m\\
       l_m
     \end{array}
   \end{array}\right) = \left(\begin{array}{c}
     \begin{array}{l}
       1\\
       0\\
       1
     \end{array}
   \end{array}\right) . \]
\begin{tabular}{l}
  
\end{tabular}

{\tmstrong{Remark}} The existence of such a matrix $A$ can also be established
through an application of part c) in Lemma 4.2.

\begin{tabular}{l}
  
\end{tabular}

It is not imperative to choose the second column of the matrix
$\left(\begin{array}{ccc}
  \mathfrak{c} & \mathfrak{m} & \mathfrak{a}\\
  k_c & k_m & k_a\\
  l_c & l_m & l_a
\end{array}\right)$. Part of the discus-

sion in the present section could be based on the other two choices as well.
The expediency of choosing the second column will become clear in Proposition
6.5 and Section 11, however. \ If we write
\[ \Psi (A) \left(\begin{array}{ccc}
     \mathfrak{c} & \mathfrak{m} & \mathfrak{a}\\
     k_c & k_m & k_a\\
     l_c & l_m & l_a
   \end{array}\right) = \left(\begin{array}{ccc}
     x_1 & 1 & x_2\\
     y_1 & 0 & y_2\\
     z_1 & 1 & z_2
   \end{array}\right) ; x_i, y_i, z_i  \begin{array}{l}
     \epsilon
   \end{array} \mathbb{Z}; i \begin{array}{l}
     \epsilon
   \end{array} \begin{array}{l}
     \{ 1, 2 \}
   \end{array}, \]
then,

(6.4)
\[ \left(\begin{array}{ccc}
     x_1 & 1 & x_2\\
     y_1 & 0 & y_2\\
     z_1 & 1 & z_2
   \end{array}\right) \begin{array}{l}
     \tmop{adj}\\
     \\
     
   \end{array} = \left(\begin{array}{c}
     \begin{array}{l}
       \begin{array}{l}
         
       \end{array} - y_2\\
       {y_2}  z_1 - y_1 z_2 \begin{array}{l}
         
       \end{array}\\
       \begin{array}{l}
         
       \end{array} \begin{array}{l}
         
       \end{array} \begin{array}{l}
         
       \end{array} y_1 \begin{array}{l}
         
       \end{array}
     \end{array} \begin{array}{l}
       \begin{array}{l}
         
       \end{array} x_2 - z_2\\
       x_1 z_2 - x_2 z_1 \begin{array}{l}
         
       \end{array}\\
       \begin{array}{l}
         
       \end{array} z_1 - x_1
     \end{array} \begin{array}{l}
       \begin{array}{l}
         
       \end{array} \begin{array}{l}
         
       \end{array} y_2\\
       x_2 y_1 - x_1 y_2\\
       \begin{array}{l}
         
       \end{array} \begin{array}{l}
         
       \end{array} - y_1
     \end{array}\\
     
   \end{array}\right), \]
and by (5.26)

(6.5)
\[ \left(\begin{array}{ccc}
     x_1 & 1 & x_2\\
     y_1 & 0 & y_2\\
     z_1 & 1 & z_2
   \end{array}\right) \begin{array}{l}
     \tmop{adj}\\
     \\
     
   \end{array} = \left(\begin{array}{ccc}
     l_c + 3 k_c & - (2 k_c + 3\mathfrak{c}) & \mathfrak{c}\\
     - (l_b + 3 \nu k_b) & 2 k_b + 3 \nu \mathfrak{b} & -\mathfrak{b}\\
     l_a - 3 k_a & - (2 k_a - 3\mathfrak{a}) & \mathfrak{a}
   \end{array}\right) \Psi (A)^{- 1} \]
Since an application from the left of the matrix $\Psi (A)^{- 1}$ to a row
vector corresponds to the transformation of the affiliated binary quadratic
form by the matrix $A^{- 1}$, and since the first and the third row of the
matrix on the right hand side of (6.4) correspond to symmetric forms, we
conclude that the quadratic forms corresponding to the first and the third row
of the matrix

$\left(\begin{array}{ccc}
  l_c + 3 k_c & - (2 k_c + 3\mathfrak{c}) & \mathfrak{c}\\
  - (l_b + 3 \nu k_b) & 2 k_b + 3 \nu \mathfrak{b} & -\mathfrak{b}\\
  l_a - 3 k_a & - (2 k_a - 3\mathfrak{a}) & \mathfrak{a}
\end{array}\right)$ are equivalent to a symmetric form each. This leads to the

following significant conclusion.

\begin{tabular}{l}
  
\end{tabular}

{\tmstrong{6.2 Proposition}} Every cycle of reduced binary quadratic forms
including a Markoff form also includes a symmetric form.

\begin{tabular}{l}
  
\end{tabular}

It has been known for a long time that every Markoff form $F$ is equivalent to
$- F$. Since a symmetric form $H$ is obviously equivalent to $- H$, that
statement follows immediately from Proposition 6.2. However, we should like to
point out that Proposition 6.2 can be derived directly from the equivalence of
the forms $F$ and $- F$. This can be seen as follows. A matrix $M \epsilon
\tmop{SL} (2, \mathbb{Z})$ which transforms a form $Q$ into $- Q$ has trace
zero. All matrices in $\tmop{SL} (2, \mathbb{Z})$ with trace zero are similar
over $\tmop{SL} (2, \mathbb{Z})$ to a matrix with zeros in the diagonal. If
$T^{- 1} M T$ has zeros in the diagonal for some $T \epsilon \tmop{SL} (2,
\mathbb{Z})$, and $Q$ is transformed by $T$ into $H$, then $H$ is symmetric.

We turn now to the second row of the matrix on the right hand side of (6.4).
First we note that,

(6.6)
\[ \tmop{entry} (2, 1) + \tmop{entry} (2, 3) = {(y_2}  z_1 - y_1 z_2) + (x_2
   y_1 - x_1 y_2) \]
\[ = \det \left(\begin{array}{ccc}
     x_1 & 1 & x_2\\
     y_1 & 0 & y_2\\
     z_1 & 1 & z_2
   \end{array}\right) = \det \left(\begin{array}{ccc}
     \mathfrak{c} & \mathfrak{m} & \mathfrak{a}\\
     k_c & k_m & k_a\\
     l_c & l_m & l_a
   \end{array}\right) = 2 . \]
It follows from (6.4) and (6.5) that the discriminant of this row is equal to
$9\mathfrak{b}^2 - 4$. Hence,
\[ (x_1 z_2 - x_2 z_1)^2 - 4 {(y_2}  z_1 - y_1 z_2) (x_2 y_1 - x_1 y_2) \]
\[ = (x_1 z_2 - x_2 z_1)^2 - 4 {(y_2}  z_1 - y_1 z_2) {(2 - (y_2}  z_1 - y_1
   z_2)) \]
\[ = (x_1 z_2 - x_2 z_1)^2 + 4 {(y_2}  z_1 - y_1 z_2 - 1)^2 - 4 =
   9\mathfrak{b}^2 - 4, \]
which leads to

(6.7)
\[ (x_1 z_2 - x_2 z_1)^2 + 4 {(y_2}  z_1 - y_1 z_2 - 1)^2 = 9\mathfrak{b}^2 .
\]
Since the sum of two squares is divisible by 3 if and only if each summand
shares this property, we conclude that
\[ \mathfrak{p}= \frac{1}{3} (x_1 z_2 - x_2 z_1) \epsilon \mathbb{Z},
   \begin{array}{l}
     
   \end{array} \mathfrak{q}= \frac{1}{3} {(y_2}  z_1 - y_1 z_2 - 1) \epsilon
   \mathbb{Z}, \]
so that (6.7) becomes

(6.8)
\[ \mathfrak{p}^2 + 4\mathfrak{q}^2 =\mathfrak{b}^2 . \]
As we shall see below (Proposition 6.5), the Pythagorean triple
$(\mathfrak{p}, 2\mathfrak{q}, \mathfrak{b})$ is primitive in case
$\mathfrak{p}$ is odd, and the Pythagorean triple $\left(
\frac{\mathfrak{p}}{2}, \mathfrak{q}, \frac{\mathfrak{b}}{2} \right)$ is
primitive in case $\mathfrak{p}$ is even. Hence, the standard parametrization
for Pythagorean triples ensures the existence of two integers, $\mathfrak{f}
\tmop{and} \mathfrak{g}$, such that

(6.9)
\[ \mathfrak{p}=\mathfrak{f}^2 -\mathfrak{g}^2, \begin{array}{l}
     
   \end{array} \mathfrak{q}=\mathfrak{f}\mathfrak{g}, \begin{array}{l}
     
   \end{array} \mathfrak{b}=\mathfrak{f}^2 +\mathfrak{g}^2 . \]
In conclusion, the second row vector of the matrix on the right hand side of
(6.4) takes the form,
\[ {(y_2}  z_1 - y_1 z_2, \begin{array}{l}
     
   \end{array} x_1 z_2 - x_2 z_1, \begin{array}{l}
     
   \end{array} x_2 y_1 - x_1 y_2) = (1 + 3\mathfrak{f}\mathfrak{g},
   \begin{array}{l}
     
   \end{array} 3 (\mathfrak{f}^2 -\mathfrak{g}^2), \begin{array}{l}
     
   \end{array} 1 - 3\mathfrak{f}\mathfrak{g}) . \]
To summarize, we have arrived at the following situation. Given a Markoff
number $\mathfrak{m}$, there exist three equivalent quadratic forms,

(6.10)
\[ F (s, t) =\mathfrak{m}s^2 - (2 k - 3\mathfrak{m}) s t + (l - 3 k) t^2,
   \begin{array}{l}
     
   \end{array} k^2 + 1 =\mathfrak{m}l, \begin{array}{l}
     
   \end{array} 0 < 2 k <\mathfrak{m} \]
\[ G (s, t) = (1 - 3\mathfrak{f}\mathfrak{g}) s^2 + 3 (\mathfrak{g}^2
   -\mathfrak{f}^2) s t + (1 + 3\mathfrak{f}\mathfrak{g}) t^2, \]
\[ H (s, t) =\mathfrak{u}s^2 +\mathfrak{v}s t -\mathfrak{u}t^2, \]
all of which, after having been subjected to a transformation by the matrix
$\left(\begin{array}{c}
  \begin{array}{l}
    \begin{array}{l}
      
    \end{array} 0\\
    - 1 \begin{array}{l}
      
    \end{array}
  \end{array} \begin{array}{l}
    1\\
    0
  \end{array}
\end{array}\right)$, if necessary, may be assumed to be reduced. Standard
theory for binary quadratic forms (cf. [L], Satz 202) ensures that there exist
(fundamental) automorphs of these quadratic forms, which in each of these
three particular cases take the form, in the order of their appearance above,

(6.11)
\[ \mathfrak{F}= \left(\begin{array}{c}
     \begin{array}{l}
       3\mathfrak{m}- k\\
       \begin{array}{l}
         
       \end{array} \begin{array}{l}
         -
       \end{array} \mathfrak{m}
     \end{array} \begin{array}{l}
       - 3 k + l\\
       \begin{array}{l}
         
       \end{array} \begin{array}{l}
         
       \end{array} k
     \end{array}
   \end{array}\right), \]
\[ \mathfrak{G}= \left(\begin{array}{c}
     \begin{array}{l}
       \begin{array}{l}
         
       \end{array} 3\mathfrak{g}^2\\
       3\mathfrak{f}\mathfrak{g}- 1
     \end{array} \begin{array}{l}
       3\mathfrak{f}\mathfrak{g}+ 1\\
       \begin{array}{l}
         
       \end{array} 3\mathfrak{f}^2
     \end{array}
   \end{array}\right), \]
\[ \mathfrak{H}= \left(\begin{array}{c}
     \begin{array}{l}
       \frac{3\mathfrak{m}-\mathfrak{v}}{2} \begin{array}{l}
         
       \end{array}\\
       \begin{array}{l}
         
       \end{array} \begin{array}{l}
         
       \end{array} \mathfrak{u}
     \end{array} \begin{array}{l}
       \begin{array}{l}
         
       \end{array} \begin{array}{l}
         
       \end{array} \mathfrak{u}\\
       \frac{3\mathfrak{m}+\mathfrak{v}}{2}
     \end{array}
   \end{array}\right) . \]

Recall from [F], Gesammelte Abhandlungen, Band III, p.606 (IV) that a number
$\mathfrak{m}$ is Markoff if and only if $\mathfrak{m}$ is representable by a
quadratic form $Q$ which is equivalent to $- Q$, and which has the
discriminant $9\mathfrak{m}^2 - 4$. We can now give an alternative
characterization of Markoff numbers.

\

{\tmstrong{6.3 Proposition}} An integer $\mathfrak{m} \geqslant 1$ is a
Markoff number if and only if $\mathfrak{m}= \frac{1}{3} \tmop{tr}
(\mathfrak{G})$, where

$\mathfrak{G}= \left(\begin{array}{c}
  \begin{array}{l}
    \begin{array}{l}
      
    \end{array} 3\mathfrak{g}^2\\
    3\mathfrak{f}\mathfrak{g}+ 1
  \end{array} \begin{array}{l}
    3\mathfrak{f}\mathfrak{g}- 1\\
    \begin{array}{l}
      
    \end{array} 3\mathfrak{f}^2
  \end{array}
\end{array}\right)$ for some integers $\mathfrak{f}, \mathfrak{g}$; and
$\mathfrak{G}$ is equivalent to a symmetric matrix.

\

{\tmstrong{Proof}} A unimodular 2x2 matrix $A$ which is equivalent to a
symmetric matrix is also equivalent to $A^t$. But a matrix which conjugates
$A$ to $A^t$ has to be symmetric too. This entails that $A$ is the product of
two symmetric matrices. Applied to a matrix of the form $\mathfrak{G}=
\left(\begin{array}{c}
  \begin{array}{l}
    \begin{array}{l}
      
    \end{array} 3\mathfrak{g}^2\\
    3\mathfrak{f}\mathfrak{g}+ 1
  \end{array} \begin{array}{l}
    3\mathfrak{f}\mathfrak{g}- 1\\
    \begin{array}{l}
      
    \end{array} 3\mathfrak{f}^2
  \end{array}
\end{array}\right)$, which is unimodular for any $\mathfrak{f}$ and
$\mathfrak{g}$, this means that $\mathfrak{G}=\mathfrak{S}\mathfrak{T}$, where
$\mathfrak{S}$ and $\mathfrak{T}$ are symmetric and unimodular. Since
\[ (\mathfrak{S}\mathfrak{T})^t =\mathfrak{T}\mathfrak{S}=
   \left(\begin{array}{c}
     \begin{array}{l}
       \begin{array}{l}
         
       \end{array} 3\mathfrak{g}^2\\
       3\mathfrak{f}\mathfrak{g}- 1
     \end{array} \begin{array}{l}
       3\mathfrak{f}\mathfrak{g}+ 1\\
       \begin{array}{l}
         
       \end{array} 3\mathfrak{f}^2
     \end{array}
   \end{array}\right), \]
and hence
\[ \mathfrak{S}^{- 1} \mathfrak{T}^{- 1} = \left(\begin{array}{c}
     \begin{array}{l}
       \begin{array}{l}
         
       \end{array} 3\mathfrak{f}^2\\
       - 3\mathfrak{f}\mathfrak{g}+ 1
     \end{array} \begin{array}{l}
       - 3\mathfrak{f}\mathfrak{g}- 1\\
       \begin{array}{l}
         \begin{array}{l}
           
         \end{array}
       \end{array} 3\mathfrak{g}^2
     \end{array}
   \end{array}\right), \]
we get
\[ \mathfrak{S}\mathfrak{T}\mathfrak{S}^{- 1} \mathfrak{T}^{- 1} =
   \left(\begin{array}{c}
     \begin{array}{l}
       \begin{array}{l}
         
       \end{array} 3\mathfrak{g}^2\\
       3\mathfrak{f}\mathfrak{g}+ 1
     \end{array} \begin{array}{l}
       3\mathfrak{f}\mathfrak{g}- 1\\
       \begin{array}{l}
         
       \end{array} 3\mathfrak{f}^2
     \end{array}
   \end{array}\right)  \left(\begin{array}{c}
     \begin{array}{l}
       \begin{array}{l}
         
       \end{array} 3\mathfrak{f}^2\\
       - 3\mathfrak{f}\mathfrak{g}+ 1
     \end{array} \begin{array}{l}
       - 3\mathfrak{f}\mathfrak{g}- 1\\
       \begin{array}{l}
         \begin{array}{l}
           
         \end{array}
       \end{array} 3\mathfrak{g}^2
     \end{array}
   \end{array}\right) = \left(\begin{array}{c}
     \begin{array}{l}
       6\mathfrak{f}\mathfrak{g}- 1\\
       \begin{array}{l}
         
       \end{array} 6\mathfrak{f}^2
     \end{array} \begin{array}{l}
       \begin{array}{l}
         
       \end{array} \begin{array}{l}
         
       \end{array} - 6\mathfrak{g}^2\\
       \begin{array}{l}
         
       \end{array} - 6\mathfrak{f}\mathfrak{g}  - 1
     \end{array}
   \end{array}\right) . \]
Hence $\tmop{tr} (\mathfrak{S}\mathfrak{T}\mathfrak{S}^{- 1} \mathfrak{T}^{-
1}) = - 2$. Now Fricke's identity implies,
\[ 9\mathfrak{m}^2 + (\tmop{tr} (\mathfrak{S}))^2 + (\tmop{tr}
   (\mathfrak{T}))^2 = 3\mathfrak{m} \tmop{tr} (\mathfrak{S}) \tmop{tr}
   (\mathfrak{T}), \]
which means that $\mathfrak{m}$ has to be a Markoff number. Thus we have shown
that the enunciated condition is sufficient. That it is also necessary follows
from Proposition 6.2 and the fact that each form of type $F$ is equivalent to
a form of type $G$ (cf. (6.10); see also Remark 4 following Proposition 6.4
below). $\Box$

\

{\tmstrong{Remarks}} 1) If $\mathfrak{G}=\mathfrak{S}\mathfrak{T}$ is a
factorization as in the proof of Proposition 6.3, with $\mathfrak{S}$ and
$\mathfrak{T}$ being symmetric and unimodular, then
$\mathfrak{G}=\mathfrak{S}_n \mathfrak{T}_n$, where $\mathfrak{S}_n =
\mathfrak{G}^n \mathfrak{S}$ and $\mathfrak{T}_n =\mathfrak{T}\mathfrak{G}^{-
n}$, is also a factorization of $\mathfrak{G}$ into symmetric and unimodular
matrices. Since $\tmop{tr} (\mathfrak{T}_{n + 1}) = 3\mathfrak{m} \tmop{tr}
(\mathfrak{T}_n) - \tmop{tr} (\mathfrak{T}_{n - 1})$ this shows that, given
any Markoff triple $(\mathfrak{a}, \mathfrak{m}, \mathfrak{c})$, there exists
a factorization $\mathfrak{G}=\mathfrak{S}\mathfrak{T}$ of $\mathfrak{G}$ into
symmetric and unimodular matrices $\mathfrak{S}$ and $\mathfrak{T}$ such that
$\left\{ \frac{1}{3} \tmop{tr} (\mathfrak{S}) \nocomma, \frac{1}{3} \tmop{tr}
(\mathfrak{T}) \right\} = \{ \mathfrak{a}, \mathfrak{c} \}$. Moreover, there
are no other factorizations of $\mathfrak{G}$ into two symmetric and
unimodular matrices.

2) The line of reasoning in the proof of Proposition 6.3 can be adapted to
show that for any positive integer $n$ the quadratic forms $Q (s, t) = s^2 + n
s t + t^2$ and $- Q (s, t)$ are equivalent if and only if $n = 3$. This is of
course not new. \ The point of the proof below is to highlight the close
connection between quadratic forms of discriminant $n^2 - 4$, symmetric forms,
and Markoff numbers. \ Since $- Q (s, t)$ is equivalent to the form $- s^2 + n
s t - t^2$, our claim is equivalent to showing that the matrices
\[ A = \left(\begin{array}{c}
     \begin{array}{l}
       0\\
       1 \begin{array}{l}
         
       \end{array}
     \end{array} \begin{array}{l}
       - 1\\
       \begin{array}{l}
         
       \end{array} n
     \end{array}
   \end{array}\right) \begin{array}{l}
     
   \end{array} \tmop{and} \begin{array}{l}
     
   \end{array} A^t = \left(\begin{array}{c}
     \begin{array}{l}
       \begin{array}{l}
         
       \end{array} 0\\
       - 1 \begin{array}{l}
         
       \end{array}
     \end{array} \begin{array}{l}
       \begin{array}{l}
         1
       \end{array}\\
       \begin{array}{l}
         n
       \end{array}
     \end{array}
   \end{array}\right) \begin{array}{l}
     
   \end{array} \]
are equivalent. But if these two matrices are equivalent, then it follows
exactly as in the proof of Proposition 6.3 that $A$ is a product of two
symmetric unimodular matrices. Now the same manipulations as in the proof of
Proposition 6.3 show that $n$ together with the traces of the two symmetric
unimodular matrices appearing in such a factorization form a Markoff triple.
By Remark 1 all Markoff triples which include $n$ as a member can be realized
in that manner. In particular the factorization can be arranged in such a way
that $n$ is not smaller then the other two members of the Markoff triple.
Therefore, suppose that
\[ A = R S \nocomma, \begin{array}{l}
     
   \end{array} \tmop{where} \begin{array}{l}
     
   \end{array} R, S \epsilon \tmop{Sl} (2, \mathbb{Z}), \begin{array}{l}
     
   \end{array} R = R^t \nocomma, \begin{array}{l}
     
   \end{array} S = S^t, \begin{array}{l}
     
   \end{array} n \geqslant \max \{ \tmop{tr} (R), \tmop{tr} (S) \} . \]
Since $A^t = S R$, it follows that
\[ A R = R A^t, \begin{array}{l}
     
   \end{array} S A = A^t S. \]
Letting
\[ R = \left(\begin{array}{c}
     \begin{array}{l}
       r_1\\
       r_2 \begin{array}{l}
         
       \end{array}
     \end{array} \begin{array}{l}
       r_2\\
       r_3
     \end{array}
   \end{array}\right), \begin{array}{l}
     
   \end{array} S = \left(\begin{array}{c}
     \begin{array}{l}
       s_1\\
       s_2 \begin{array}{l}
         
       \end{array}
     \end{array} \begin{array}{l}
       s_2\\
       s_3
     \end{array}
   \end{array}\right), \]
these two identities are equivalent to
\[ r_1 + n r_2 + r_3 = 0, \begin{array}{l}
     
   \end{array} \tmop{and} \begin{array}{l}
     
   \end{array} s_1 - n s_2 + s_3 = 0, \begin{array}{l}
     
   \end{array} \tmop{respectively} . \]
Since $r_2$ and $s_2$ cannot be equal to zero, we conclude that
\[ n \leqslant r_1 + r_3, \begin{array}{l}
     
   \end{array} n \leq s_1 + s_3 . \]
But since $n \geqslant \max \{ \tmop{tr} (R), \tmop{tr} (S) \}$, and since
$(\tmop{tr} (R), n, \tmop{tr} (S))$ is a Markoff triple, it follows that $n =
3$, as claimed. On the other hand, we do have a factorization into two
symmetric unimodular matrices in case $n = 3$, namely
\[ \left(\begin{array}{c}
     \begin{array}{l}
       0\\
       1 \begin{array}{l}
         
       \end{array}
     \end{array} \begin{array}{l}
       - 1\\
       \begin{array}{l}
         
       \end{array} 3
     \end{array}
   \end{array}\right) = \left(\begin{array}{c}
     \begin{array}{l}
       \begin{array}{l}
         
       \end{array} 1\\
       - 1
     \end{array} \begin{array}{l}
       - 1\\
       \begin{array}{l}
         
       \end{array} 2
     \end{array}
   \end{array}\right) \left(\begin{array}{c}
     \begin{array}{l}
       1\\
       1
     \end{array} \begin{array}{l}
       1\\
       2
     \end{array}
   \end{array}\right) . \]
Note though, that in general a trival factorization of the matrix $A$ into two
symmetric matrices belonging to $\tmop{GL} (2, \mathbb{Z})$ does always exist,
namely
\[ A = \left(\begin{array}{c}
     \begin{array}{l}
       - 1\\
       \begin{array}{l}
         
       \end{array} 0
     \end{array} \begin{array}{l}
       0\\
       1
     \end{array}
   \end{array}\right) \left(\begin{array}{c}
     \begin{array}{l}
       0\\
       1 \begin{array}{l}
         
       \end{array}
     \end{array} \begin{array}{l}
       1\\
       n
     \end{array}
   \end{array}\right) . \]

Proposition 6.2 in combination with the special form of the discriminant
allows us to give a more incisive characterization of the cycles of reduced
forms containing a Markoff form.

\begin{tabular}{l}
  
\end{tabular}

{\tmstrong{6.4 Proposition}} Any cycle of reduced forms containing a Markoff
form associated with a Markoff number $\mathfrak{m} \geqslant 5$ contains two
(distinct) symmetric forms $H_1 (s, t) =\mathfrak{u}_1 s^2 +\mathfrak{v}_1 s t
-\mathfrak{u}_1 t^2$ and $H_2 (s, t) =\mathfrak{u}_2 s^2 +\mathfrak{v}_2 s t
-\mathfrak{u}_2 t^2$ such that $\mathfrak{u}_1 \neq -\mathfrak{u}_2$.

\begin{tabular}{l}
  
\end{tabular}

{\tmstrong{Proof}} By Proposition 6.2, any cycle of reduced forms containing a
Markoff form $F$ with discriminant $9\mathfrak{m}^2 - 4$ contains a symmetric
form $H $ with a fundamental automorph
\[ \mathfrak{H}= \left(\begin{array}{c}
     \begin{array}{l}
       \frac{3\mathfrak{m}-\mathfrak{v} }{2} \begin{array}{l}
         
       \end{array}\\
       \begin{array}{l}
         
       \end{array} \begin{array}{l}
         
       \end{array} \mathfrak{u}
     \end{array} \begin{array}{l}
       \begin{array}{l}
         
       \end{array} \begin{array}{l}
         
       \end{array} \mathfrak{u}\\
       \frac{3\mathfrak{m}+\mathfrak{v}}{2}
     \end{array}
   \end{array}\right), \]
where

(6.12)
\[ 4\mathfrak{u}^2 +\mathfrak{v}^2 = 9\mathfrak{m}^2 - 4 \text{.} \]
First we deal with the case when $\mathfrak{m}$ is odd. Employing the standard
parametrization for Pythagorean quadruples (cf. [M2], p.14) we conclude that
there exist integers $n, p, q, r$ such that

(6.13)
\[ 1 = n r - p q, \begin{array}{l}
     
   \end{array} \begin{array}{l}
     
   \end{array} \mathfrak{u}= n q + p r, \]
\[ \mathfrak{v}= - n^2 - p^2 + q^2 + r^2, \begin{array}{l}
     
   \end{array} \begin{array}{l}
     
   \end{array} 3\mathfrak{m}= n^2 + p^2 + q^2 + r^2 . \]
Hence,

(6.14)
\[ \mathfrak{T}  \mathfrak{T}^t =\mathfrak{H}, \begin{array}{l}
     
   \end{array} \tmop{where} \begin{array}{l}
     
   \end{array} \mathfrak{T}= \left(\begin{array}{c}
     \begin{array}{l}
       n\\
       q
     \end{array} \begin{array}{l}
       p\\
       r
     \end{array}
   \end{array}\right), \det (\mathfrak{T}) = 1 \]
We claim that $n \neq r$. Suppose this were not true. Then,
\[ \mathfrak{u}= n (q + p), \begin{array}{l}
     
   \end{array} \mathfrak{v}= (q + p) (q - p) . \]
It follows from (6.12)

(6.15)
\[ q + p \text{ divides the discriminant } 9\mathfrak{m}^2 - 4 . \]
Moreover, since the three coefficients of $H$ are divisible by $q + p$, every
integer represented by $H$ is divisible by $q + p$. Since $H \tmop{and} F$ are
equivalent, and since $\mathfrak{m}$ is represented by $F$, $\mathfrak{m}$ is
represented by $H$ as well. It follows that $\mathfrak{m}$ is divisible by $q
+ p$. Combined with (6.15) this implies that $q + p$ is odd and that 4 is
divisible by $q + p$. Hence,

(6.16)
\[ q + p = 1 \begin{array}{l}
     
   \end{array} \tmop{or} \begin{array}{l}
     
   \end{array} q + p = - 1 \begin{array}{l}
     
   \end{array} \]
Combining the first identity in (6.16) with the first identity in (6.13)
yields,

(6.17)
\[ n^2 + p^2 = 1 + p . \]
This diophantine equation has four solutions, namely $(n, p) = (\pm 1, 0)$ \
and $(n, p) = (\pm 1, 1)$. Combining the second identity in (6.16) with the
first identity in (6.13) yields,

(6.18)
\[ n^2 + p^2 = 1 - p . \]
This diophantine equation has four solutions, namely $(n, p) = (\pm 1, - 1)$ \
and $(n, p) = (\pm 1, 0)$. In conclusion, for all solutions of (6.17) and
(6.18) we get $4\mathfrak{u}^2 +\mathfrak{v}^2 = 5$, which implies
$\mathfrak{m}= 1$, and therefore $n \neq r$ as claimed, in case $\mathfrak{m}$
is odd and $\mathfrak{m} \geqslant 5$.

We turn to the case of an even Markoff number $\mathfrak{m}$. In this case
(6.12) implies that $\mathfrak{v}$ is even. Letting $\mathfrak{w}=
\frac{\mathfrak{v}}{2}$, (6.12) turns into

(6.19)
\[ \mathfrak{u}^2 +\mathfrak{w}^2 = 9\mathbb{m}^2 - 1 . \]
But this implies that $\mathfrak{u} \tmop{and} \mathfrak{v}$ have to be even.
Hence, letting $\mathbb{u}= \frac{\mathfrak{u}}{2}, \begin{array}{l}
  
\end{array} \mathbb{v}= \frac{\mathfrak{v}}{2}$, (6.19) turns into

(6.20)
\[ (2\mathbb{u})^2 + (2\mathbb{v})^2 = 9\mathbb{m}^2 - 1 . \]
This time the parametrization takes the following form. There exist integers
$n, p, q, r$, such that

(6.21)
\[ \mathbb{u}= n r - p q, \begin{array}{l}
     
   \end{array} \begin{array}{l}
     
   \end{array} \mathbb{v}= n q + p r, \]
\[ 1 = n^2 + p^2 - q^2 - r^2, \begin{array}{l}
     
   \end{array} \begin{array}{l}
     
   \end{array} 3\mathbb{m}= n^2 + p^2 + q^2 + r^2 . \]
Hence,

(6.22)
\[ \mathfrak{T}  \mathfrak{T}^t =\mathfrak{H}, \begin{array}{l}
     
   \end{array} \tmop{where} \begin{array}{l}
     
   \end{array} \mathfrak{T}= \left(\begin{array}{c}
     \begin{array}{l}
       n - q\\
       r + p
     \end{array} \begin{array}{l}
       r - p\\
       n + q
     \end{array}
   \end{array}\right), \det (\mathfrak{T}) = 1. \]
We claim that $n - q \neq n + q$. Suppose this were not true. Then $q = 0$,
and therefore,

(6.23)
\[ \mathbb{u}= n r, \begin{array}{l}
     
   \end{array} \begin{array}{l}
     
   \end{array} \mathbb{v}= p r. \]
It follows from (6.12),

(6.24)
\[ r \text{ divides the discrminant } 9\mathfrak{m}^2 - 4 . \]
Moreover, since the three coefficients of $H$ are divisible by $r$, every
integer represented by $H$ is divisible by $r$. Since $H \tmop{and} F$ are
equivalent, and since $\mathfrak{m}$ is represented by $F$, $\mathfrak{m}$ is
represented by $H$ as well. It follows that $\mathfrak{m}$ is divisible by
$r$. Combined with (6.20) this implies that that 2 is divisible by $r$. If $|
r | = 2$, then the third identity in (6.21) yields $n^2 + p^2 = 5$.
Substituting this into the fourth identity in (6.21) leads to $\mathbb{m}= 3$,
hence to $\mathfrak{m}= 6$, which is not a Markoff number. If $| r | = 1$,
then the third identity in (6.21) implies that $n^2 + p^2 = 2$. Substituting
this into the fourth identity in (6.21) leads to $\mathbb{m}= 1$, hence to
$\mathfrak{m}= 2$. In conclusion, $n - q \neq n + q$ as claimed, in case
$\mathfrak{m} \geqslant 5$ is an even Markoff number.

Combining the two separate cases for $\mathfrak{m}$, we have shown that there
always exists a matrix $\mathfrak{T} \epsilon \tmop{SL} (2, \mathbb{Z})$ with
distinct diagonal entries such that $\mathfrak{T}\mathfrak{T}^t
=\mathfrak{H}$. Considering such a matrix $\mathfrak{T}_1$ for
$\mathfrak{H}_1$, let $\mathfrak{H}_2 =\mathfrak{T}^t_1 \mathfrak{T}_1$. Then,
\[ \mathfrak{H}_2 =\mathfrak{T}^{- 1}_1 \mathfrak{H}_1 \mathfrak{T}_1 . \]
Thus, $\mathfrak{H}_2$ is a fundamental automorph for a symmetric form $H_2$
which is equivalent to $H_1$. Since the diagonal entries of $\mathfrak{H}_1$
are distinct, the sum of the first coefficient of $H_1 \tmop{and} H_2$ can not
be zero $\Box$

\

{\tmstrong{Remarks}} 1) If $\mathfrak{m}= 1$ or $\mathfrak{m}= 2$, then the
(reduced) Markoff form is symmetric.

2) Proposition 6.4 places all cycles of reduced forms which include Markoff
forms with a discriminant larger than 32 among the so-called ambiguous cycles
(as defined in [BV]), or alternatively (by extension of [H-K], Definition
2.3.1) weakly ambiguous cycles. Note however, that it follows from Remark 2
following Proposition 6.3, that the class in the form class group
corresponding to these cycles is of order four. (Applying, for instance, [H-K]
Theorem 6.4.5(3) to the form $H$ composed with itself shows that the square of
the form class containing $H$ represents the number $- 1$. A more detailed
discussion of this issue will appear in Section 8) \ By contrast, the order of
the class corresponding to these cycles in the ideal class group is equal to
two. In terms of the classification scheme of cycles exhibited in [Bu], pp.
28-29, the cycles containing two ``non-affiliated'' symmetric forms (that
means that the two forms do not just differ by a \ minus sign) are being
addressed as ``Type 20''. It follows in particular that a Markoff form with a
discriminant larger than 32 can never be equivalent to a reduced ambiguous
form. For more information about \ the computational aspects of reduced forms
see [BV], Chapter 6, and for a comprehensive exposition of ambiguous classes
see [H-K], 5.6. Since the norm of the fundamental unit in our quadratic number
field is equal to one, Proposition 6.4 is a special case of [H-K], Theorem
5.6.9, which is derived through continued fraction expansions.

3) While every Markoff form is contained in a cycle of ``Type 20'', this
property does not characterize the corresponding ideal-class per se. In fact
one half of the number of ideal-classes whose square is equal to the principal
class correspond to cycles of ``Type 20''. Or put another way, the number of
ideal-classes which correspond to cycles of ``Type 20'' is equal to the number
of ideal-classes corresponding to cycles which contain at least one ambiguous
quadratic form, and hence exactly two such forms (see [H-K], Theorem 5.6.9).
Moreover, this number is determined by the prime factorization of
$9\mathfrak{m}^2 - 4$ . For a more detailed discussion see [Mo], Chapter 6,
and op. cit. pp. 352-353. But as we shall see in the next section, two Markoff
forms of discrimiant $9\mathfrak{m}^2 - 4$ have to be (properly or improperly)
equivalent nevertheless, because any Markoff form of discrimiant
$9\mathfrak{m}^2 - 4$ represents the number $\mathfrak{m}$.

4) Any unimodular integral matrix of the form \ $\mathfrak{F}=
\left(\begin{array}{c}
  \begin{array}{l}
    3\mathfrak{m}- k\\
    \begin{array}{l}
      
    \end{array} \begin{array}{l}
      
    \end{array} \mathfrak{m}
  \end{array} \begin{array}{l}
    3 k - l\\
    \begin{array}{l}
      
    \end{array} \begin{array}{l}
      
    \end{array} k
  \end{array}
\end{array}\right)$, $\mathfrak{m}$ being a positive integer, is equivalent to
a matrix of the form $\mathfrak{G}= \left(\begin{array}{c}
  \begin{array}{l}
    \begin{array}{l}
      
    \end{array} 3\mathfrak{g}^2\\
    3\mathfrak{f}\mathfrak{g}+ 1
  \end{array} \begin{array}{l}
    3\mathfrak{f}\mathfrak{g}- 1\\
    \begin{array}{l}
      
    \end{array} 3\mathfrak{f}^2
  \end{array}
\end{array}\right)$, where $\mathfrak{f}$ and $\mathfrak{g}$ are relatively
prime. To see this we note that the unimodularity of $\mathfrak{F}=
\left(\begin{array}{c}
  \begin{array}{l}
    3\\
    1
  \end{array} \begin{array}{l}
    - 1\\
    \begin{array}{l}
      
    \end{array} 0
  \end{array}
\end{array}\right) \left(\begin{array}{c}
  \begin{array}{l}
    \mathfrak{m}\\
    k
  \end{array} \begin{array}{l}
    k\\
    l
  \end{array}
\end{array}\right)$ is equivalent to the unimodularity of the matrix
$\mathfrak{S}= \left(\begin{array}{c}
  \begin{array}{l}
    \mathfrak{m}\\
    k
  \end{array} \begin{array}{l}
    k\\
    l
  \end{array}
\end{array}\right)$. Hence, by Lemma 4.2 part c) there exists a matrix
$\mathfrak{T}= \left(\begin{array}{c}
  \begin{array}{l}
    \mathfrak{g}\\
    \mathfrak{f}
  \end{array} \begin{array}{l}
    s\\
    t
  \end{array}
\end{array}\right) \epsilon \tmop{SL} (2, \mathbb{Z})$, such that
${\mathfrak{S}=\mathfrak{T}^t}  \mathfrak{T}$. Since $\mathfrak{g}t
-\mathfrak{f}s = 1$ it follows that $\mathfrak{T}\mathfrak{F}\mathfrak{T}^{-
1} =\mathfrak{T} \left(\begin{array}{c}
  \begin{array}{l}
    3\\
    1
  \end{array} \begin{array}{l}
    - 1\\
    \begin{array}{l}
      
    \end{array} 0
  \end{array}
\end{array}\right) \mathfrak{T}^t =\mathfrak{G}$. Conversely, if
$\mathfrak{f}$ and $\mathfrak{g}$ are relatively prime integers such that
$\mathfrak{m}=\mathfrak{f}^2 +\mathfrak{g}^2$, then there exist integers $s$
and $t$ such that $\mathfrak{g}t -\mathfrak{f}s = 1$. Letting $\mathfrak{T}=
\left(\begin{array}{c}
  \begin{array}{l}
    \mathfrak{g}\\
    \mathfrak{f}
  \end{array} \begin{array}{l}
    s\\
    t
  \end{array}
\end{array}\right)$, $l = s^2 + t^2$, and $k =\mathfrak{g}s +\mathfrak{f}t$,
it follows that
\[ \mathfrak{T}^{- 1} \left(\begin{array}{c}
     \begin{array}{l}
       \begin{array}{l}
         
       \end{array} 3\mathfrak{g}^2\\
       3\mathfrak{f}\mathfrak{g}+ 1
     \end{array} \begin{array}{l}
       3\mathfrak{f}\mathfrak{g}- 1\\
       \begin{array}{l}
         
       \end{array} 3\mathfrak{f}^2
     \end{array}
   \end{array}\right) \mathfrak{T}= \left(\begin{array}{c}
     \begin{array}{l}
       3\mathfrak{m}- k\\
       \begin{array}{l}
         
       \end{array} \begin{array}{l}
         
       \end{array} \mathfrak{m}
     \end{array} \begin{array}{l}
       3 k - l\\
       \begin{array}{l}
         
       \end{array} \begin{array}{l}
         
       \end{array} k
     \end{array}
   \end{array}\right) . \]
Likewise, by the same line of reasoning, any unimodular integral matrix of the
form \ $\mathfrak{F}= \left(\begin{array}{c}
  \begin{array}{l}
    3\mathfrak{m}- k\\
    \begin{array}{l}
      
    \end{array} \begin{array}{l}
      
    \end{array} -\mathfrak{m}
  \end{array} \begin{array}{l}
    - 3 k + l\\
    \begin{array}{l}
      
    \end{array} \begin{array}{l}
      
    \end{array} k
  \end{array}
\end{array}\right) = \left(\begin{array}{c}
  \begin{array}{l}
    \begin{array}{l}
      
    \end{array} 3\\
    - 1
  \end{array} \begin{array}{l}
    1\\
    0
  \end{array}
\end{array}\right) \left(\begin{array}{c}
  \begin{array}{l}
    \begin{array}{l}
      
    \end{array} \mathfrak{m}\\
    - k
  \end{array} \begin{array}{l}
    - k\\
    \begin{array}{l}
      
    \end{array} l
  \end{array}
\end{array}\right)$, $\mathfrak{m}$ being a positive integer, is equivalent to
a matrix of the form $\mathfrak{G}= \left(\begin{array}{c}
  \begin{array}{l}
    \begin{array}{l}
      
    \end{array} 3\mathfrak{g}^2\\
    3\mathfrak{f}\mathfrak{g}- 1
  \end{array} \begin{array}{l}
    3\mathfrak{f}\mathfrak{g}+ 1\\
    \begin{array}{l}
      
    \end{array} 3\mathfrak{f}^2
  \end{array}
\end{array}\right)$.

5) The matrix on the right hand side of the identity (6.5) can be written as
follows,
\[ \left(\begin{array}{ccc}
     1 & 0 & 0\\
     0 & - 1 & 0\\
     0 & 0 & 1
   \end{array}\right) \left(\begin{array}{ccc}
     \mathfrak{c} & \mathfrak{b} & \mathfrak{a}\\
     k_c + 3\mathfrak{c} & k_b & k_a\\
     l_c + 6 k_c + 9\mathfrak{c} & l_b & l_a
   \end{array}\right) \begin{array}{l}
     t\\
     \\
     
   \end{array} \left(\begin{array}{ccc}
     1 & 0 & 0\\
     0 & - 1 & 0\\
     0 & 0 & 1
   \end{array}\right) \left(\begin{array}{ccc}
     0 & - 3 & 1\\
     - 3 & 2 & 0\\
     1 & 0 & 0
   \end{array}\right) \Psi {(A)^{- 1}} , \]
in case $\nu = - 1$, and
\[ \left(\begin{array}{ccc}
     1 & 0 & 0\\
     0 & - 1 & 0\\
     0 & 0 & 1
   \end{array}\right) \left(\begin{array}{ccc}
     \mathfrak{c} & \mathfrak{b} & \mathfrak{a}\\
     k_c + 3\mathfrak{c} & k_b + 3\mathfrak{b} & k_a\\
     l_c + 6 k_c + 9\mathfrak{c} & l_b + 6 k_b + 9\mathfrak{b} & l_a
   \end{array}\right) \begin{array}{l}
     t\\
     \\
     
   \end{array} \left(\begin{array}{ccc}
     1 & 0 & 0\\
     0 & - 1 & 0\\
     0 & 0 & 1
   \end{array}\right) \left(\begin{array}{ccc}
     0 & - 3 & 1\\
     - 3 & 2 & 0\\
     1 & 0 & 0
   \end{array}\right) \Psi {(A)^{- 1}} , \]

in case $\nu = 1$. One can show that $\Psi (A) $ satisfies the identity
\[ \left(\begin{array}{ccc}
     0 & - 3 & 1\\
     - 3 & 2 & 0\\
     1 & 0 & 0
   \end{array}\right) \Psi {(A)^{- 1}}  = {(\Psi (A)^{- 1}} )^t
   \left(\begin{array}{ccc}
     0 & - 3 & 1\\
     - 3 & 2 & 0\\
     1 & 0 & 0
   \end{array}\right) \]
if and only if $A$ \ is of the from $\left(\begin{array}{c}
  \begin{array}{l}
    1\\
    r
  \end{array} \begin{array}{l}
    0\\
    1
  \end{array}
\end{array}\right)$. It is this identity that explains why we are staying
within the same equivalence class of quadratic residues on both sides of the
identity (6.5), as long as we apply matrices of the form $\Psi
(\left(\begin{array}{c}
  \begin{array}{l}
    1\\
    r
  \end{array} \begin{array}{l}
    0\\
    1
  \end{array}
\end{array}\right))$ only. That fact was implicitly instrumental in the proof
of Proposition 5.1.

\begin{tabular}{l}
  
\end{tabular}

We are now going to refine the analysis of the matrix in (6.4). Since $\Psi
(A)$ is an automorph of the quadratic form $Q$, we get,

(6.25)
\[ y^2_i + 1 = x_i z_i, \begin{array}{l}
     
   \end{array} i \begin{array}{l}
     \epsilon
   \end{array} \{ 1, 2 \} . \]
Since the first row of the matrix in (6.5) has the discriminant
$9\mathfrak{c}^2 - 4$, while the third row has the discriminant
$9\mathfrak{a}^2 - 4$, we also get,

(6.26)
\begin{eqnarray*}
  (x_1 - z_1)^2 + 4 y^2_1 = 9\mathfrak{a}^2 - 4, \begin{array}{l}
    
  \end{array} & (x_2 - z_2)^2 + 4 y^2_2 = 9\mathfrak{c}^2 - 4 & .
\end{eqnarray*}
Combining (6.25) and (6.26) for $i = 1 \tmop{and} i = 2$, respectively, leads
to,
\[ (x_1 + z_1)^2 = 9\mathfrak{a}^2, \begin{array}{l}
     
   \end{array} (x_2 + z_2)^2 = 9\mathfrak{c}^2 . \]
By Lemma 6.1

(6.27)
\[ x_1 + z_1 = 3\mathfrak{a} , \begin{array}{l}
     
   \end{array} x_2 + z_2 = 3\mathfrak{c} , \]
Letting
\[ v_a = x_1 - z_1, \begin{array}{l}
     
   \end{array} v_c = x_2 - z_2 \]
we can recast the identities (6.27) as follows,

(6.28)
\[ x_1 = \frac{1}{2} (3\mathfrak{a}  + v_a), \begin{array}{l}
     
   \end{array} z_1 = \frac{1}{2} (3\mathfrak{a}  - v_a), \begin{array}{l}
     
   \end{array} x_2 = \frac{1}{2} (3\mathfrak{c}  + v_c), \begin{array}{l}
     
   \end{array} z_2 = \frac{1}{2} (3\mathfrak{c}  - v_c) . \]
Substituting these expressions into the second row of the matrix in (6.4), and
combining the result with (6.8) yields,

(6.29)
\[ \left(\begin{array}{c}
     \begin{array}{l}
       \frac{1}{2} y_2 (3\mathfrak{a}  - v_a) - \frac{1}{2} y_1 (3\mathfrak{c}
       - v_c)
     \end{array}\\
     \frac{1}{4} (3\mathfrak{a}  + v_a) (3\mathfrak{c}  - v_c) - \frac{1}{4}
     (3\mathfrak{c}  + v_c) (3\mathfrak{a}  - v_a)\\
     \frac{1}{2} y_1 (3\mathfrak{c}  + v_c) - \frac{1}{2} y_2 (3\mathfrak{a} 
     + v_a)
   \end{array}\right) = \frac{1}{2} \left(\begin{array}{c}
     \begin{array}{l}
       y_2 (3\mathfrak{a}  - v_a) - y_1 (3\mathfrak{c}  - v_c)
     \end{array}\\
     3 (\mathfrak{c}  v_a -\mathfrak{a}  v_c)\\
     y_1 (3\mathfrak{c}  + v_c) - y_2 (3\mathfrak{a}  + v_a)
   \end{array}\right) = \left(\begin{array}{c}
     \begin{array}{l}
       1 + 3\mathfrak{q}
     \end{array}\\
     3\mathfrak{p}\\
     1 - 3\mathfrak{q}
   \end{array}\right) . \]
Writing the identities for the first and the third component as a linear
system in $y_1 \tmop{and} y_2$, and then solving for these two parameters
leads to,
\[ 3 (\mathfrak{a}  v_c -\mathfrak{c}  v_a) \left(\begin{array}{c}
     \begin{array}{l}
       y_1\\
       y_2
     \end{array}
   \end{array}\right) = \left(\begin{array}{c}
     \begin{array}{l}
       3\mathfrak{a}  + v_a\\
       3\mathfrak{c}  + v_c
     \end{array} \begin{array}{l}
       
     \end{array} \begin{array}{l}
       3\mathfrak{a}- v_a\\
       3\mathfrak{c}- v_c
     \end{array}
   \end{array}\right) \left(\begin{array}{c}
     \begin{array}{l}
       1 + 3\mathfrak{q}\\
       1 - 3\mathfrak{q}
     \end{array}
   \end{array}\right) = 6 \left(\begin{array}{c}
     \begin{array}{l}
       \mathfrak{q}v_a +\mathfrak{a}\\
       \mathfrak{q}v_c +\mathfrak{c}
     \end{array}
   \end{array}\right), \]
and finally, after invoking the second component identity in (6.29), letting
$u_a = y_1, u_c = y_2$,

(6.30)
\[ \mathfrak{p}u_a +\mathfrak{q}v_a = - 2\mathfrak{a}, \begin{array}{l}
     
   \end{array} \mathfrak{p}u_c +\mathfrak{q}v_c = - 2\mathfrak{c}. \]
In conclusion, we have shown the following.

\begin{tabular}{l}
  
\end{tabular}

{\tmstrong{{\tmem{6.5 Proposition}}}} Let $(\mathfrak{a}, \mathfrak{m},
\mathfrak{c})$ be a Markoff triple such that $\mathfrak{m}$ is dominant, and
let $\mathfrak{b}= 3\mathfrak{a}\mathfrak{c}-\mathfrak{m}$. Then the system of
diophantine equations

(6.31)
\[ p^2 + q^2 =\mathfrak{b}^2, \begin{array}{l}
     
   \end{array} u^2 + v^2 = 9\mathfrak{a}^2 - 4, \begin{array}{l}
     
   \end{array} p u + q v = - 2\mathfrak{a}, \]
has a solution with the following two properties:

(a) At least one of the two integers $u, v$ is even.

(b) If $u$ is even, the form
\[ H_a (s, t) = \frac{u}{2} s^2 + v s t - \frac{u}{2} t^2 \]
is equivalent to a form $G_a (s, t) = (1 + 3 f g) s^2 + 3 (f^2 - g^2) s t + (1
- 3 f g) t^2$. A similar statement holds if $v$ is even.

The claim of this proposition is also valid if the Markoff number
$\mathfrak{a}$ in (6.31) is replaced by $\mathfrak{c}.$

More specifically, with regard to property (b) in Proposition 6.5, (6.4) and
(6.5) show that the form $H_a$ is equivalent to a Markoff form associated with
$a$, or its opposite form, which in turn, by Remark 4 following Proposition
6.4, is equivalent to a form $G_a$, where $f$ and $g$ are relatively prime
integers.

\

{\tmstrong{Remark}} Note that the three identities in (6.31) can be expressed
in matrix form,

(6.32)
\[ \left(\begin{array}{c}
     \begin{array}{l}
       p\\
       u
     \end{array} \begin{array}{l}
       q\\
       v
     \end{array}
   \end{array}\right) \left(\begin{array}{c}
     \begin{array}{l}
       p\\
       u
     \end{array} \begin{array}{l}
       q\\
       v
     \end{array}
   \end{array}\right) \begin{array}{l}
     t\\
     \\
     
   \end{array} = \left(\begin{array}{c}
     \begin{array}{l}
       p\\
       u
     \end{array} \begin{array}{l}
       q\\
       v
     \end{array}
   \end{array}\right) \left(\begin{array}{c}
     \begin{array}{l}
       p\\
       q
     \end{array} \begin{array}{l}
       u\\
       v
     \end{array}
   \end{array}\right) = \left(\begin{array}{c}
     \begin{array}{l}
       \mathfrak{b}^2\\
       - 2\mathfrak{a}
     \end{array} \begin{array}{l}
       
     \end{array} \begin{array}{l}
       
     \end{array} \begin{array}{l}
       - 2\mathfrak{a}\\
       9\mathfrak{a}^2 - 4
     \end{array}
   \end{array}\right) . \]
Going one step further we note that,
\[ \det \left(\begin{array}{c}
     \begin{array}{l}
       \mathfrak{b}^2\\
       - 2\mathfrak{a}
     \end{array} \begin{array}{l}
       
     \end{array} \begin{array}{l}
       
     \end{array} \begin{array}{l}
       - 2\mathfrak{a}\\
       9\mathfrak{a}^2 - 4
     \end{array}
   \end{array}\right) = \frac{1}{9} ((9\mathfrak{a}^2 - 4) (9\mathfrak{b}^2 -
   4) - 16) = (3\mathfrak{a}\mathfrak{b}- 2\mathfrak{c})^2 \]
is a perfect square. Since the entries of this matrix are relatively prime, a
theorem by L. Mordell ([M1], [Ni]) on the decomposition of a binary quadratic
form into a sum of the square of two linear forms is applicable, yielding an
independent proof for the existence of a matrix $\left(\begin{array}{c}
  \begin{array}{l}
    p\\
    u
  \end{array} \begin{array}{l}
    q\\
    v
  \end{array}
\end{array}\right)$ which satisfies the identity (6.32). While this argument
shows that the system (6.31) is solvable for any Markoff triple
$(\mathfrak{a}, \mathfrak{b}, \mathfrak{c})$, it does not provide any
information regarding property (b) of Proposition 6.5, which is central to the
argument in the next section.

\

\

{\tmem{{\tmstrong{7 Proof of the Theorem}}}}

\

The proof consists of two parts. The first part will establish that, for any
Markoff number $\mathfrak{m}$, two Markoff forms (i. e. forms of type $F$ as
in (6.10)) affiliated with $\mathfrak{m}$, are either equivalent, or one of
these two forms is equivalent to the opposite of the other form. The argumants
employed in the proof of this statement are detached from the the formalism
developped so far, but indespensible for the second part to kick in. The point
of departure in the second part is the fact, implied by the first part, that
the integers $\mathfrak{u}_i$ and $\mathfrak{v}_i$ ($i \epsilon \{ 1, 2 \}$)
in Proposition 6.4 are uniquely determined by $\mathfrak{m}$. So, we can
employ Proposition 6.5 (with $\mathfrak{m}$, which from now on we assume to be
dominant, playing the part of $\mathfrak{a}$) to conclude that $\mathfrak{a}$
uniquely determines $\mathfrak{b}$, and hence the triple $(\mathfrak{a},
\mathfrak{b}, \mathfrak{c})$.

Before we turn to outlining the first part of the proof, we need to establish
the proper settings.

Recall that $\mathbb{m}=\mathfrak{m}$ if $\mathfrak{m}$ is odd, and
$\mathbb{m}= \frac{\mathfrak{m}}{2}$ if $\mathfrak{m}$ is even. Let
\[ \Delta = \left\{ \begin{array}{c}
     9\mathbb{m}^2 - 4 \begin{array}{l}
       
     \end{array} \tmop{if} \mathfrak{m} \tmop{is} \tmop{odd}\\
     9\mathbb{m}^2 - 1 \begin{array}{l}
       
     \end{array} \tmop{if} \mathfrak{m} \tmop{is} \tmop{even}
   \end{array} \right. \]
From now on $\mathbb{E}$ denotes the prinipal form with discriminant $\Delta$.
Note that a binary quadratic form $G (s, t) = (1 + 3\mathfrak{f}\mathfrak{g})
s^2 + 3 (\mathfrak{f}^2 -\mathfrak{g}^2) s t + (1 - 3\mathfrak{f}\mathfrak{g})
t^2$ is primitive if and only if $\mathfrak{f}^2 +\mathfrak{g}^2$ is odd,
regardless of the specific properties of the integers $\mathfrak{f}$ and
$\mathfrak{g}$. In case $\mathfrak{f}^2 +\mathfrak{g}^2$ is even, the
coefficients of $G$ are divisible by 2, in which case the form $\frac{1}{2} G$
is primitive. The equivalence of the forms in (6.10) implies that the
coefficients of $F$ must be divisible by 2 as well, and that $\frac{1}{2} F$
is primitive. Let
\[ \mathbb{F}= \left\{ \begin{array}{c}
     F \quad \tmop{if} \mathfrak{m} \tmop{is} \tmop{odd}\\
     \frac{1}{2} \begin{array}{l}
       F
     \end{array} \tmop{if} \mathfrak{m} \tmop{is} \tmop{even}
   \end{array} \right. \]
There are two basic observations which will be instrumental in the first part
of the proof. First, the principal form $E$ of discriminant $9\mathfrak{m}^2 -
4$ represents $-\mathfrak{m}^2$ (this is exactly how Frobenius defined a
Markoff number at the beginning of his paper). Any representation of
$-\mathbb{m}^2$ by the principal form $\mathbb{E}$ of discriminant $\Delta$ is
proper. The latter property will be used, by switching to the settings of
ideal factorization in quadratic number fields ([He], [H-K], or more generally
in Dedekind domains, [Ma]), to show that the form class $[-\mathbb{E}]$ has an
essentially unique factorization into classes of forms representing the
squares of prime factors of $\mathbb{m}$, and that such a factor can never
have the order 2 in the form class group, and it can have the order 4 only in
case $\mathbb{m}$ is a prime number. The second basic observation is the fact
that the square of a form class containing the form $\mathbb{F}$ is equal to
$[-\mathbb{E}]$. Since $\mathbb{F}$ represents $\mathbb{m}$, the (essentially)
unique factorization of $[-\mathbb{E}]$ indicated above, entails the
(essential) uniqueness of a Markoff form. We turn now to the details.

\

{\tmstrong{7.1}} {\tmstrong{Lemma}} If $q$ is a proper factor of
$\mathbb{m}$, then $- q^2$ cannot be represented by $\mathbb{E}$.

\

{\tmstrong{Proof}} The identity $\mathbb{E} (x, y) = -\mathbb{m}^2$ holds for
some integers $x$ and $y$ if and only if $(x, y, \mathfrak{m})$ is a Markoff
triple. In particular, $\gcd (x, y) = 1$, i. e. all such representations of
$-\mathbb{m}^2$ are proper. Suppose that $q$ is a proper factor of
$\mathbb{m}$, such that $\mathbb{E} (x, y) = - q^2$ for some integers $x, y$.
Let $q' = \frac{\mathbb{m}}{q}$. Then,
\[ \mathbb{E} (q' x, q' y) = -\mathbb{m}^2 . \]
Since this is an improper representation of the number $-\mathbb{m}^2$, we
have reached a contradiction. $\Box$

\

Given any (proper) representation of $-\mathbb{m}^2$,
\[ \mathbb{E} (x, y) = -\mathbb{m}^2, \]
there exist integers $u, v$, such that $x u - y v$=1, and the matrix \
$\left(\begin{array}{c}
  \begin{array}{l}
    x\\
    u
  \end{array} \begin{array}{l}
    y\\
    v
  \end{array}
\end{array}\right)$ \ transforms $-\mathbb{E}$ into a quadratic form
\[ K (x, y) =\mathbb{m}^2 x^2 + n x y + \tilde{n} y ^2, \]
for certain integers $n$ and $\tilde{n}$. In the sequel we will frequently
shift back and forth between form classes and ideal classes in the narrow
sense, for which there is a one-to-one correspondence. Therefore, from now on,
only ideal classes in the narrow sense will be considered. \ In the present
case this means that we associate with the form $K$ an ideal $\mathcal{J}_K$
as follows,
\[ \mathcal{J}_K = \left\{ \mathbb{m}^2 \alpha + \frac{n - \sqrt{\Delta}}{2}
   \beta \quad / \begin{array}{l}
     
   \end{array} \alpha, \beta \tmop{algebraic} \tmop{integers} \right\} . \]
This ideal has the norm
\[ \mathcal{N} (\mathcal{J}_K) =\mathbb{m}^2 . \]
Let $\mathbb{m}= p_1 \ldots .p_r$ be a prime factorization of $\mathbb{m}$.

\

{\tmstrong{7.2}} {\tmstrong{Lemma}} If $\mathcal{J}_1$ and $\mathcal{J}_2$
are ideals belonging to the ideal class corresponding to $[-\mathbb{E}]$, such
that $\mathcal{N} (\mathcal{J}_1) =\mathcal{N} (\mathcal{J}_2) =\mathbb{m}^2$,
then either $\mathcal{J}_1 =\mathcal{J}_2$ or $\mathcal{J}_1 =
\overline{\mathcal{J}_2}$, where $\overline{\mathcal{J}_2}$ denotes the
conjugate ideal of $\mathcal{J}_2$.

\

{\tmstrong{Proof}} Suppose that there are two ideals, $\mathcal{J}_1$ and
$\mathcal{J}_2$, both belonging to the ideal class corresponding to
$[-\mathbb{E}]$, such that $\mathcal{J}_1 \neq \mathcal{J}_2$ and
$\mathcal{J}_1 \neq \overline{\mathcal{J}_2}$. For each of these two ideals \
there is a unique factorization into squares of prime ideals related to the
prime factorization of $\mathbb{m}$,
\[ \mathcal{J}_1 = \dot{\mathcal{J}_{p_{_1}}^2} . . .
   \dot{\mathcal{J}_{p_{_r}}^2}, \begin{array}{l}
     
   \end{array} \begin{array}{l}
     
   \end{array} \mathcal{J}_2 = \ddot{\mathcal{J}_{p_{_1}}^2} . . .
   \ddot{\mathcal{J}_{p_{_r}}^2}, \]
where $\dot{\mathcal{J}_{p_{_{\iota}}}^2}, \ddot{\mathcal{J}_{p_{\iota}}^2}
\epsilon \left\{ \mathcal{I}_{p_{\iota}}, \overline{\mathcal{I}_{p_{\iota}}}
\right\}$, and $\mathcal{I}_{p_{\iota}}, \overline{\mathcal{I}_{p_{\iota}}}$
are the only prime ideals containing $p_{\iota}$, for $\iota = 1 \ldots$r, \
([He], Satz 90, [H-K], Theorem 5.8.8(1)). Here it is important to note that
$\mathbb{m}$ and $\Delta$ are relatively prime. We rearrange the primes in the
factorization of $\mathbb{m}$ in such a way that
\[ \dot{\mathcal{J}_{p_{_{\iota}}} } = \ddot{\mathcal{J}_{p_{_{\iota}}} }
   \tmop{for} 1 \leqslant \iota \leqslant s, \begin{array}{l}
     
   \end{array} \dot{\mathcal{J}_{p_{_{\iota}}} } \neq
   \ddot{\mathcal{J}_{p_{_{\iota}}} } \tmop{for} s < \iota \leqslant r. \]
For any algebraic integer $\kappa$ we denote by $\langle \kappa \rangle$ the
principal ideal generated by $\kappa$. Since $\mathcal{J}_1$ and
$\mathcal{J}_2$ are equivalent ideals, there exist algebraic integers $\alpha$
and $\beta$ such that
\[ \langle \alpha \rangle \mathcal{J}_1 = \langle \beta \rangle \mathcal{J}_2
   . \]
We assume that the norm $\mathcal{N} (\langle \alpha \rangle)$ (and hence also
$\mathcal{N} (\langle \beta \rangle)$ is minnimal. Since the factorization
into prime ideals on either side of this equation is unique, it follows that
\[ \langle \alpha \rangle = \ddot{{\mathcal{J}^2_{p_{}}} }_{s + 1} . . .
   \ddot{{\mathcal{J}^2_{p_{_r}}} }, \begin{array}{l}
     
   \end{array} \begin{array}{l}
     
   \end{array} \langle \beta \rangle = \dot{{\mathcal{J}^2_{p_{}}} }_{s + 1} .
   . . \dot{{\mathcal{J}^2_{p_{_r}}} }, \]
which implies a factorization of the principal form class $[\mathbb{E}]$ into
squares of form classes representing the corresponding primes in the
fatorization of $\mathbb{m}$,

\[ [\mathbb{E}] = \left[ \ddot{F_{p_{}} }_{s + 1}  \right]^2 . . . \left[
   \ddot{F_{{p_r}_{}} }  \right]^2 = \left[ \dot{F_{p_{}} }_{s + 1} \right]^2
   . . . \left[ \dot{F_{p_{ r}}} \right]^2 . \]
Hence
\[ [-\mathbb{E}] = \left[ \dot{F_{p_{}} }_1 \right]^2 . . . \left[
   \dot{F_{p_r}} \right]^2 = \left[ \dot{F_{p_{}} }_1 \right]^2 . . . \left[
   \dot{F_{p_s}} \right]^2 . \]
Let $q = p_1 . . . p_s $. Then the last equation means that the form
$-\mathbb{E}$ represents $q^2$ ([H-K], Corollary 6.4.8(1)), or equivalently
that $\mathbb{E}$ represents $- q^2$. Since, by assumption, $\mathcal{J}_1
\neq \mathcal{J}_2$ and $\mathcal{J}_1 \neq \overline{\mathcal{J}_2}$, the
integer $q$ is has to be a proper factor of $\mathbb{m}$, thus contradicting
Lemma 7.1. $\Box$

\

For each prime number $p$ there exist exactly two form classes of
discriminant $\Delta$ representing $p$ ([H-K],Theorem 6.4.13(2)). For each
such prime we choose one of those two classes, and denote it by $[F_p]$. Then
the other class which represents $p$ is equal to $[F_p ]^{- 1}$. The ideal
class corresponding to $[F_p]$ contains exactly one prime ideal, which we
denote by $\mathcal{J}^{(1)}_p$, such that $\mathcal{N} (\mathcal{J}^{(1)}_p)
= p$ ([H-K], Theorem 5.8.8(1)). Likewise, we denote by $\mathcal{J}^{(- 1)}_p$
the unique prime ideal belonging to the ideal class corresponding to $[F_p
]^{- 1}$.

\

{\tmstrong{7.3 Corollary}} The form class $[-\mathbb{E}]$ admits exactly two
facorizations into squares of classes containing forms which represent the
prime factors of $\mathbb{m}$. One of those two factorizations can be obtained
from the other one by simply taking the inverses of each factor. Moreover,
none of these squares of classes is equal to the principal form, and none of
them is ambiguous, unless $\mathbb{m}$ is a prime number.

\

{\tmstrong{Proof}} First, such a factorization exists, because there is a
prime factorization of the ideal $\mathcal{J}_K$, which gives rise to a
factorization of the corresponding form classes. Suppose that
\[ [-\mathbb{E}] = [\nobracket F _{p_1}]^{2 \varepsilon_1} \nobracket . . .
   [\nobracket F _{p_r}]^{2 \varepsilon_r} \nobracket, \begin{array}{l}
     
   \end{array} \varepsilon_{\iota} \epsilon \{ 1, - 1 \} \tmop{for} \iota = 1,
   . . . r \]
is any factorization of this type. Then the ideal
\[ \mathcal{J}= (\nobracket \mathcal{J}^{(\varepsilon_1)}_{p_1} . .
   .\mathcal{J}^{(\varepsilon_r)}_{p_r})^2 \nobracket \]
has the norm $\mathcal{N} (\mathcal{J}) = \mathfrak{m}^2$, and it belongs to
the ideal class corresponding to $[-\mathbb{E}]$. By Lemma 7.2, $\mathcal{J}$
and $\bar{\mathcal{J}} = (\nobracket \mathcal{J}^{(- \varepsilon_1)}_{p_1} . .
.\mathcal{J}^{(- \varepsilon_r)}_{p_r})^2 \nobracket$ are the only ideals
sharing these properties. Since the prime factorizations of those two ideals
are unique, this in turn implies that $[-\mathbb{E}]$ admits only two
facorizations into squares of classes containing forms which represent the
prime factors of $\mathbb{m}$, namely
\[ [-\mathbb{E}] = [\nobracket F _{p_1}] \nobracket^{2 \varepsilon_1} . . .
   [\nobracket F _{p_r}]^{2 \varepsilon_r} \nobracket \infixand  [-\mathbb{E}]
   = [\nobracket F _{p_1}]^{- 2 \varepsilon_1} \nobracket . . . [\nobracket F
   _{p_r}]^{- 2 \varepsilon_r} \nobracket . \]
None of the classes $[F_{p_{\iota}} ]^2$ is equal to the principal class,
because the opposite statement would once again conflict with Lemma 7.2. Also,
the proof of Lemma 7.2 specialized to the case $s = r - 1$, shows that
$[F_{p_{\iota}} ]^2$ cannot be ambiguous. Equivalently, $[F_{p_{\iota}} ]$ can
never have the order 4 in the form class group. $\Box$
\[ \  \]
{\tmstrong{7.4 Proposition}} \ If $\mathbb{G}$ is a form of discrminant
$\Delta$ which represents $\mathbb{m}$ properly, such that $[\mathbb{G}]^2 =
[-\mathbb{E}]$, then there is only one other form sharing these two
properties, namely the opposite form $\bar{\mathbb{G}}$. Hence, $\{
\mathbb{G}, \bar{\mathbb{G}} \} = \{ \mathbb{F}, \bar{\mathbb{F}} \}$.

\

{\tmstrong{Proof}} Let $\mathbb{G}$ be a form with the two stated properties.
Since $\mathbb{G}$ represents $\mathbb{m}$ properly, there exists an ideal
$\mathcal{I}$ in the ideal class corresponding to $[\mathbb{G}]$, such that
$\mathcal{N} (\mathcal{I}) = \mathbb{m}$. Let
\[ \mathcal{I}=\mathcal{J}^{(\varepsilon_1)}_{p_1} . .
   .\mathcal{J}^{(\varepsilon_r)}_{p_r}, \begin{array}{l}
     
   \end{array} \varepsilon_{\iota} \epsilon \{ 1, - 1 \} \quad \tmop{for}
   \iota = 1, . . . r, \]
the prime factorization of $\mathcal{I}$, and let

(7.1)
\[ [\mathbb{G}] = [\nobracket F _{p_1}] \nobracket^{\varepsilon_1} . . .
   [\nobracket F _{p_r}] \nobracket^{\varepsilon_r} \]
be the the corresponding factorization of $[\mathbb{G}]$ into form classes
representing the primes in the prime factorization of $\mathbb{m}$. By
assumption,
\[ [-\mathbb{E}] = [\nobracket F _{p_1}]^{2 \varepsilon_1} \nobracket . . .
   [\nobracket F _{p_r}]^{2 \varepsilon_r} \nobracket . \]
Since, by the second statement in Corollary 7.3 none of the factors on the
right hand side of this identity can be an ambiguous form class, thus ensuring
that $[\nobracket F _{p_{\iota}}]^{\varepsilon_{\iota}} \nobracket$ is always
uniquely determined by $[\nobracket F _{p_{\iota}}]^{2 \varepsilon_{\iota}}
\nobracket$, it follows from the first statement in Corollary 7.3 that each
factor on the right hand side of (7.1) is uniquely determined by
$[-\mathbb{E}]$, save for a uniform change for all factors on the right hand
side of (7.1) from $\varepsilon_{\iota}$ to $- \varepsilon_{\iota}$. The
impact of effecting this change, however, means nothing but switching from
$[\mathbb{G}]$ to $[\bar{\mathbb{G}}]$. In conclusion, $[\mathbb{F}]$ and
$[\bar{\mathbb{F}}]$ are indeed the only classes sharing the two properties
enuciated in Proposition 7.3. $\Box$

\

{\tmstrong{Remarks}} 1) Note that, in the line of reasoning above, we did not
encounter any complications due to the possible occurrence of multiple prime
facors in the factorization of $\mathbb{m}$ due to cancellation of the
corresponding form classes in the form class group, because by Lemma 7.1, such
cancellations are prohibited.

2) In view of the line of reasoning leading up to the conclusion in
Proposition 7.4 it is tempting to speculate about the structure of the
subgroup $\mathcal{G}$ of the form class group generated by the form classes
$[\nobracket F _{p_{\iota}}]  \nobracket$. The fact that the smallest non-zero
absolute value of integers represented by $\mathbb{F}$ is equal to
$\mathbb{m}$, entails that the smallest non-zero absolute value of integers
represented by $[\nobracket F _{p_{\iota}}] \nobracket$ is equal to
$p_{\iota}$. It follows that the set of those generators of $\mathcal{G}$ has
the cardinality $r$. By contrast, does the set of squares of these generators
and the squares of their inverses contain only two elements ? If this were the
case, then $\mathcal{G}$ would turn out to be isomorphic to a direct product
of $r$ cyclic groups of order $4 r$. In particular, the number of odd prime
factors of the Markoff number $\mathfrak{m}$ would determine the structure of
the group $\mathcal{G}$ completely. Moreover, Corollary 7.3 would follow
trivially.

2) Let $\mathbb{m}$ be any integer which is the sum of two perfect squares.
Then there exist integers $b$ and $n$, such that $b^2 + 1 = (4 n + 1)
\mathbb{m}^2$, or equivalently, $b^2 + 4\mathbb{m}= (4 n + 5) \mathbb{m}^2 -
1$. The quadratic form $F (x, y) =\mathbb{m}x^2 + b x y -\mathbb{m}y^2$, which
is symmetric, has discriminant $(4 n + 5) \mathbb{m}^2 - 1$, and it represents
$\mathbb{m}$. The integer $4 n + 5$ is a perfect square if and only if $n =
1$, and this happens to be true if and only if $2\mathbb{m}$ is a Markoff
number. The general setting adopted here extends certain features of the
situation considered above. What's missing, though, is the property that all
representations of $-\mathbb{m}^2$ by the principal form of discriminant $(4 n
+ 5) \mathbb{m}^2 - 1$ have to be proper. In the context of Markoff numbers,
this is the situation where a Markoff form is equal to one of the two
symmetric forms to which it is equivalent, and which occurs 15 times among the
first 40 Markoff numbers. Are there infinitely many Markoff numbers with this
property ?

\

The combination of Proposition 6.3, Proposition 6.4 and Proposition 7.4 leads
to the following conclusion.

\

{\tmstrong{7.5 Corollary}} For every Markoff number $\mathfrak{m} \geqslant
5$, there exist exactly four integers $u_1, v_1, u_2, v_2$ with the following
properties:

(7.2) The integers $u_1$ and $u_2$ are positive and even.

(7.3)
\[ u^2_1 + v^2_1 = u^2_2 + v^2_2 = 9\mathfrak{m}^2 - 4 \]
(7.4) If $f$ and $g$ are integers such that $f^2 + g^2 =\mathfrak{m}$, and the
form
\[ G (s, t) = (1 + 3 f g) s^2 + 3 (f^2 - g^2) s t + (1 - 3 f g) t^2 \]
is equivalent to a symmetric form $Q$, then $Q \epsilon \{ \pm H_1, \pm H_2,
\pm \bar{H}_1, \pm \bar{H}_2 \}$, where
\[ H_i (s, t) = \frac{u_i}{2} s^2 + v_i s t - \frac{u_i}{2} t^2,
   \begin{array}{l}
     
   \end{array} \text{} \bar{H}_i (s, t) = \frac{u_i}{2} s^2 - v_i s t -
   \frac{u_i}{2} t^2, \begin{array}{l}
     
   \end{array} i \epsilon \{ 1, 2 \} . \]

Equipped with Proposition 6.5 and Corollary 7.5 we are now in a position to
prove the following statement.

\

{\tmstrong{Proposition 7.6}} If $(\mathfrak{a}, \mathfrak{m}, \mathfrak{c})$
is a Markoff triple such that $\mathfrak{m} \geqslant 5$ and $\mathfrak{a},
\mathfrak{c} \leqslant \mathfrak{m}$, then there exist integers
$\tmmathbf{u}_a, \tmmathbf{v}_a, \tmmathbf{u}_c, \tmmathbf{v}_c$ such that

(7.5)
\[ \{ \tmmathbf{u}_a, \tmmathbf{v}_a \} = \{ u_i, v_i \}, \begin{array}{l}
     
   \end{array} \{ \tmmathbf{u}_c, \tmmathbf{v}_c \} = \{ u_j, v_j \},
   \begin{array}{l}
     
   \end{array} \tmop{where} \begin{array}{l}
     
   \end{array} \{ i, j \} = \{ 1, 2 \} . \]
($u_1, v_1, u_2, v_2$ as in Corollary 7.5)

(7.6)
\[ \mathfrak{a}^2 =\tmmathbf{p}^2_a +\tmmathbf{q}^2_a, \begin{array}{l}
     
   \end{array} \mathfrak{c}^2 =\tmmathbf{p}^2_c +\tmmathbf{q}^2_c, \]
where $\tmmathbf{p}_a, \tmmathbf{q}_a, \tmmathbf{p}_c, \tmmathbf{q}_c$ are
integers whose absolute values are uniquely determined by the properties

(7.7)
\[ | \tmmathbf{p}_a \tmmathbf{u}_a +\tmmathbf{q}_a \tmmathbf{v}_a | =
   2\mathfrak{m}, \begin{array}{l}
     
   \end{array} \tmop{and} \begin{array}{l}
     {}[\nobracket
   \end{array} | \tmmathbf{p}_a | < \frac{| \tmmathbf{v}_a |}{2},
   \begin{array}{l}
     
   \end{array} \tmop{if} \begin{array}{l}
     
   \end{array} | \tmmathbf{u}_a | \leqslant | \tmmathbf{v}_a |
   \begin{array}{l}
     ] \nobracket
   \end{array} \begin{array}{l}
     
   \end{array} \tmop{or} \begin{array}{l}
     {}[\nobracket
   \end{array} | \tmmathbf{q}_a | < \frac{| \tmmathbf{u}_a |}{2},
   \begin{array}{l}
     
   \end{array} \tmop{if} \begin{array}{l}
     
   \end{array} | \tmmathbf{v}_a | \leqslant | \tmmathbf{u}_a |
   \begin{array}{l}
     ] \nobracket
   \end{array} \]
(7.8)
\[ | \tmmathbf{p}_c \tmmathbf{u}_c +\tmmathbf{q}_c \tmmathbf{v}_c | =
   2\mathfrak{m}, \begin{array}{l}
     
   \end{array} \tmop{and} \begin{array}{l}
     {}[\nobracket
   \end{array} | \tmmathbf{p}_c | < \frac{| \tmmathbf{v}_c |}{2},
   \begin{array}{l}
     
   \end{array} \tmop{if} \begin{array}{l}
     
   \end{array} | \tmmathbf{u}_c | \leqslant | \tmmathbf{v}_c |
   \begin{array}{l}
     ] \nobracket
   \end{array} \begin{array}{l}
     
   \end{array} \tmop{or} \begin{array}{l}
     {}[\nobracket
   \end{array} | \tmmathbf{q}_c | < \frac{| \tmmathbf{u}_c |}{2},
   \begin{array}{l}
     
   \end{array} \tmop{if} \begin{array}{l}
     
   \end{array} | \tmmathbf{v}_c | \leqslant | \tmmathbf{u}_c |
   \begin{array}{l}
     ] \nobracket
   \end{array} \]

{\tmstrong{Proof}} The essence of the argument is the following. A linear
diophantine equation
\[ A x + B y = C \]
has at most one solution $(x_0, y_0)$ such that $| x_0 | < \frac{| B |}{2}$.
We are going to apply this to the linear diophantine equation which we obtain
from the third identity in (6.31) by letting $A = u, B = v,$ and $C = -
2\mathfrak{a}$ for appropriately chosen Markoff triples. Since in each of the
cases we consider the existence of such a solution guaranteed, we can infer
the uniqueness of the Markoff numbers $\mathfrak{a}$ and $\mathfrak{c}$ from
the uniqueness of these solutions. Specifically, in order to deal with the
case addressed in (7.7) we shall employ Proposition 6.5, using the following
correspondence of parameters, in the order listed below.
\[ \begin{array}{l}
     \tmop{Proposition} 6.5\\
     \begin{array}{l}
       
     \end{array} \begin{array}{l}
       
     \end{array} \begin{array}{l}
       
     \end{array} \begin{array}{l}
       
     \end{array} \mathfrak{a}\\
     \begin{array}{l}
       
     \end{array} \begin{array}{l}
       
     \end{array} \begin{array}{l}
       
     \end{array} \begin{array}{l}
       
     \end{array} \mathfrak{b}\\
     \begin{array}{l}
       
     \end{array} \begin{array}{l}
       
     \end{array} \begin{array}{l}
       
     \end{array} \begin{array}{l}
       
     \end{array} \mathfrak{c}\\
     \begin{array}{l}
       
     \end{array} \begin{array}{l}
       
     \end{array} \begin{array}{l}
       
     \end{array} \begin{array}{l}
       
     \end{array} \mathfrak{m}
   \end{array} \begin{array}{l}
     
   \end{array} \begin{array}{l}
     
   \end{array} \begin{array}{l}
     
   \end{array} \begin{array}{l}
     
   \end{array} \begin{array}{l}
     \tmop{Proposition} 7.6\\
     \begin{array}{l}
       
     \end{array} \begin{array}{l}
       
     \end{array} \begin{array}{l}
       
     \end{array} \begin{array}{l}
       
     \end{array} \mathfrak{m}\\
     \begin{array}{l}
       
     \end{array} \begin{array}{l}
       
     \end{array} \begin{array}{l}
       
     \end{array} \begin{array}{l}
       
     \end{array} \mathfrak{a}\\
     \begin{array}{l}
       
     \end{array} \begin{array}{l}
       
     \end{array} \begin{array}{l}
       
     \end{array} \begin{array}{l}
       
     \end{array} \mathfrak{c}\\
     \begin{array}{l}
       
     \end{array} \begin{array}{l}
       
     \end{array} 3\mathfrak{m}\mathfrak{c}-\mathfrak{a}
   \end{array} \]
It is important to note that, since $\mathfrak{m}$ in the right column is
dominant (among the first three entries only), $\mathfrak{m}$ in the left
column has to be dominant, thus rendering Proposition 6.5 applicable to our
present situation. By Proposition 6.5 there exists $i \epsilon \{ 1, 2 \}$,
and there exist integers $\tmmathbf{u}_a, \tmmathbf{v}_a$, $\tmmathbf{p}_a,
\tmmathbf{q}_a,$ such that

(7.9)
\[ \{ \tmmathbf{u}_a, \tmmathbf{v}_a \} = \{ u_i, v_i \}, \begin{array}{l}
     
   \end{array} \mathfrak{a}^2 =\tmmathbf{p}^2_a +\tmmathbf{q}^2_a, \text{}
   \begin{array}{l}
     
   \end{array} | \tmmathbf{p}_a \tmmathbf{u}_a +\tmmathbf{q}_a \tmmathbf{v}_a
   | = 2\mathfrak{m}. \]
If $| \tmmathbf{u}_a | \leqslant | \tmmathbf{v}_a |$ then, since
$\tmmathbf{u}^2_a +\tmmathbf{v}^2_a = 9\mathfrak{m}^2 - 4$, and since
$\mathfrak{m} \geqslant 5$ by assumption,
\[ | \tmmathbf{v}_a | > 2\mathfrak{m}. \]
Hence, by the second property of (7.9),
\[ | \tmmathbf{p}_a | \leqslant \mathfrak{a}<\mathfrak{m}< \frac{|
   \tmmathbf{v}_a |}{2}, \]
which, in light of our introductory remark, means that the number $|
\tmmathbf{p}_a |$ is uniquely determined by the integers $\tmmathbf{u}_a$ and
$\tmmathbf{v}_a$ via the the third property of (7.9). This in turn implies
that the integer $| \tmmathbf{q}_a |$ is uniquely determined by the third
property in (7.9) as well. Finally, the second property in (7.9) implies that
the Markoff number $\mathfrak{a}$ is uniquely determined be the integers
$\tmmathbf{u}_a$ and $\tmmathbf{v}_a$. If $| \tmmathbf{v}_a | < |
\tmmathbf{u}_a |$, then the same line of reasoning, \ simply by replacing
$\tmmathbf{v}_a$ by $\tmmathbf{u}_a$, leads to the same conclusion.

Next, to deal with the case addressed in (7.8) we employ Proposition 6.5 once
again, using the following correspondence of parameters, in the order listed
below.
\[ \begin{array}{l}
     \tmop{Proposition} 6.5\\
     \begin{array}{l}
       
     \end{array} \begin{array}{l}
       
     \end{array} \begin{array}{l}
       
     \end{array} \begin{array}{l}
       
     \end{array} \mathfrak{a}\\
     \begin{array}{l}
       
     \end{array} \begin{array}{l}
       
     \end{array} \begin{array}{l}
       
     \end{array} \begin{array}{l}
       
     \end{array} \mathfrak{b}\\
     \begin{array}{l}
       
     \end{array} \begin{array}{l}
       
     \end{array} \begin{array}{l}
       
     \end{array} \begin{array}{l}
       
     \end{array} \mathfrak{c}\\
     \begin{array}{l}
       
     \end{array} \begin{array}{l}
       
     \end{array} \begin{array}{l}
       
     \end{array} \begin{array}{l}
       
     \end{array} \mathfrak{m}
   \end{array} \begin{array}{l}
     
   \end{array} \begin{array}{l}
     
   \end{array} \begin{array}{l}
     
   \end{array} \begin{array}{l}
     
   \end{array} \begin{array}{l}
     \tmop{Proposition} 7.6\\
     \begin{array}{l}
       
     \end{array} \begin{array}{l}
       
     \end{array} \begin{array}{l}
       
     \end{array} \begin{array}{l}
       
     \end{array} \mathfrak{m}\\
     \begin{array}{l}
       
     \end{array} \begin{array}{l}
       
     \end{array} \begin{array}{l}
       
     \end{array} \begin{array}{l}
       
     \end{array} \mathfrak{c}\\
     \begin{array}{l}
       
     \end{array} \begin{array}{l}
       
     \end{array} \begin{array}{l}
       
     \end{array} \begin{array}{l}
       
     \end{array} \mathfrak{a}\\
     \begin{array}{l}
       
     \end{array} \begin{array}{l}
       
     \end{array} 3\mathfrak{m}\mathfrak{a}-\mathfrak{c}
   \end{array} \]
Since $\mathfrak{m}$ in the right column is the dominant (among the first
three entries only), $\mathfrak{m}$ in the left column has to be dominant. By
Proposition 6.5 there exists $j \epsilon \{ 1, 2 \}$, and there exist integers
$\tmmathbf{u}_c, \tmmathbf{v}_c$, $\tmmathbf{p}_c, \tmmathbf{q}_c,$ such that

(7.10)
\[ \{ \tmmathbf{u}_c, \tmmathbf{v}_c \} = \{ u_j, v_j \}, \begin{array}{l}
     
   \end{array} \mathfrak{c}^2 =\tmmathbf{p}^2_c +\tmmathbf{q}^2_c, \text{}
   \begin{array}{l}
     
   \end{array} | \tmmathbf{p}_c \tmmathbf{u}_c +\tmmathbf{q}_c \tmmathbf{v}_c
   | = 2\mathfrak{m}. \]
If $| \tmmathbf{u}_c | \leqslant | \tmmathbf{v}_c |$ then, since
$\tmmathbf{u}^2_c +\tmmathbf{v}^2_c = 9\mathfrak{m}^2 - 4$, and since
$\mathfrak{m} \geqslant 5$ by assumption,
\[ | \tmmathbf{v}_c | > 2\mathfrak{m}. \]
Hence, by the second property of (7.10),
\[ | \tmmathbf{p}_c | \leqslant \mathfrak{c}<\mathfrak{m}< \frac{|
   \tmmathbf{v}_c |}{2}, \]
which, means that the number $| \tmmathbf{p}_c |$ is uniquely determined by
the integers $\tmmathbf{u}_c$ and $\tmmathbf{v}_c$ via the third property of
(7.10). This in turn implies that the integer $| \tmmathbf{q}_c |$ is uniquely
determined by the third property in (7.10) as well. Finally, the second
property in (7.9) implies that the Markoff number $\mathfrak{c}$ is uniquely
determined by the integers $\tmmathbf{u}_c$ and $\tmmathbf{v}_c$. If $|
\tmmathbf{v}_c | < | \tmmathbf{u}_c |$, then the same line of reasoning, \
simply by replacing $\tmmathbf{v}_c$ by $\tmmathbf{u}_c$, leads to the same
conclusion. $\Box$

\

\

{\tmstrong{Remark }}The reason why we have been so painstakingly repetitive
in dealing with the very similar cases (7.7) and (7.8) is the following. Since
$\mathfrak{m} \geqslant 5$ by assumption, the Markoff numbers $\mathfrak{a}$
and $\mathfrak{c}$ have to be distinct. Thus, disregarding Proposition 6.4 in
the formulation of Corollary 7.5, that is disregarding the fact that we
already know that there are two ``non-affiliated'' symmetric forms which are
equivalent to a given Markoff form (or its opposite form), and departing
instead from the weaker assumption that there exists at least one such form,
the arguments above provide us with an independent proof that there exist
actually two ``non-affiliated'' symmetric forms which are equivalent to a
given Markoff form (or its opposite form). To summarize, we have arrived at
the following conclusion. Given a Markoff triple $(\mathfrak{a}, \mathfrak{m},
\mathfrak{c})$, such that $\mathfrak{m} \geqslant \max \{ \mathfrak{a},
\mathfrak{c}, 5 \}$, the unique pair of equivalence classes of binary
quadratic forms of discrimiant $9\mathfrak{m}^2 - 4$, one class being the
opposite of the other, which contain at least one symmetric form each, and for
which each form belonging to either one of these two classes represents the
integer $\mathfrak{m}$, contain each exactly two distinct symmetric forms,
$H_1 (s, t) =\mathfrak{u}_1 s^2 +\mathfrak{v}_1 s t -\mathfrak{u}_1 t^2$,
$\mathfrak{u}_1 > 0$, $H_2 (s, t) =\mathfrak{u}_2 s^2 +\mathfrak{v}_2 s t
-\mathfrak{u}_2 t^2$, $\mathfrak{u}_2 > 0$, and $\bar{H}_1 (s, t)
=\mathfrak{u}_1 s^2 -\mathfrak{v}_1 s t -\mathfrak{u}_1 t^2$, $\bar{H}_2 (s,
t) =\mathfrak{u}_2 s^2 -\mathfrak{v}_2 s t -\mathfrak{u}_2 t^2$, respectively.
One of these two forms determines the Markoff number $\mathfrak{a}$, while the
other one determines the Markoff number $\mathfrak{c}$. The modus of this
relationship between the set $\{ \mathfrak{a}, \mathfrak{c} \}$ and the set
$\{ H_1, H_2 \}$, or alternatively the set $\{ \bar{H}_1, \bar{H}_2 \}$, \ has
been described in the proof of Proposition 7.6.

\

From the previous remark it is clear that Proposition 7.6 implies the Theorem
for $\mathfrak{m} \geq 5$. Since the claim is trivially true for $\mathfrak{m}
\epsilon \{ 1, 2 \}$, the Theorem has been proved.

\

\

{\tmstrong{8 Principal class and Gauss composition}}

\

In the present section we revisit the settings of Proposition 6.3 and,
following the outline of Olga Taussky's paper [T] characterizing integral $2
\times 2$ matrices which can be factored into two symmetric integral matrices,
achieve a more specific description of the forms and their composition which
were instrumental in the proof of Proposition 7.1. Proposition 6.3
characterized Markoff numbers as the trace of integer matrices $\mathfrak{G}=
\left(\begin{array}{c}
  \begin{array}{l}
    \begin{array}{l}
      
    \end{array} 3\mathfrak{g}^2\\
    3\mathfrak{f}\mathfrak{g}+ 1
  \end{array} \begin{array}{l}
    3\mathfrak{f}\mathfrak{g}- 1\\
    \begin{array}{l}
      
    \end{array} 3\mathfrak{f}^2
  \end{array}
\end{array}\right)$ which are equivalent to symmetric matrices, or
equivalently, for which there exists a (symmetric) integral unimodular \
matrix $S$ such that

(8.1)
\[ \mathfrak{G}S = S\mathfrak{G}^t . \]
The matrix $S = \left(\begin{array}{c}
  \begin{array}{l}
    x\\
    y
  \end{array} \begin{array}{l}
    y\\
    z
  \end{array}
\end{array}\right)$ solves the matrix identity (8.1) if and only if

(8.2)
\[ (1 + 3\mathfrak{f}\mathfrak{g}) x + 3 (\mathfrak{f}^2 -\mathfrak{g}^2) y +
   (1 - 3\mathfrak{f}\mathfrak{g}) z = 0. \]
Since $S$ is assumed to be unimodular,

(8.3)
\[ y^2 + 1 = x z. \]
We are now going to derive a binary quadratic form $Q$ of discriminant $4
(9\mathfrak{m}^2 - 4)$ such that the combination of (8.1) and (8.3) is
equivalent to the statement that $Q$ represents the number 4. In order to
facilitate the notation in the subsequent manipulations we employ (6.9) to
recast the linear diophantine equation (8.2),

(8.4)
\[ (1 + 3\mathfrak{q}) x + 3\mathfrak{p}y + (1 - 3\mathfrak{q}) z = 0. \]
Our first task is to obtain a characterization of the bases of the two
dimensional lattice of solutions to the equation (8.2) which is suitable for
our present purpose. The following criterion was stated without proof in [Sm].

\

{\tmstrong{8.1 Lemma}} Two integral vectors in the two dimensional lattice of
solutions of a homogeneous linear diophantine equation in three independent
variables form a basis if and only if the greatest common divisor of the
components of their cross product is equal to one.

\

{\tmstrong{Proof}} Consider the general linear diophantine equation in three
variables,

(8.5)
\[ \tmmathbf{a}x +\tmmathbf{b}y +\tmmathbf{c}z = 0, \]
as well as two vectors
\[ \varphi_i = \left(\begin{array}{c}
     u_i\\
     v_i\\
     w_i
   \end{array}\right) \epsilon \mathbb{Z}^3 \nocomma, \begin{array}{l}
     
   \end{array} i \epsilon \{ 1, 2 \} . \]
For any $2 \times 2$ matrix $R = \left(\begin{array}{c}
  \begin{array}{l}
    r_{11}\\
    r_{21}
  \end{array} \begin{array}{l}
    r_{12}\\
    r_{22}
  \end{array}
\end{array}\right)$ we have

(8.6)
\[ (r_{11} \varphi_1 + r_{12} \varphi_2) \times (r_{21} \varphi_1 + r_{22}
   \varphi_2) = \det (R) (\varphi_1 \times \varphi_2) . \]
By [G], Lemma 279 there exist vectors $\psi_1$, $\psi_2 \epsilon \mathbb{Z}^3$
such that

(8.7)
\[ \psi_1 \times \psi_2 = \gcd (\tmmathbf{a}, \tmmathbf{b}, \tmmathbf{c})^{-
   1} \left(\begin{array}{c}
     \tmmathbf{a}\\
     \tmmathbf{b}\\
     \tmmathbf{c}
   \end{array}\right) . \]
Now suppose that $\varphi_1$ and $\varphi_2$ form a basis of the two
dimensional lattice $\mathcal{L}$ of integral solutions of (8.5). Then there
exists a matrix $R \epsilon M_2 (\mathbb{Z})$ such that
\[ \psi_1 = r_{11} \varphi_1 + r_{12} \varphi_2, \begin{array}{l}
     
   \end{array} \psi_2 = r_{21} \varphi_1 + r_{22} \varphi_2 . \]
Hence by (8.6) and (8.7),
\[ \det (R) (\varphi_1 \times \varphi_2) = \gcd (\tmmathbf{a}, \tmmathbf{b},
   \tmmathbf{c})^{- 1} \left(\begin{array}{c}
     \tmmathbf{a}\\
     \tmmathbf{b}\\
     \tmmathbf{c}
   \end{array}\right) . \]
Since the vectors $\varphi_1 \times \varphi_2$ and \ $\left(\begin{array}{c}
  \tmmathbf{a}\\
  \tmmathbf{b}\\
  \tmmathbf{c}
\end{array}\right)$ are collinear, it follows that $| \det (R) | = 1$, which
means that the greatest common divisor of the components of $\varphi_1 \times
\varphi_2$ is equal to one.

Now suppose that $\varphi_1$ and $\varphi_2$ are linearly independent vectors
in $\mathcal{L}$ such that the greatest common divisor of the components of
$\varphi_1 \times \varphi_2$ is equal to one. Let \ $\{ \psi_1, \psi_2 \}$ be
a basis of $\mathcal{L}$. Then there exists a matrix $R \epsilon M_2
(\mathbb{Z})$ such that
\[ \varphi_1 = r_{11} \psi_1 + r_{12} \psi_2, \begin{array}{l}
     
   \end{array} \varphi_2 = r_{21} \psi_1 + r_{22} \psi_2 . \]
By (8.6),
\[ \varphi_1 \times \varphi_2 = \det (R) (\psi_1 \times \psi_2) . \]
Since by the first part of the proof the greatest common divisor of the
components of $\psi_1 \times \psi_2$ is equal to one, while the same is true,
by assumption, for the vector $\varphi_1 \times \varphi_2$, we must have $|
\det (R) | = 1$, which means that $\{ \varphi_1, \varphi_2 \}$ is a basis as
well. \ $\Box$

\

The following vectors belong to the lattice $\mathcal{L}$ of solutions of the
equation (8.4),

(8.8)
\[ \varphi_1 = \left(\begin{array}{c}
     3\mathfrak{p}\\
     - 2\\
     3\mathfrak{p}
   \end{array}\right), \begin{array}{l}
     
   \end{array} \varphi_2 = \left(\begin{array}{c}
     1 - 3\mathfrak{q}\\
     0\\
     - (1 + 3\mathfrak{q})
   \end{array}\right) . \]
Moreover,

(8.8)
\[ \varphi_1 \times \varphi_2 = 2 \left(\begin{array}{c}
     1 + 3\mathfrak{q}\\
     3\mathfrak{p}\\
     1 - 3\mathfrak{q}
   \end{array}\right) . \]
Now suppose that $\mathfrak{m}$ is an even Markoff number. Then $\mathfrak{p}$
is even, while $\mathfrak{q}$ is odd. Since the greatest common divisor of the
components of the vector $\frac{1}{4} (\varphi_1 \times \varphi_2)$ is equal
to one, it follows from Lemma 7.1 that the vectors

(8.9)
\[ \psi_1 = \frac{1}{2} \varphi_1, \begin{array}{l}
     
   \end{array} \psi_2 = \frac{1}{2} \varphi_2 \]
from a basis of $\mathcal{L}$. By Proposition 6.3 this implies that there
exist integers $s$ and $t$ such that the components of the vector
\[ \left(\begin{array}{c}
     x\\
     y\\
     z
   \end{array}\right) = s \psi_1 + t \psi_2 \]
solve the equation (8.3). Hence,
\[ s^2 + 1 = \frac{1}{4} (3\mathfrak{p}s + (1 - 3  \mathfrak{q}) t)
   (3\mathfrak{p}s - (1 + 3\mathfrak{q}) t) = \frac{1}{4} [9\mathfrak{p}^2 s^2
   + 9\mathfrak{q}^2 t^2 - t^2 - 18\mathfrak{p}\mathfrak{q}s t], \]
or, letting

(8.10)
\[ Q (s, t) = (9\mathfrak{p}^2 - 4) s^2 - 18\mathfrak{p}\mathfrak{q}s t +
   (9\mathfrak{q}^2 - 1) t^2 = 9 (\mathfrak{p}s -\mathfrak{q}t)^2 - (4 s^2 +
   t^2), \]
(8.11)
\[ Q (s, t) = 4. \]
For the discriminant $\Delta$ of the quadratic form $Q$ we obtain,

(8.12)
\[ \Delta = 18^2 \mathfrak{p}^2 \mathfrak{q}^2 - 4 (9\mathfrak{p}^2 - 4)
   (9\mathfrak{p}^2 - 1) = 4 [9 (\mathfrak{p}^2 + (2\mathfrak{q})^2) - 4] = 4
   (9\mathfrak{m}^2 - 4) . \]
Next suppose that $\mathfrak{m}$ is an odd Markoff number. Then $\mathfrak{p}$
is odd and $\mathfrak{q}$ is even. Let

(8.13)
\[ \psi_1 = \frac{1}{2} (\varphi_1 + \varphi_2), \begin{array}{l}
     
   \end{array} \psi_2 = \varphi_2 . \]
Since $\psi_1 \times \psi_2 = \left(\begin{array}{c}
  1 + 3\mathfrak{q}\\
  3\mathfrak{p}\\
  1 - 3\mathfrak{q}
\end{array}\right)$, it follows from Lemma 8.1 that $\{ \psi_1, \psi_2 \}$ is
a basis of $\mathcal{L}$. By Proposition 6.3 this implies that there exist
integers $s$ and $t$ such that the components of the vector
\[ \left(\begin{array}{c}
     x\\
     y\\
     z
   \end{array}\right) = s \psi_1 + t \psi_2 \]
solve the equation (8.3). Hence,
\[ s^2 + 1 = \left[ (3 (\mathfrak{p}-\mathfrak{q}) + 1) \frac{s}{2} + (1 -
   3\mathfrak{q}) t \right] \left[ (3 (\mathfrak{p}-\mathfrak{q}) - 1)
   \frac{s}{2} - (1 + 3\mathfrak{q}) t \right] = \]
\[ (9 (\mathfrak{p}-\mathfrak{q})^2 - 1) \frac{s^2}{4} + (9\mathfrak{q}^2 - 1)
   t^2 + [(1 - 3\mathfrak{q}) (3 (\mathfrak{p}-\mathfrak{q}) - 1) - (1 +
   3\mathfrak{q}) (3 (\mathfrak{p}-\mathfrak{q}) + 1)] \frac{s t}{2} = \]
\[ (9 (\mathfrak{p}-\mathfrak{q})^2 - 1) \frac{s^2}{4} + (9\mathfrak{q}^2 - 1)
   t^2 - (9\mathfrak{q} (\mathfrak{p}-\mathfrak{q}) + 1) s t, \]
which is equivalent to,

(8.14)
\[ (9 (\mathfrak{p}-\mathfrak{q})^2 - 5) s^2 - 4 (9\mathfrak{q}
   (\mathfrak{p}-\mathfrak{q}) + 1) s t + 4 (9\mathfrak{q}^2 - 1) t^2 = 4. \]
We claim that $s$ can be chosen to be even. Since $\mathfrak{p}$ is odd and
$\mathfrak{q}$ is even, there exists an integer $n$, as well as odd integers
$\mathfrak{n}_1 \nocomma, \mathfrak{n}_{2,} \mathfrak{n}_3$ such that
\[ 9 (\mathfrak{p}-\mathfrak{q})^2 - 5 = 9 (2 n + 1)^2 - 5 = 4 [9 n (n + 1) +
   1] = 4\mathfrak{n}_1, \]
\[ 4 (9\mathfrak{q} (\mathfrak{p}-\mathfrak{q}) + 1) = 4 [9\mathfrak{q} (2 n +
   1) + 1] = 4\mathfrak{n}_2, \]
\[ 4 (9\mathfrak{q}^2 - 1) = 4\mathfrak{n}_3 . \]
Thus, dividing both sides in (8.14) by 4 leads to

(8.15)
\[ \mathfrak{n}_1 s^2 -\mathfrak{n}_2 s t +\mathfrak{n}_3 t^2 = 1. \]
Since the discriminant of the quadratic form on the left hand side of (8.15)
is equal to $9\mathfrak{m}^2 - 4$, a fundamental isomorph has the form
\[ \left(\begin{array}{c}
     \begin{array}{l}
       \frac{3\mathfrak{m}+\mathfrak{n}_2}{2}\\
       \\
       \begin{array}{l}
         
       \end{array} \begin{array}{l}
         
       \end{array} \mathfrak{n}_1
     \end{array} \begin{array}{l}
       \begin{array}{l}
         
       \end{array} -\mathfrak{n}_3\\
       \\
       \begin{array}{l}
         
       \end{array} \frac{3\mathfrak{m}-\mathfrak{n}_2}{2}
     \end{array}
   \end{array}\right) \begin{array}{l}
     
   \end{array} \tmop{or} \begin{array}{l}
     
   \end{array} \left(\begin{array}{c}
     \begin{array}{l}
       \frac{3\mathfrak{m}-\mathfrak{n}_2}{2}\\
       \\
       \begin{array}{l}
         
       \end{array} -\mathfrak{n}_1
     \end{array} \begin{array}{l}
       \begin{array}{l}
         
       \end{array} \begin{array}{l}
         
       \end{array} \mathfrak{n}_3\\
       \\
       \begin{array}{l}
         
       \end{array} \frac{3\mathfrak{m}+\mathfrak{n}_2}{2}
     \end{array}
   \end{array}\right) \begin{array}{l}
     
   \end{array} . \]
Since $\mathfrak{m}- 1 = 0 \begin{array}{l}
  (\tmop{mod} 4)
\end{array}$, this implies that in terms of parity a fundamental isomorph
displays only one of the following two \ patterns
\[ \left(\begin{array}{c}
     \begin{array}{l}
       \tmop{even}\\
       \tmop{odd}
     \end{array} \begin{array}{l}
       \tmop{odd}\\
       \tmop{odd}
     \end{array}
   \end{array}\right) \begin{array}{l}
     
   \end{array} \tmop{or} \begin{array}{l}
     
   \end{array} \left(\begin{array}{c}
     \begin{array}{l}
       \tmop{odd}\\
       \tmop{odd}
     \end{array} \begin{array}{l}
       \tmop{odd}\\
       \tmop{even}
     \end{array}
   \end{array}\right) . \]
One of these patterns can be obtained from the other by considering the
inverse of the respective isomorph. If $s$ is odd and $t$ is even, then we
apply the fundamental automorph with the parity pattern \
$\left(\begin{array}{c}
  \begin{array}{l}
    \tmop{even}\\
    \tmop{odd}
  \end{array} \begin{array}{l}
    \tmop{odd}\\
    \tmop{odd}
  \end{array}
\end{array}\right) \begin{array}{l}
  
\end{array}$to the vector $\left(\begin{array}{c}
  s\\
  t
\end{array}\right)$, yielding an even number in the first entry of the
resulting vector. If $s$ is odd and $t$ is odd, then we apply the fundamental
automorph with the parity pattern $\left(\begin{array}{c}
  \begin{array}{l}
    \tmop{odd}\\
    \tmop{odd}
  \end{array} \begin{array}{l}
    \tmop{odd}\\
    \tmop{even}
  \end{array}
\end{array}\right)$ to the vector $\left(\begin{array}{c}
  s\\
  t
\end{array}\right)$, yielding again an even number in the first entry of the
resulting vector. In conclusion, we may indeed assume that the number $s$ in
(8.14) is even. Hence, letting $s = 2 r$, and writing once again $s$ in place
of $r$, we obtain,
\[ 9 [\mathfrak{p}s -\mathfrak{q} (s + t)]^2 - [4 s^2 + (s + t)^2] = 1, \]
and finally, letting \ $r = s + t$, but writing $t$ in place of $r$,

(8.16)
\[ Q (s, t) = 1, \]
where $Q$ is defined as in (8.10). The following statement summarizes what has
been accomplished so far.

\

{\tmstrong{8.2 Proposition}} If $Q (s, t) = (9\mathfrak{p}^2 - 4) s^2 -
18\mathfrak{p}\mathfrak{q}s t + (9\mathfrak{q}^2 - 1) t^2$, then $Q$ belongs
to the principal class of discriminant $4 (9\mathfrak{m}^2 - 4)$ in case
$\mathfrak{m}$ is odd, while $\frac{1}{4} Q$ belongs to the principal class of
discriminant $\frac{9\mathfrak{m}^2 - 4}{4} = 9\mathbb{m}^2 - 1$ \ in case
$\mathfrak{m}$ is even.

\

{\tmstrong{Proof}} In case $\mathfrak{m}$ is odd the claim is an immediate
consequence of (8.16). In case $\mathfrak{m}$ is even we note that, since
$\mathfrak{p}$ is even, while $\mathfrak{q}$ is odd, the coefficients of $Q$
are divisible by 4. Hence the claim follows from (8.11) in this
case.$\begin{array}{l}
  
\end{array} \Box$

\

Noting that all the steps leading up to Proposition 8.2 are reversible, we
obtain the following characterization of Markoff numbers.

\

{\tmstrong{8.3 Corollary}} A positive integer $\mathfrak{m}$ is a Markoff
number if and only if there exist integers $\mathfrak{p}$ and $\mathfrak{q}$
such that $\mathfrak{p}^2 + 4\mathfrak{q}^2 =\mathfrak{m}^2$, and the
quadratic form $Q (s, t) = (9\mathfrak{p}^2 - 4) s^2 -
18\mathfrak{p}\mathfrak{q}s t + (9\mathfrak{q}^2 - 1) t^2$ belongs to the
principal class of discriminant $4 (9\mathfrak{m}^2 - 4)$ in case
$\mathfrak{m}$ is odd, while the quadratic form $\frac{1}{4} Q$ belongs to the
principal class of discriminant $9\mathbb{m}^2 - 1$ in case $\mathfrak{m}$ is
even.

\

For arbitrary integers $\mathfrak{p}$ and $\mathfrak{q}$ let
\[ G (x, y) = (1 + 3\mathfrak{q}) x^2 + 3\mathfrak{p}x y + (1 - 3\mathfrak{q})
   y^2, \begin{array}{l}
     
   \end{array} G^{\sharp} = (2 + 3\mathfrak{p}) x^2 + 12\mathfrak{q}x y + (2 -
   3\mathfrak{p}) y^2 . \]
Then we have the following identities,
\[ Q (u v, v^2 - u^2) = - G (u, v) G (v, - u), \begin{array}{l}
     
   \end{array} Q (v^2 - u^2, 4 u v) = - G^{\sharp} (u, v) G^{\sharp} (v, - u)
   . \]
Both of these two identities are particular manifestations of Gauss
compositions of binary quadratic forms. More specifically, we obtain the
following (see for instance [Sp], pp. 382-383)

\

8.4 {\tmstrong{Proposition}} (a) If $\mathfrak{m}$ is odd, and we define the
bilinear substitution
\[ x_1 = (2 + 3\mathfrak{p}) y_1 z_2 + (2 - 3\mathfrak{p}) y_2 z_1, \]
\[ x_2 = - y_1 z_1 - 3\mathfrak{q}y_1 z_2 + 3\mathfrak{q}y_2 z_1 + y_2 z_2, \]
then,
\[ Q (x_2, x_1) = - G^{\sharp} (y_1, y_2) G^{\sharp} (z_2, - z_1) \]
(b) If $\mathfrak{m}$ is even, and we define the bilinear substitution
\[ x_1 = \frac{1 + 3\mathfrak{q}}{2} y_1 z_2 + \frac{1 - 3\mathfrak{q}}{2} y_2
   z_1, \]
\[ x_2 = - y_1 z_1 - \frac{3\mathfrak{p}}{2} y_1 z_2 + \frac{3\mathfrak{p}}{2}
   y_2 z_1 + y_2 z_2, \]
then,
\[ Q (x_1, x_2) = - G (y_1, y_2) G (z_2, - z_1) . \]

Denoting Gauss composition of forms by $\ast$, Proposition 8.4 implies,
\[ G^{\sharp} \ast G^{\sharp}  \begin{array}{l}
     
   \end{array} \tmop{is} \tmop{equivalent} \tmop{to} \begin{array}{l}
     
   \end{array} - \tilde{Q}, \begin{array}{l}
     
   \end{array} \tmop{where} \tilde{Q} (s, t) = Q (t, s), \]
in case $\mathfrak{m}$ is odd, and
\[ \left( \frac{1}{2} G \right) \ast \left( \frac{1}{2} G \right)
   \begin{array}{l}
     
   \end{array} \tmop{is} \tmop{equivalent} \tmop{to} \begin{array}{l}
     
   \end{array} - \frac{1}{4} Q \]
in case $\mathfrak{m}$ is even. The appearance of a seemingly different form
in the composition in case $\mathfrak{m}$ \ is odd can be explained as
follows. In the following discussion we adopt the notation introduced in the
proof of Proposition 7.1. By [H-K], Theorem 6.4.14 there exists a canonical
group homomorphism
\[ \Theta_{\Delta, 2} : \begin{array}{l}
     
   \end{array} \mathfrak{F}_{4 \Delta} \rightarrow \mathfrak{F}_{\Delta} . \]
If $\Delta$ is a positive discriminant such that $\Delta = 1 \begin{array}{l}
  
\end{array} \text{}$(mod 4), which is the case when $\mathfrak{m}$ is odd,
then the kernel of this homomorphism is trivial, that is $\Theta_{\Delta, 2}$
is actually an isomorphism, unless the following two conditions are met,

(8.17)
\[ \Delta = 5 \begin{array}{l}
     
   \end{array} (\tmop{mod} 8) \begin{array}{l}
     
   \end{array} \tmop{and} \begin{array}{l}
     
   \end{array} \varepsilon_{\Delta} \epsilon \mathcal{O}_{4 \Delta}, \]
where $\varepsilon_{\Delta}$ is the fundamental unit of discriminant $\Delta$,
and $\mathcal{O}_{4 \Delta}$ is the maximal order of the quadratic number
field $\mathbb{Q} \left( \sqrt{4 \Delta} \right)$. If $\mathfrak{m}$ is an odd
Markoff number, then the first of those two conditions is met, while the
second condition fails: $\varepsilon_{\Delta} = \frac{3\mathfrak{m}+
\sqrt{9\mathfrak{m}^2 - 4}}{2} \nin \mathcal{O}_{4 \Delta}$. Hence,
$\Theta_{\Delta, 2}$ is an isomorphism. Since

(8.18)
\[ G^{\sharp} (\mathfrak{f}, \mathfrak{g}) =\mathfrak{m} (3\mathfrak{m}- 2),
\]
and since either $\mathfrak{f}$ or $\mathfrak{g}$ is even, each form in the
form class $\Theta_{\Delta, 2} ([G^{\sharp}])$ represents $\mathfrak{m}
(3\mathfrak{m}- 2)$. Since $\mathfrak{m}$ and $3\mathfrak{m}- 2$ are
relatively prime, it follows from [H-K], Corollary 6.4.9 that
\[ \Theta_{\Delta, 2} ([G^{\sharp}]) =\mathcal{G}\mathcal{K}, \]
where each form in $\mathcal{G}$ represents $\mathfrak{m}$, and each form in
$\mathcal{K}$ represents $3\mathfrak{m}- 2$. Since $3\mathfrak{m}- 2$ divides
the discriminant $9\mathfrak{m}^2 - 4$, the class $\mathcal{K}$ is ambiguous,
i.e. $\mathcal{K}^2 =\mathcal{E}$. If $\mathfrak{m}$ is an odd Markoff number,
then by Corollary 8.3 the form $Q$ is in the principal class of discriminant
$4 (9\mathfrak{m}^2 - 4)$, and therefore $Q$ is equivalent to $\tilde{Q}$.
Hence,
\[ \Theta_{\Delta, 2} ([G^{\sharp}])^2 = \Theta_{\Delta, 2} ([G^{\sharp}]^2) =
   \Theta_{\Delta, 2} ([- \tilde{Q}]) = \Theta_{\Delta, 2} ([- Q]) =
   \bar{\mathcal{E}} . \]
In combination with Proposition 7.1 this implies that $\mathcal{G}$ is either
equal to $[G]$ or the opposite class $[\bar{G}]$ of $[G]$.

\

\

{\tmstrong{{\tmem{9. Markoff triples and the norm form equation}}}}

\begin{tabular}{l}
  
\end{tabular}

Having established the uniqueness of a dominant Markoff number in Section 7,
there are two aspects that will be touched upon in the remainder of this work.
First, a description of the data that are being determined by a single Markoff
number $\mathfrak{m}$ in a way that reflects its dominance, and second, in
consideration of the multitude of identities that led to the conclusion of the
uniqueness of a dominant Markoff number, to highlight the purely algebraic
side of the formalism. To deal with the former, we shall adopt as a framework
a norm from equation that uses no data other than $\mathfrak{m}$ and the
discriminant $9\mathfrak{m}^2 - 4$. To appreciate the need for the latter, it
suffices to point out, that for any pair of non-zero rational numbers $u
\tmop{and} v$ the triple $(a, b, c)$ of rational numbers, where
\[ a = \frac{u^2 + v^2 + 1}{u}, \begin{array}{l}
     
   \end{array} b = \frac{u^2 + v^2 + 1}{u}, \begin{array}{l}
     
   \end{array} c = \frac{u^2 + v^2 + 1}{u v}, \]
solves the Markoff equation $a^2 + b^2 + c^2 = a b c$, a fact that hints at a
lack of depth of the whole formalism when considered within this broader
setting. The two aspects turn out to be linked to each other in some way. To
begin with, we need to introduce the necessary framework for the discussion
and switch to a more expedient notation. Let $(\mathfrak{m}, \mathfrak{a}_0,
\mathfrak{a}_1)$ be a Markoff triple such that $\mathfrak{m} \geqslant
\mathfrak{a}_1 \geqslant \mathfrak{a}_0$, and define recursively

(9.1)
\[ \mathfrak{a}_{n + 1} = 3\mathfrak{m}\mathfrak{a}_n -\mathfrak{a}_{n - 1}
   \tmop{for} n \geqslant 1, \begin{array}{l}
     
   \end{array} \mathfrak{a}_{n - 1} = 3\mathfrak{m}\mathfrak{a}_n
   -\mathfrak{a}_{n + 1} \tmop{for} n \leqslant 0 \]
Then the uniqueness of the dominant Markoff number $\mathfrak{m}$ implies
that, up to permutations, the two-sided sequence of triples $(\mathfrak{m},
\mathfrak{a}_n, \mathfrak{a}_{n + 1}), n \epsilon \mathbb{Z}$, represents
exactly all those Markoff triples which contain $\mathfrak{m}$ as a member.
Notice, however, that in case $\mathfrak{m}= 1$ \ or $\mathfrak{m}= 2$ the
recursion is essentially only one-sided, leading to a duplication of Markoff
numbers if the recursion is two-sided. Let $\lambda$ be the following
fundamental unit and its inverse, respectively, in the quadratic field
$\mathbb{Q} (\sqrt{9\mathfrak{m}^2 - 4})$,

(9.2)
\[ \lambda = \frac{3\mathfrak{m}}{2} + \frac{\sqrt{9\mathfrak{m}^2 - 4}}{2},
   \begin{array}{l}
     
   \end{array} \lambda^{- 1} = \frac{3\mathfrak{m}}{2} -
   \frac{\sqrt{9\mathfrak{m}^2 - 4}}{2} \]
For $x = r + s \sqrt{9\mathfrak{m}^2 - 4}  (r, s \epsilon \mathbb{Q})$ in
$\mathbb{Q} (\sqrt{9\mathfrak{m}^2 - 4})$ we denote by $x^{\ast} = r - s
\sqrt{9\mathfrak{m}^2 - 4}$ its conjugate. Let

(9.3)
\[ \omega = \frac{\mathfrak{a}_1 -\mathfrak{a}_0 \lambda^{-
   1}}{\sqrt{9\mathfrak{m}^2 - 4}} = \frac{\mathfrak{a}_0}{2} +
   \frac{3\mathfrak{a}_0 \mathfrak{m}- 2\mathfrak{a}_1}{2
   \sqrt{9\mathfrak{m}^2 - 4}} . \]
Then

(9.4)
\[ \mathfrak{a}_n = \omega \lambda^n + \omega^{\ast} \lambda^{- n}
   \begin{array}{l}
     
   \end{array} \tmop{for} \tmop{all} n \begin{array}{l}
     \epsilon
   \end{array} \mathbb{Z}, \]
and $\omega \omega^{\ast}$ solves the norm form equation,

(9.5)
\[ (9\mathfrak{m}^2 - 4) \omega \omega^{\ast} =\mathfrak{m}^2, \]
or, written as a diophantine equation,

(9.6)
\[ x^2 - D y^2 = - 4\mathfrak{m}^2, \begin{array}{l}
     
   \end{array} \tmop{where} \begin{array}{l}
     
   \end{array} D = 9\mathfrak{m}^2 - 4, \begin{array}{l}
     
   \end{array} \]
where,

(9.7)
\[ \begin{array}{l}
     
   \end{array} x = 3\mathfrak{a}_0 \mathfrak{m}- 2\mathfrak{a}_1,
   \begin{array}{l}
     
   \end{array} y =\mathfrak{a}_0 . \]
A solution \ $(x, y) = (u, v)$ of the norm form equation (9.6) is called a
fundamental solution ([St]), if the following two inequalities hold,

(9.8)
\[ 0 < v \leqslant \frac{\mathfrak{m}}{\sqrt{3\mathfrak{m}- 2}},
   \begin{array}{l}
     
   \end{array} \begin{array}{l}
     | u | \leqslant \mathfrak{m} \sqrt{3\mathfrak{m}- 2}
   \end{array} . \]
For any solution $(x, y)$ of (9.6) there exists a fundamental solution $(u,
v)$ and an integer $n$ such that
\[ x + y \sqrt{D} = (u + v \sqrt{D}) \lambda^n, \]
and two solutions for which such a relation holds with a common fundamental
solution $u + v \sqrt{D}$ are called equivalent. In the general theory of norm
form equations it is shown, that the first equation in (9.6), with more
general parameters on either side of the equation, has only finitely many
fundamental solutions. The solution (9.7) is a fundamental solution. Since
$\mathfrak{a}_1 \geqslant \mathfrak{a}_0$ by assumption, the first inequality
in (9.8) trivially implies the second one. Switching the roles of
$\mathfrak{a}_0 \tmop{and} \mathfrak{a}_1$ leads to the conjugate equivalence
class of solutions, which in this particular case is always distinct from the
former in case $\mathfrak{m} \geqslant 5$. The uniqueness of the dominant
Markoff number $\mathfrak{m}$ is equivalent to the statement that there are no
other equivalence classes of solutions.

\begin{tabular}{l}
  
\end{tabular}

\

{\tmstrong{{\tmem{10 Recursions for the discriminant}}}}

\begin{tabular}{l}
  
\end{tabular}

We return now to the settings of Proposition 6.5, while retaining the notation
in (9.1) for Markoff numbers which belong to a triple that includes
$\mathfrak{m}$, to show that the three diophantine equations in (6.31) fit the
scheme of a recursion akin to the one in (9.1). By Proposition 6.5 there exist
four integers $u_0, u_1, v_0, v_1$ such that

(10.1 )
\[ u^2_n + v^2_n = 9\mathfrak{a}^2_n - 4 \]
(10.2)
\[ \mathfrak{p}u_n +\mathfrak{q}v_n = - 2\mathfrak{a}_n, \]
for $\left. n \begin{array}{l}
  \epsilon
\end{array} \{ 0, 1 \right\}$, and

(10.3)
\[ \mathfrak{p}^2 +\mathfrak{q}^2 =\mathfrak{m}^2 . \]
The second component identity in (6.29) yields,

(10.4)
\[ \mathfrak{a}_1 v_0 -\mathfrak{a}_0 v_1 = 2\mathfrak{p}. \]
Moreover, the Pythagorean triple $(\mathfrak{m}, \mathfrak{p}, \mathfrak{q})$
is (uniquely) affiliated with the residue classes
\[ \pm \frac{\mathfrak{a}_0}{\mathfrak{a}_1}  \begin{array}{l}
     
   \end{array} (\tmop{mod} \mathfrak{m}) . \]
By Proposition 6.5 and the second component identity in (6.29) there exist
integers $\mathbb{u}_1, \mathbb{u}_2, \mathbb{v}_1, \mathbb{v}_2$ such that

(10.5)
\[ \left. \mathbb{u}^2_i +\mathbb{v}^2_i = 9\mathfrak{a}^2_i - 4,
   \begin{array}{l}
     
   \end{array} \mathfrak{p}\mathbb{u}_i +\mathfrak{q}\mathbb{v}_i = -
   2\mathfrak{a}_i, i \epsilon \{ 1, 2 \right\} ; \begin{array}{l}
     
   \end{array} \mathfrak{a}_2 \mathbb{v}_1 -\mathfrak{a}_1 \mathbb{v}_2 =
   2\mathfrak{p}. \]
\begin{tabular}{l}
  
\end{tabular}

{\tmstrong{10.1 Lemma}} The following identities hold true,

(10.6)
\[ \mathbb{u}_2 = 3\mathfrak{m}u_1 - u_0, \begin{array}{l}
     
   \end{array} \mathbb{v}_2 = 3\mathfrak{m}v_1 - v_0 . \]
\begin{tabular}{l}
  
\end{tabular}

{\tmstrong{Proof}} \ Let $\mathfrak{u}_2 = 3\mathfrak{m}u_1 - u_0$,
$\mathfrak{v}_2 = 3\mathfrak{m}v_1 - v_0$. Then
\[ \mathfrak{p}\mathfrak{u}_2 +\mathfrak{q}\mathfrak{v}_2 =\mathfrak{p}
   (3\mathfrak{m}u_1 - u_0) +\mathfrak{q} (3\mathfrak{m}v_1 - v_0) =
   3\mathfrak{m} (\mathfrak{p}u_1 +\mathfrak{q}v_1) - (\mathfrak{p}u_0
   +\mathfrak{q}v_0) = - 2 (3\mathfrak{m}\mathfrak{a}_1 -\mathfrak{a}_0) = -
   2\mathfrak{a}_2 . \]
It follows from this and the second identity in (10.5) for $i = 2$ that there
exists an integer $x$ such that $\mathbb{v}_2 =\mathfrak{v}_2 +\mathfrak{p}x$,
and it follows from (10.2) for $n = 1$ as well as the second identity in
(10.5) for $i = 1$, that there exists an integer $y$ such that $\mathbb{v}_1 =
v_1 +\mathfrak{p}y$. The third identity in (10.5), and (10.4), together with
(9.1) yield,
\[ \begin{array}{l}
     
   \end{array} \mathfrak{a}_2 \mathbb{v}_1 -\mathfrak{a}_1 \mathbb{v}_2
   =\mathfrak{a}_2 (v_1 +\mathfrak{p}y) -\mathfrak{a}_1 (\mathfrak{v}_2
   +\mathfrak{p}x) =\mathfrak{a}_2 (v_1 +\mathfrak{p}y) -\mathfrak{a}_1
   (3\mathfrak{m}v_1 - v_0 +\mathfrak{p}x) \]
\[ = (\mathfrak{a}_2 - 3\mathfrak{m}\mathfrak{a}_1) v_1 +\mathfrak{a}_1 v_0 +
   (\mathfrak{a}_2 y -\mathfrak{a}_1 x) \mathfrak{p}= -\mathfrak{a}_0 v_1
   +\mathfrak{a}_1 v_0 + (\mathfrak{a}_2 y -\mathfrak{a}_1 x) \mathfrak{p}=
   2\mathfrak{p}+ (\mathfrak{a}_2 y -\mathfrak{a}_1 x) \mathfrak{p}=
   2\mathfrak{p}. \]
It follows that $\mathfrak{a}_2 y -\mathfrak{a}_1 x = 0$. Hence, since
$\mathfrak{a}_1 \tmop{and} \mathfrak{a}_2$ are relatively prime,
$\mathfrak{a}_2 \tmop{divides} x$, and $\mathfrak{a}_1 \tmop{divides} y$. Now
suppose that $x, \tmop{and} \tmop{hence} y$ are non-zero. Then, by (10.1) for
$i = 1$, and by the first identity in (10.5),
\[ \mathfrak{p}\mathfrak{a}_1 \leqslant | \mathbb{v}_1 - v_1 | \leqslant \max
   (| \mathbb{v}_1 |, | v_1 |) \leqslant \sqrt{9\mathfrak{a}^2_1 - 4} <
   3\mathfrak{a}_1 . \]
This implies that $\mathfrak{p} \leqslant 2$,and hence either $\mathfrak{p}=
1$ or $\mathfrak{p}= 2$. If $\mathfrak{p}= 1$, then by (10.3), $\mathfrak{m}=
1 \tmop{and} \mathfrak{q}= 0$, which is impossible. If $\mathfrak{p}= 2$, then
again by (10.3), either $\mathfrak{m}= 2 \tmop{and} \mathfrak{q}= 0$, or
$\mathfrak{m}$ is not an integer, which is also impossible. In conclusion $x =
y = 0$, thus settling the claim. $\Box$

\begin{tabular}{l}
  
\end{tabular}

Replacing $\left. i \begin{array}{l}
  \epsilon
\end{array} \{ 1, 2 \right\}$ in (10.5) by $\left. i \begin{array}{l}
  \epsilon
\end{array} \{ - 1, 0 \right\}$ and repeating the arguments in the proof of
Lemma 10.1 yields the identities $\mathbb{u}_{- 1} = 3\mathfrak{m}u_0 - u_1,
\begin{array}{l}
  
\end{array} \mathbb{v}_{- 1} = 3\mathfrak{m}v_0 - v_1 $. \ A simple induction
argument that uses nothing but Lemma 10.1, and this modified version,
establish the following.

\begin{tabular}{l}
  
\end{tabular}

{\tmstrong{10.2 Proposition}} There exist two uniquely determined two-sided
sequences of integers $\{ u_n \}$ and $\{ v_n \}$ such that
\[ u_{n + 1} = 3\mathfrak{m}u_n - u_{n - 1}, \begin{array}{l}
     
   \end{array} v_{n + 1} = 3\mathfrak{m}v_n - v_{n - 1} \]
\[ u^2_n + v^2_n = 9\mathfrak{a}^2_n - 4, \begin{array}{l}
     
   \end{array} \mathfrak{p}u_n +\mathfrak{q}v_n = - 2\mathfrak{a}_n,
   \begin{array}{l}
     
   \end{array} \]
\begin{tabular}{l}
  
\end{tabular}

{\tmstrong{Remark}} If
\[ \mathcal{A}_n = \left(\begin{array}{c}
     \begin{array}{l}
       \begin{array}{l}
         
       \end{array} 9\mathfrak{a}_{n + 1}^2 - 4\\
       9\mathfrak{a}_n \mathfrak{a}_{n + 1} - 6\mathfrak{m} \begin{array}{l}
         
       \end{array}
     \end{array} \begin{array}{l}
       9\mathfrak{a}_n \mathfrak{a}_{n + 1} - 6\mathfrak{m}\\
       \begin{array}{l}
         
       \end{array} \begin{array}{l}
         
       \end{array} 9\mathfrak{a}_n^2 - 4
     \end{array}
   \end{array}\right), \begin{array}{l}
     n \epsilon \mathbb{Z};
   \end{array} \begin{array}{l}
     
   \end{array} \mathcal{B}= \left(\begin{array}{c}
     \begin{array}{l}
       3\mathfrak{m}\\
       - 1
     \end{array} \begin{array}{l}
       1\\
       0
     \end{array}
   \end{array}\right), \]
then
\[ \det (\mathcal{A}_n) = 16, \begin{array}{l}
     
   \end{array} \mathcal{B}^t \mathcal{A}_n \mathcal{B}= \left(\begin{array}{c}
     \begin{array}{l}
       \begin{array}{l}
         
       \end{array} 9\mathfrak{a}_{n + 2}^2 - 4\\
       9\mathfrak{a}_{n + 1} \mathfrak{a}_{n + 2} - 6\mathfrak{m}
     \end{array} \begin{array}{l}
       \begin{array}{l}
         
       \end{array} 9\mathfrak{a}_{n + 1} \mathfrak{a}_{n + 2} -
       6\mathfrak{m}\\
       \begin{array}{l}
         
       \end{array} \begin{array}{l}
         
       \end{array} \begin{array}{l}
         
       \end{array} 9\mathfrak{a}_{n + 1}^2 - 4
     \end{array}
   \end{array}\right) =\mathcal{A}_{n + 1} . \]
Since $\det (\mathcal{A}_n)$ is a perfect square, and since the greatest
common divisor of the entries of $\mathcal{A}_n$ is equal to 1, it follows
from Mordell's theorem that there exists an integral 2x2 matrix
$\mathcal{C}_n$ such that
\[ \det (\mathcal{C}_n) = 4, \begin{array}{l}
     
   \end{array} \mathcal{C}_n^t \mathcal{C}_n =\mathcal{A}_n . \]
Comparison with Proposition 10.2 yields the following identity

(10.7)
\[ u_n v_n + u_{n + 1} v_{n + 1} = 9\mathfrak{a}_n \mathfrak{a}_{n + 1} -
   6\mathfrak{m}, \]
which is of some interest in its own right.

\

\begin{tabular}{l}
  
\end{tabular}

{\tmstrong{{\tmem{11 Recursions for the quadratic residues}}}}

\begin{tabular}{l}
  
\end{tabular}

In this final section the algebraic framework for the recursions involving the
parameters $k \tmop{and} l$ will be described. The major purpose is to
highlight the role of the matrix $\mathfrak{F}$ in this context. Returning to
the settings at the beginning of Section 6, especially (6.4) and (6.5), we
will employ Proposition 9.2 to convert the recursions for the $u_n \tmop{and}
v_n$ into recursions for the quadratic residues. In order to remain consistent
with the notation introduced in Section 6, the Markoff number around which the
recursion is to be developed will be denoted by $\mathfrak{b}$ rather than
$\mathfrak{m}$. In the applications of the formalism of Section 9 the letter
$\mathfrak{m}$ has to be replaced throughout by the letter $\mathfrak{b}$.
First, transcribing (9.1), the point of departure are the recursions for all
Markoff triples which include $\mathfrak{b}$,

(11.1)
\[ \mathfrak{a}_{n + 1} = 3\mathfrak{b}\mathfrak{a}_n -\mathfrak{a}_{n - 1}
   \tmop{for} n \geqslant 1, \begin{array}{l}
     
   \end{array} \mathfrak{a}_{n - 1} = 3\mathfrak{b}\mathfrak{a}_n
   -\mathfrak{a}_{n + 1} \tmop{for} n \leqslant 0 . \]
under the proviso that,

(11.2)
\[ \mathfrak{b} \geqslant \max (\mathfrak{a}_{- 1}, \mathfrak{a}_0) . \]
This implies that for some $\omega \epsilon \mathbb{Q} (\sqrt{9\mathfrak{b}^2
- 4})$ we have $\mathfrak{a}_n = \omega \lambda^n + \omega^{\ast} \lambda^{-
n} \begin{array}{l}
  
\end{array} \tmop{for} \tmop{all} n \epsilon \mathbb{Z}.$ Next, we need to
adapt the identity (5.26) to our current needs. Since the parameter $\nu$ is
going to change to its opposite sign as we pass through the triple which
$\mathfrak{b}$ dominates in the recursion (11.1), we define a one-sided
recursion, Opting for $\nu = - 1$, we note that the case for $\nu = 1$ can be
handled in a similar way. Thus the identity (5.26) for values of $n \geqslant
1$ takes the form

(11.3)
\[ \left(\begin{array}{ccc}
     {\mathfrak{a}_{n - 1}}  & \mathfrak{m}_n & \mathfrak{a}_n\\
     k_{n - 1} & \mathbb{k}_n & k_n\\
     l_{n - 1} & \mathbb{l}_n & l_n
   \end{array}\right) \begin{array}{l}
     - 1\\
     \\
     
   \end{array} = \frac{1}{2} \left(\begin{array}{ccc}
     l_{n - 1} + 3 k_{n - 1} & - (2 k_{n - 1} + 3\mathfrak{c}_{n - 1}) &
     \mathfrak{a}_{n - 1}\\
     - (l_b - 3 k_b) & 2 k_b - 3\mathfrak{b} & -\mathfrak{b}\\
     l_n - 3 k_n & - (2 k_n - 3\mathfrak{a}_n) & \mathfrak{a}_n
   \end{array}\right), \]
From Section 6, in particular (6.28), we know that for every $n$ there exists
$A_n \begin{array}{l}
  \epsilon
\end{array} \tmop{SL} (2, \mathbb{Z})$ such that
\[ \Psi (A_n) \left(\begin{array}{ccc}
     {\mathfrak{a}_{n - 1}}  & \mathfrak{m}_n & \mathfrak{a}_n\\
     k_{n - 1} & \mathbb{k}_n & k_n\\
     l_{n - 1} & \mathbb{l}_n & l_n
   \end{array}\right) = \left(\begin{array}{ccc}
     \frac{1}{2} (3\mathfrak{a}_n + v_n) & \begin{array}{l}
       
     \end{array} 1 \begin{array}{l}
       
     \end{array} & \frac{1}{2} (3\mathfrak{a}_{n - 1} + v_{n - 1})\\
     u_n & \begin{array}{l}
       
     \end{array} 0 \begin{array}{l}
       
     \end{array} & u_{n - 1}\\
     \frac{1}{2} {(3\mathfrak{a}_n}  - v_n) & 1 & \frac{1}{2}
     {(3\mathfrak{a}_{n - 1}}  - v_{n - 1})
   \end{array}\right) \]
\[ = \left(\begin{array}{ccc}
     \frac{1}{2} (3\mathfrak{a}_1 + v_1) & \begin{array}{l}
       
     \end{array} 1 \begin{array}{l}
       
     \end{array} & \frac{1}{2} (3\mathfrak{a}_0 + v_0)\\
     u_1 & \begin{array}{l}
       
     \end{array} 0 \begin{array}{l}
       
     \end{array} & u_0\\
     \frac{1}{2} {(3\mathfrak{a}_1}  - v_1) & 1 & \frac{1}{2}
     {(3\mathfrak{a}_0}  - v_0)
   \end{array}\right)  \left(\begin{array}{ccc}
     \mathfrak{b} & 0 & 1\\
     0 & 1 & 0\\
     - 1 & 0 & 0
   \end{array}\right) \begin{array}{l}
     n - 1\\
     \\
     
   \end{array} \]
\[ = \Psi (A_1) \left(\begin{array}{ccc}
     {\mathfrak{a}_0}  & \mathfrak{m}_1 & \mathfrak{a}_1\\
     k_0 & \mathbb{k}_1 & k_1\\
     l_0 & \mathbb{l}_1 & l_1
   \end{array}\right) \left(\begin{array}{ccc}
     3\mathfrak{b} & 0 & 1\\
     0 & 1 & 0\\
     - 1 & 0 & 0
   \end{array}\right) \begin{array}{l}
     n - 1\\
     \\
     
   \end{array} \]
Let $\mathcal{F}_n = \Psi (A_1)^{- 1} \Psi (A_n) $, and let $\mathcal{B}=
\left(\begin{array}{ccc}
  3\mathfrak{b} & 0 & 1\\
  0 & 1 & 0\\
  - 1 & 0 & 0
\end{array}\right)$. Then for every $n \geqslant 1$,
\[ \mathcal{F}_n \left(\begin{array}{ccc}
     {\mathfrak{a}_{n - 1}}  & \mathfrak{m}_n & \mathfrak{a}_n\\
     k_{n - 1} & \mathbb{k}_n & k_n\\
     l_{n - 1} & \mathbb{l}_n & l_n
   \end{array}\right) = \left(\begin{array}{ccc}
     {\mathfrak{a}_0}  & \mathfrak{m}_1 & \mathfrak{a}_1\\
     k_0 & \mathbb{k}_1 & k_1\\
     l_0 & \mathbb{l}_1 & l_1
   \end{array}\right) \mathcal{B}^{n - 1}, \]
and (11.3) implies
\[ \mathcal{F}_n^t \left(\begin{array}{c}
     \begin{array}{l}
       l_b - 3 k_b
     \end{array}\\
     3\mathfrak{b}- 2 k_b\\
     \mathfrak{b}
   \end{array}\right) = \left(\begin{array}{c}
     \begin{array}{l}
       l_b - 3 k_b
     \end{array}\\
     3\mathfrak{b}- 2 k_b\\
     \mathfrak{b}
   \end{array}\right) . \]
It follows that there exists an integer $j_n$ such that
\[ \mathcal{F}_n  = \left(\begin{array}{ccc}
     1 & 0 & 0\\
     0 & \frac{1}{2} & 0\\
     1 & 0 & 1
   \end{array}\right) \Psi (\mathfrak{F})^{j_n} \left(\begin{array}{ccc}
     1 & 0 & 0\\
     0 & 2 & 0\\
     0 & 0 & 1
   \end{array}\right), \begin{array}{l}
     
   \end{array} \mathfrak{F}= \left(\begin{array}{c}
     \begin{array}{l}
       3\mathfrak{m}- k\\
       \begin{array}{l}
         
       \end{array} \begin{array}{l}
         
       \end{array} \mathfrak{m}
     \end{array} \begin{array}{l}
       3 k - l\\
       \begin{array}{l}
         
       \end{array} \begin{array}{l}
         
       \end{array} k
     \end{array}
   \end{array}\right) \]
Some further considerations show that $j_n = - n + 1$, and hence a closer look
at the structure of the following matrix is desirable,

(11.4)
\[ \mathcal{F}= \left(\begin{array}{ccc}
     1 & 0 & 0\\
     0 & \frac{1}{2} & 0\\
     1 & 0 & 1
   \end{array}\right) \Psi (\mathfrak{F}) \left(\begin{array}{ccc}
     1 & 0 & 0\\
     0 & 2 & 0\\
     1 & 0 & 1
   \end{array}\right) \]
Let

(11.5)
\[ \varrho_{\pm} = \frac{2 k_b \pm 3\mathfrak{b}}{2\mathfrak{b}} +
   \frac{\sqrt{9\mathfrak{b}^2 - 4}}{2\mathfrak{b}} \]
The coefficients of the corresponding quadratic form $F$ can be recovered from
this quantity by noting,

(11.6)
\[ \frac{\mathfrak{b}}{2} (\varrho_{\pm} + \varrho_{\pm} \ast) = 2 k_b \pm
   3\mathfrak{b}, \begin{array}{l}
     
   \end{array} \mathfrak{b} \varrho_{\pm} \varrho_{\pm}^{\ast} = l_b \pm 3 k_b
\]
The diagonalization of the matrix $\mathcal{F}$ is the content of the next
statement.

\begin{tabular}{l}
  
\end{tabular}

{\tmstrong{11.1 Lemma}} The following identity holds true, letting $\varrho =
\varrho_-$,

(11.7)
\[ \frac{1}{(\varrho - \varrho \ast)^2} \left(\begin{array}{ccc}
     (\varrho^{\ast})^2 & - 2 \varrho^{\ast} & 1\\
     - 2 \varrho \varrho^{\ast} & 2 (\varrho + \varrho \ast) & - 2\\
     \varrho^2 & - 2 \varrho & 1
   \end{array}\right) \mathcal{F} \left(\begin{array}{ccc}
     1 & 1 & 1\\
     \varrho & \frac{1}{2} (\varrho + \varrho \ast) & \varrho^{\ast}\\
     \varrho^2 & \varrho \varrho^{\ast} & (\varrho^{\ast})^2
   \end{array}\right) \]
\[ \frac{1}{\omega \omega^{\ast}} \left(\begin{array}{ccc}
     (\varrho^{\ast})^2 & - 2 \varrho^{\ast} & 1\\
     - 2 \varrho \varrho^{\ast} & 2 (\varrho + \varrho \ast) & - 2\\
     \varrho^2 & - 2 \varrho & 1
   \end{array}\right) \mathcal{F} \left(\begin{array}{ccc}
     1 & 1 & 1\\
     \varrho & \frac{1}{2} (\varrho + \varrho \ast) & \varrho^{\ast}\\
     \varrho^2 & \varrho \varrho^{\ast} & (\varrho^{\ast})^2
   \end{array}\right) = \left(\begin{array}{ccc}
     \lambda^2 & 0 & 1\\
     0 & 1 & 0\\
     1 & 0 & \lambda^{- 2}
   \end{array}\right) . \]

\begin{tabular}{l}
  
\end{tabular}

The proof of Lemma 11.1 is obtained through manipulations involving the
identities (11.6). At this point a comment about the general pattern of the
eigenvalues of a matrix of the form $\Psi (\left(\begin{array}{c}
  \begin{array}{l}
    p\\
    r
  \end{array} \begin{array}{l}
    q\\
    s
  \end{array}
\end{array}\right))$ is in order. The characteristic polynomial is always of
the form
\[ - x^3 + (\mathfrak{t}^2 - 1) x^2 - (\mathfrak{t}^2 - 1) x + 1 = - (x^2 -
   (\mathfrak{t}^2 - 2) x + 1) (x - 1), \begin{array}{l}
     
   \end{array} \mathfrak{t}= \tmop{tr} (\left(\begin{array}{c}
     \begin{array}{l}
       p\\
       r
     \end{array} \begin{array}{l}
       q\\
       s
     \end{array}
   \end{array}\right)), \]
which leads to the eigenvalues
\[ 1, \begin{array}{l}
     
   \end{array} \frac{\mathfrak{t}^2 - 2 \pm \mathfrak{t} \sqrt{\mathfrak{t}^2
   - 4}}{2} = (\frac{\mathfrak{t} \pm \sqrt{\mathfrak{t}^2 - 4}}{2})^2 . \]
In particular, the eigenvalues are squares of numbers in the quadratic number
field affiliated with the discriminant. This has the interesting consequence
that the matrix $\Psi (\mathfrak{F})$ has a square root in the ring of
matrices with entries from that number field. The identity (11.3) can now be
recast as follows,
\[ \left(\begin{array}{c}
     \begin{array}{l}
       \omega \lambda^n + \omega^{\ast} \lambda^{- n}\\
       \omega \varrho \lambda^n + \omega^{\ast} \varrho^{\ast} \lambda^{- n}\\
       \omega \varrho^2 \lambda^n + \omega^{\ast} (\varrho^{\ast})^2
       \lambda^{- n} \begin{array}{l}
         
       \end{array}
     \end{array} \begin{array}{l}
       3 \omega^2 \lambda^{2 n + 1} + 3 (\omega^{\ast})^2 \lambda^{- (2 n +
       1)} + \tau\\
       3 \omega^2 \varrho \lambda^{2 n + 1} + 3 (\omega^{\ast})^2
       \varrho^{\ast} \lambda^{- (2 n + 1)} + \frac{1}{2} (\varrho +
       \varrho^{\ast}) \tau \begin{array}{l}
         
       \end{array}\\
       3 \omega^2 \varrho^2 \lambda^{2 n + 1} + 3 (\omega^{\ast})^2
       (\varrho^{\ast})^2 \lambda^{- (2 n + 1)} + \varrho \varrho^{\ast} \tau
     \end{array} \begin{array}{l}
       \omega \lambda^{n + 1} + \omega^{\ast} \lambda^{- (n + 1)}\\
       \omega \varrho \lambda^{n + 1} + \omega^{\ast} \varrho^{\ast}
       \lambda^{- (n + 1)}\\
       \omega \varrho^2 \lambda^{n + 1} + \omega^{\ast} (\varrho^{\ast})^2
       \lambda^{- (n + 1)}
     \end{array}
   \end{array}\right) \begin{array}{l}
     - 1\\
     \\
     
   \end{array} \]
\[ = \frac{1}{2} \left(\begin{array}{c}
     \begin{array}{l}
       \omega \varrho^2 \lambda^n + \omega^{\ast} (\varrho^{\ast})^2
       \lambda^{- n}\\
       \begin{array}{l}
         
       \end{array} \begin{array}{l}
         
       \end{array} - \varrho \varrho^{\ast} \mathfrak{b}\\
       \omega \varrho^2 \lambda^{n + 1} + \omega^{\ast} (\varrho^{\ast})^2
       \lambda^{- (n + 1)} \begin{array}{l}
         
       \end{array} \begin{array}{l}
         
       \end{array} \begin{array}{l}
         
       \end{array}
     \end{array} \begin{array}{l}
       - 2 (\omega \varrho \lambda^n + \omega^{\ast} \varrho^{\ast} \lambda^{-
       n})\\
       \begin{array}{l}
         
       \end{array} \begin{array}{l}
         
       \end{array} \begin{array}{l}
         
       \end{array} \begin{array}{l}
         
       \end{array} (\varrho + \varrho^{\ast}) \mathfrak{b}\\
       - (\omega \varrho \lambda^{n + 1} + \omega^{\ast} \varrho^{\ast}
       \lambda^{- (n + 1)}) \begin{array}{l}
         
       \end{array} \begin{array}{l}
         
       \end{array}
     \end{array} \begin{array}{l}
       \omega \lambda^n + \omega^{\ast} \lambda^{- n}\\
       \begin{array}{l}
         
       \end{array} \begin{array}{l}
         
       \end{array} \begin{array}{l}
         
       \end{array} -\mathfrak{b}\\
       \omega \lambda^{n + 1} + \omega^{\ast} \lambda^{- (n + 1)}
     \end{array}
   \end{array}\right) \]
\[ + \frac{3}{2} \left(\begin{array}{c}
     \begin{array}{l}
       \begin{array}{l}
         
       \end{array} \begin{array}{l}
         
       \end{array} \omega \varrho \lambda^n + \omega^{\ast} \varrho^{\ast}
       \lambda^{- n}\\
       \begin{array}{l}
         
       \end{array} \begin{array}{l}
         
       \end{array} \begin{array}{l}
         
       \end{array} \begin{array}{l}
         
       \end{array} \begin{array}{l}
         
       \end{array} \begin{array}{l}
         
       \end{array} 0\\
       - (\omega \varrho \lambda^{n + 1} + \omega^{\ast} \varrho^{\ast}
       \lambda^{- (n + 1)}) \begin{array}{l}
         
       \end{array}
     \end{array} \begin{array}{l}
       \begin{array}{l}
         
       \end{array} \begin{array}{l}
         
       \end{array} \omega \lambda^n + \omega^{\ast} \lambda^{- n}\\
       \begin{array}{l}
         
       \end{array} \begin{array}{l}
         
       \end{array} \begin{array}{l}
         
       \end{array} \begin{array}{l}
         
       \end{array} \begin{array}{l}
         
       \end{array} 0\\
       \begin{array}{l}
         
       \end{array} - (\omega \lambda^{n + 1} + \omega^{\ast} \lambda^{- (n +
       1)})
     \end{array} \begin{array}{l}
       \begin{array}{l}
         
       \end{array} \begin{array}{l}
         
       \end{array} 0\\
       \begin{array}{l}
         
       \end{array} \begin{array}{l}
         
       \end{array} 0\\
       \begin{array}{l}
         
       \end{array} \begin{array}{l}
         
       \end{array} 0
     \end{array}
   \end{array}\right), \]
where $\tau = \frac{4\mathfrak{b}}{9\mathfrak{b}^2 - 4}$, and it can be proved
algebraically (disregarding the integrality of the entries involved) by
employing the norm form equation, as well as an identity relating the
quantities $\omega \tmop{and} \varrho$,
\[ \omega \omega^{\ast} = \frac{\mathfrak{b}^2}{9\mathfrak{b}^2 - 4} =
   \frac{1}{(\varrho - \varrho^{\ast})^2} . \]
In the context of the example in the third remark at the end of Section 1, if
$(A, A B, B)$ is an admissible triple such that
\[ A B = \left(\begin{array}{c}
     \begin{array}{l}
       k_b\\
       \mathfrak{b}
     \end{array} \begin{array}{l}
       3 k_b - l_b\\
       3\mathfrak{b}- k_b
     \end{array}
   \end{array}\right) = \left(\begin{array}{c}
     \begin{array}{l}
       k_b\\
       \mathfrak{b}
     \end{array} \begin{array}{l}
       - l_b\\
       - k_b
     \end{array}
   \end{array}\right) \left(\begin{array}{c}
     \begin{array}{l}
       1\\
       0
     \end{array} \begin{array}{l}
       3\\
       1
     \end{array}
   \end{array}\right), \]
then,

\[ A = \left(\begin{array}{c}
     \begin{array}{l}
       \omega \varrho_+ \lambda + \omega^{\ast} \varrho_+^{\ast} \lambda^{-
       1}\\
       \begin{array}{l}
         
       \end{array} \omega \lambda + \omega^{\ast} \lambda^{- 1}
     \end{array} \begin{array}{l}
       3 (\omega \varrho_+ \lambda + \omega^{\ast} \varrho_+^{\ast} \lambda^{-
       1}) - (\omega \varrho^2_+ \lambda + \omega^{\ast} (\varrho_+^{\ast})^2
       \lambda^{- 1})\\
       \begin{array}{l}
         
       \end{array} \begin{array}{l}
         
       \end{array} 3 (\omega \lambda + \omega^{\ast} \lambda^{- 1}) - (\omega
       \varrho_+ \lambda + \omega^{\ast} \varrho_+^{\ast} \lambda^{- 1})
     \end{array}
   \end{array}\right) \]
\[ = \left(\begin{array}{c}
     \begin{array}{l}
       \omega \varrho_+ \lambda + \omega^{\ast} \varrho_+^{\ast} \lambda^{-
       1}\\
       \begin{array}{l}
         
       \end{array} \omega \lambda + \omega^{\ast} \lambda^{- 1}
     \end{array} \begin{array}{l}
       - (\omega \varrho^2_+ \lambda + \omega^{\ast} (\varrho_+^{\ast})^2
       \lambda^{- 1})\\
       \begin{array}{l}
         
       \end{array} - (\omega \varrho_+ \lambda + \omega^{\ast}
       \varrho_+^{\ast} \lambda^{- 1})
     \end{array}
   \end{array}\right) \left(\begin{array}{c}
     \begin{array}{l}
       1\\
       0
     \end{array} \begin{array}{l}
       3\\
       1
     \end{array}
   \end{array}\right), \]
and
\[ B = \left(\begin{array}{c}
     \begin{array}{l}
       \omega \varrho_- + \omega^{\ast} \varrho_-^{\ast}\\
       \begin{array}{l}
         
       \end{array} \omega + \omega^{\ast}
     \end{array} \begin{array}{l}
       3 (\omega \varrho_- + \omega^{\ast} \varrho_-^{\ast}) - (\omega
       \varrho^2_- + \omega^{\ast} (\varrho_-^{\ast})^2)\\
       \begin{array}{l}
         
       \end{array} \begin{array}{l}
         
       \end{array} 3 (\omega + \omega^{\ast}) - (\omega \varrho_- +
       \omega^{\ast} \varrho_-^{\ast})
     \end{array}
   \end{array}\right) \]
\[ = \left(\begin{array}{c}
     \begin{array}{l}
       \omega \varrho_- + \omega^{\ast} \varrho_-^{\ast}\\
       \begin{array}{l}
         
       \end{array} \omega + \omega^{\ast}
     \end{array} \begin{array}{l}
       - (\omega \varrho^2_- + \omega^{\ast} (\varrho_-^{\ast})^2)\\
       \begin{array}{l}
         
       \end{array} - (\omega \varrho_- + \omega^{\ast} \varrho_-^{\ast})
     \end{array}
   \end{array}\right) \left(\begin{array}{c}
     \begin{array}{l}
       1\\
       0
     \end{array} \begin{array}{l}
       3\\
       1
     \end{array}
   \end{array}\right) \]
This identity shows in a very explicit way the change of sign, here encoded in
the term $\varrho_{\pm}$, from $'' +''$ on the left, to $'' -''$ on the right,
where $\mathfrak{b}= \frac{1}{3} \tmop{tr} (A B) $is the dominant Markoff
number.

\begin{tabular}{l}
  
\end{tabular}

\begin{tabular}{l}
  
\end{tabular}

\begin{tabular}{l}
  
\end{tabular}

{\tmem{{\tmstrong{References}}}}

[AO] H. Appelgate, H. Onishi, ``The similarity problem for 3x3 integer
matrices'', Linear Algebra Appl. 42 (1982), 159-174

[Ba] P. Bachmann, ``Zahlentheorie'', B. G. Teubner, Leipzig 1898, Vierter
Teil, Erste Abteilung

[Bo] E. Bombieri, ``Continued fractions and the Markoff tree'', Expositiones
Mathematicae 25 (3) (2007), 187-213

[Bu] D. A. Buell, ``Binary quadratic forms'', Springer Verlag, 1989

[BV] J. Buchmann, U. Vollmer, ``Binary quadratic forms. An algorithmic
approach'', Springer Verlag, 2007

[Ca] J.W.S. Cassels, ``An introduction to Diophantine Approximation'',
Cambridge Univ. Press, 1957 (Chapter II)

[CV] S. Cecotti, C. Vafa, ``On the classification of N=2 Supersymmetric
Theories'', Commun. Math. Phys. 158 (1993), 569-644

[Co] H. Cohn, ``Markoff Forms and Primitive Words'', Math. Ann. 196 (1972),
8-22

[CF] T.W. Cusick, M.E. Flahive, ``The Markoff and Lagrange spectra'',
Mathematical Surveys and Monographs, 30, American Mathematical Society (1989)

[F] F.G. Frobenius ``Ueber die Markoffschen Zahlen'', Sitzungsberichte der
Koeniglichen Preus- \ \ \ sischen Akademie der Wissenschaften zu Berlin
(1913), 458-487 \ \ \ \ \ \ \ \ \ \ \ \ \ \ \ \ \ \ \ \ \ \ \ \ \ \ \ \ \ \ \
\ \ [Gesammelte Abhandlungen, Band III, Springer Verlag]

[G] C. F. Gauss, ``Disquisitiones Arithmeticae'', Yale University Press, 1966

[H-K] F. Halter-Koch, ``Quadratic Irrationals. An Introduction to Classical
Number Theory'', CRC Press, 2013

[He] E. Hecke, ``Vorlesungen ueber die Theorie der algebraischen Zahlen''
Akademische Verlagsgesellschaft Leipzig, 1923 \foreignlanguage{german}{}

[HZ] F. Hirzebruch, D. Zagier, ``The Atiyah-Singer Theorem and Elementary
Number Theory'', Publish or Perish (1974)

[Ma] D. A. Marcus, ``Number Fields'', Springer Verlag, 2018

[L] E. Landau, ``Vorlesungen ueber Zahlentheorie'', Verlag von S. Hirzel in
Leipzig (1927), Erster Band

[Mo] R. A. Mollin, ``Quadratics'', CRC Press, 1996

[M1] L. J. Mordell, ``On the representation of a binary quadratic form as a
sum of squares of linear forms'', Math. Z. 35 (1932), 1-15

[M2] L. J. Mordell, ``Diophantine Equations'', Academic Press, 1969

[Ne] M. Newman, ``Integral Matrices'', Academic Press, 1972

[Ni] I. Niven, ``Integers of quadratic fields as sum of squares'' Trans. Amer.
Math. Soc. 48 (1940), 405-417

[Pe] S. Perrine, ``L'interpretation matricielle de la theory de Markoff
classique'', Int. J. Math. Math. Sci. 32 (2002),no.4, 193-262

[Pn] O. Perron, ``Die Lehre von den Kettenbruechen'', Band I, Dritte Auflage,
B. G. Teubner (1954)

[Po] J.Popp, ``The combinatorics of frieze patterns and Markoff numbers'',
arXiv:math/0511633

[R] R. Remak ``Ueber indefinite binaere quadratische Minimalformen'', Math.
Ann. 92, 3-4 (1924), 155-182

[Re] C. Reutenauer, ``On Markoff's property and Sturmian words'', Math. Ann.
336 (1) (2006),

1-12

[Ru] A. N. Rudakov, ``The Markov numbers and exceptional bundles on $P^2$ ``,
Math. USSR. Izv., 32(1) (1989), 99-102

[Se] C. Series, ``The geometry of Markoff numbers'', Math. Intelligencer 7
(1985), no.3, 20-29

[Sm] H. J. S. Smith, ``On systems of indeterminate linear equations'', Report
of the British Association for 1860. Sectional Proceedings p.6 (Coll. Math.
Papers, I, 365-366)

[Sp] A. Speiser, ``Ueber die Komposition der binaeren quadratischen Formen''
in ``Festschrift Heinrich Weber'' (1912), 375-395

[St] B. Stolt, ``On the Diophantine equation $u^2 - D v^2 = \pm 4 N$, Part
II'', Arkiv foer Matematik 2 (2-3) (1952), 251-268

[T] O. Taussky, ``The factorization of an integral matrix into a product of
two integral symmetric matrices, II. The general case $n = 2$ ``, Commun. Pure
Appl. Math. 26 (1973), 847-854

[W] M. Waldschmidt, ``Open Diophantine Problems'', Moscow Mathematical Journal
4 (2004), no.1, 245-305

[Za] D. Zagier, ``On the Number of Markoff Numbers Below a Given Bound'',
Mathematics of Computation, 39 (1982), no.160, 709-723

[Zh1] Y. Zhang, ``An elementary proof of Markoff conjecture for prime powers''
\ \ \ \ \ \ \ \ \ \ \ \ \ \ \ \ \ \ \ \ \ \ \ \ \ arXiv:math.NT/0606283

\

\begin{tabular}{l}
  
\end{tabular}

\begin{tabular}{l}
  
\end{tabular}

$\begin{array}{l}
  \tmop{Department} \tmop{of} \tmop{Mathematics}\\
  \text{Tulane University}\\
  \text{New Orleans, LA 70118}\\
  \text{e-mail: nriedel@tulane.edu}
\end{array}$

\

\

\

\ \

\

\end{document}